\documentclass[a4paper, 10pt, final]{amsart}

\usepackage[T1]{fontenc}
\usepackage[english]{babel}
\usepackage{amsmath}
\usepackage{amsfonts}
\usepackage{amssymb}
\usepackage{graphicx}
\usepackage{amscd}
\usepackage{mathrsfs}
\usepackage{lmodern}
\usepackage{graphicx}
\usepackage{pgfplots}
\usepgfplotslibrary{groupplots}
\usepgfplotslibrary{units}

\usetikzlibrary{arrows}
\usetikzlibrary{graphs,shapes.geometric,arrows,fit,matrix,positioning}

\usepackage[OT2,T1]{fontenc}

\usepackage{amsmath, amsthm, amsfonts}

\usepackage[top=2.5cm, bottom=2.5cm, left=2.5cm, right=2.5cm]{geometry}

\makeatletter
\newcounter {subsubsubsection}[subsubsection]
\renewcommand\thesubsubsubsection{\thesubsubsection .\@alph\c@subsubsubsection}
\newcommand\subsubsubsection{\@startsection{subsubsubsection}{4}{\z@}%
                                     {-3.25ex\@plus -1ex \@minus -.2ex}%
                                     {1.5ex \@plus .2ex}%
                                     {\normalfont\normalsize\bfseries}}
\newcommand*\l@subsubsubsection{\@dottedtocline{4}{10.0em}{4.1em}}
\newcommand*{\subsubsubsectionmark}[1]{}
\makeatother

\makeatletter
\newcounter{subsubsubsubsection}[subsubsubsection]
\renewcommand\thesubsubsubsubsection{\thesubsubsubsection .\@alph\c@subsubsubsubsection}
\newcommand\subsubsubsubsection{\@startsection{subsubsubsubsection}{4}{\z@}%
 {-3.25ex\@plus -1ex \@minus -.2ex}%
 {1.5ex \@plus .2ex}%
 {\normalfont\normalsize\bfseries}}
\newcommand*\l@subsubsubsubsection{\@dottedtocline{3}{10.0em}{4.1em}}
\newcommand*{\subsubsubsubsectionmark}[1]{}
\makeatother

\usepackage{graphicx}

\DeclareMathOperator{\Reg}{Reg}
\DeclareMathOperator{\Lie}{Lie}

\DeclareMathOperator{\tder}{tder}
\DeclareMathOperator{\TAut}{TAut}

\DeclareMathOperator{\weight}{weight}
\DeclareMathOperator{\id}{id}

\DeclareMathOperator{\pr}{pr}
\DeclareMathOperator{\PGL}{PGL}
\DeclareMathOperator{\Sol}{Sol}
\DeclareMathOperator{\GRT}{GRT}
\DeclareMathOperator{\crys}{crys}

\DeclareMathOperator{\Ad}{Ad}

\DeclareMathOperator{\har}{har}

\DeclareMathOperator{\KV}{KV}
\DeclareMathOperator{\depth}{depth}
\DeclareMathOperator{\M}{M}

\DeclareMathOperator{\KRV}{KRV}

\DeclareMathOperator{\ad}{ad}
\DeclareMathOperator{\un}{un}
\DeclareMathOperator{\Stab}{Stab}

\DeclareMathOperator{\DS}{DS}

\DeclareMathOperator{\Spec}{Spec}
\DeclareMathOperator{\lie}{lie}
\DeclareMathOperator{\shft}{shft}
\DeclareMathOperator{\B}{B}
\DeclareMathOperator{\Spf}{Spf}

\DeclareMathOperator{\inv}{inv}

\DeclareMathOperator{\Gal}{Gal}
\DeclareMathOperator{\Aut}{Aut}
\DeclareMathOperator{\an}{an}

\DeclareMathOperator{\Li}{Li}

\DeclareMathOperator{\dR}{dR}

\DeclareMathOperator{\comp}{comp}

\DeclareMathOperator{\KZ}{KZ}

\DeclareMathOperator{\dec}{dec}

\DeclareMathOperator{\Moy}{Moy}

\theoremstyle{definition}

\newtheorem{Theorem}{Theorem}[subsection]
\newtheorem{Proposition}[Theorem]{Proposition}

\newtheorem{Lemma}[Theorem]{Lemma}

\newtheorem{Definition}[Theorem]{Definition}

\newtheorem{Proposition-Definition}[Theorem]{Proposition-Definition}

\newtheorem{Theorem-Definition}[Theorem]{Theorem-Definition}
\newtheorem{Lemma-Definition}[Theorem]{Lemma-Definition}

\newtheorem{Notation}[Theorem]{Notation}

\newtheorem{Corollary-Definition}[Theorem]{Corollary-Definition}
\newtheorem{Conjecture}[Theorem]{Conjecture}

\newtheorem{Question}[Theorem]{Question}
\newtheorem{Corollary}[Theorem]{Corollary}

\newtheorem{Remark}[Theorem]{Remark}
\newtheorem{Nota Bene}[Theorem]{Nota Bene}

\newcommand{\simlra}{\buildrel \sim \over \longrightarrow}

\numberwithin{equation}{subsection}

\input cyracc.def
\DeclareFontFamily{U}{russian}{}
\DeclareFontShape{U}{russian}{m}{n}
        { <5><6> wncyr5
        <7><8><9> wncyr7
        <10><10.95><12><14.4><17.28><20.74><24.88> wncyr10 }{}
\DeclareSymbolFont{Russian}{U}{russian}{m}{n}
\DeclareSymbolFontAlphabet{\mathcyr}{Russian}
\makeatletter
\let\@math@cyr\mathcyr
\renewcommand{\mathcyr}[1]{\@math@cyr{\cyracc #1}}
\makeatother
\newcommand{\sh}{\mathcyr{sh}} 

\author{David Jarossay}

\address{Ben-Gurion University of the Negev, Mathematics Department, Beer-Sheva, Israel}
\email{jarossay@post.bgu.ac.il}

\title{Adjoint cyclotomic multiple zeta values and cyclotomic multiple harmonic values}

\begin{document}

\maketitle

\begin{abstract}
We introduce adjoint cyclotomic multiple zeta values and cyclotomic multiple harmonic values. They are two variants of cyclotomic multiple zeta values, closely related to each other. They arise as key tools for the study of $p$-adic cyclotomic multiple zeta values. Moreover, cyclotomic multiple harmonic values provide an adelic lift to a cyclotomic generalization of finite multiple zeta values. We establish certain standard properties of these two objects. We consider two types of properties : some related to double shuffle relations, and some related to associator and Kashiwara-Vergne relations.

This is Part II-1 of \emph{$p$-adic cyclotomic multiple zeta values and $p$-adic pro-unipotent harmonic actions}.
\end{abstract}

\tableofcontents

\section{Introduction}

\subsection{The algebraic theory of cyclotomic multiple zeta values}

Let $d$ be a positive integer, $(n_{i})_{d}=(n_{1},\ldots,n_{d})$ be a sequence of positive integers and  $(\xi_{i})_{d}=(\xi_{1},\ldots,\xi_{d})$ be a sequence of roots of unity in $\mathbb{C}$, with $(\xi_{d},n_{d}) \not= (1,1)$. The following complex number is called a cyclotomic multiple zeta value (abbreviated as MZV$\mu_{N}$, where $N$ is a positive integer such that $\xi_{i}^{N}=1$ for all $i$) :
\begin{equation} \label{eq:multizetas} \zeta \big( (n_{i})_{d};(\xi_{i})_{d} \big)  = \sum_{0<m_{1}<\ldots<m_{d}} \frac{\big( \frac{\xi_{2}}{\xi_{1}}\big)^{m_{1}} \ldots \big( \frac{1}{\xi_{d}}\big)^{m_{d}}}{m_{1}^{n_{1}} \ldots m_{d}^{n_{d}}} .
\end{equation}
The adjective cyclotomic and the notation $\mu_{N}$ are omitted if $N=1$. One says that $n=n_{1}+\ldots+n_{d}$ is the weight of $\big( (n_{i})_{d};(\xi_{i})_{d} \big)$ and that $d$ is its depth. Denoting by $(\epsilon_{n},\ldots,\epsilon_{1})=(\underbrace{0,\ldots,0}_{n_{d}-1},\xi_{d},\ldots,\underbrace{0,\ldots,0}_{n_{1}-1},
\xi_{1})$, we have
\begin{equation} \label{eq:multizetas integral} \zeta \big( (n_{i})_{d};(\xi_{i})_{d} \big) = (-1)^{d}\int_{0}^{1} \frac{dt_{n}}{t_{n}-\epsilon_{n}}  \int_{0}^{t_{n-1}} \ldots \int_{0}^{t_{2}} \frac{dt_{1}}{t_{1}-\epsilon_{1}}.
\end{equation}
\indent This shows that MZV$\mu_{N}$'s are Betti - de Rham periods of the pro-unipotent fundamental groupoid of $\mathbb{P}^{1} - \{0,\mu_{N},\infty\}$ (\cite{Deligne Goncharov}, \S5.16). 
\newline\indent Let $p$ be a prime number which does not divide $N$. $p$-adic cyclotomic multiple zeta values ($p$MZV$\mu_{N}$'s) are defined as $p$-adic analogues of the integrals in (\ref{eq:multizetas integral}) 
\cite{Deligne Goncharov} \cite{Furusho 1} \cite{Furusho 2} \cite{Yamashita} \cite{Unver MZV} \cite{U2}, \cite{I-1}, \cite{I-3}. They are elements $\zeta_{p,\alpha}\big((n_{i})_{d};(\xi_{i})_{d}\big)$ of the field $K_{p}$ generated by $\mathbb{Q}_{p}$ and a primitive $N$-th root of unity, where $\alpha \in \mathbb{Z} \cup \{\pm \infty\} - \{0\}$ is the number of iterations of the Frobenius of the crystalline pro-unipotent fundamental groupoid of $\mathbb{P}^{1} - \{0,\mu_{N},\infty\}$. They are reductions of $p$-adic periods by \cite{Yamashita}.
\newline\indent In the framework of $\pi_{1}^{\un}(\mathbb{P}^{1} - \{0,\mu_{N},\infty\})$, MZV$\mu_{N}$'s resp. $p$MZV$\mu_{N}$'s often appear via their non-commutative generating series $\Phi_{\KZ}$ resp. $\Phi_{p,\alpha}$, which is an element of the $\mathbb{C}$-algebra, resp. $K_{p}$-algebra of non-commutative power series over the formal variables $e_{x}$ where $x$ is $0$ or a $N$-th root of unity. The coefficient of $e_{0}^{n_{d}-1}e_{\xi_{d}}\ldots e_{0}^{n_{1}-1}e_{\xi_{1}}$ in $\Phi_{\KZ}$, resp. in $\Phi_{p,\alpha}$ is $(-1)^{d}\zeta\big((n_{i})_{d};(\xi_{i})_{d}\big)$, resp. $(-1)^{d}\zeta_{p,\alpha}\big((n_{i})_{d};(\xi_{i})_{d}\big)$.
\newline\indent According to the philosophy of periods, one wants to study the algebra generated by MZV$\mu_{N}$'s resp. $p$MZV$\mu_{N}$'s over the $N$-th cyclotomic field $k_{N}$, using the motivic nature of $\pi_{1}^{\un}(\mathbb{P}^{1} - \{0,\mu_{N},\infty\})$ established in \cite{Deligne Goncharov}. There are many known examples of families of polynomial equations in that algebra. Three among them are particularly meaningful and well-known : these are the regularized double shuffle equations, the associator equations, and the Kashiwara-Vergne equations. We will refer to them as the standard equations. At least in the $N=1$ case, it is conjectured that each of them generates all polynomial equations satisfied by MZV$\mu_{N}$'s, resp. $p$MZV$\mu_{N}$'s. Moreover, one can show that the set of solutions to these equations have structures of torsors for the Ihara action which is a byproduct of the motivic Galois action on $\pi_{1}^{\un}(\mathbb{P}^{1} - \{0,\mu_{N},\infty\})$ (\cite{Deligne Goncharov} \S5.12).

\subsection{A computation of $p$-adic cyclotomic multiple zeta values which keeps track of the motivic Galois action}

In \cite{I-1} \cite{I-2} \cite{I-3}, we have found formulas for $p$-adic cyclotomic multiple zeta values, which are analogues of the expression of sums of series (\ref{eq:multizetas}).
\newline\indent In these formulas, $\Phi_{p,\alpha}$ is involved via $\Ad_{\Phi_{p,\alpha}}(e_{1}) = \Phi_{p,\alpha}^{-1}e_{1}\Phi_{p,\alpha}$ and its images $\Ad_{\Phi_{p,\alpha}^{(\xi)}}(e_{\xi})$, $\xi \in \mu_{N}(K_{p})$, by the automorphisms $(z \mapsto \xi z)_{\ast}$ of $\pi_{1}^{\un,\dR}(\mathbb{P}^{1} - \{0,\mu_{N},\infty\})$. The correspondence between the coefficients of  $\{\Phi_{p,\alpha}^{(\xi)}\}^{-1}e_{\xi}\Phi^{(\xi)}_{p,\alpha}$ and the coefficients of $\Phi_{p,\alpha}$ will be discussed in a subsequent paper.
\newline\indent Concretely, computing $p$MZV$\mu_{N}$'s means computing the Frobenius structure of the KZ differential equation (equation (\ref{eq: nabla KZ})) associated with $\pi_{1}^{\un,\dR}(\mathbb{P}^{1} - \{0,\mu_{N},\infty\})$, whose solutions are called multiple polylogarithms (equation (\ref{eq:MPL})) and admit the following power series expansion at $0$ :
\begin{equation} \label{eq:Li 0} \Li\big((n_{i})_{d};(\xi_{i})_{d} \big)(z) = \sum_{0<m_{1}<\ldots<m_{d}} \frac{\big( \frac{\xi_{2}}{\xi_{1}} \big)^{m_{1}} \ldots \big( \frac{z}{\xi_{d}} \big)^{m_{d}}}{m_{1}^{n_{1}} \ldots m_{d}^{n_{d}}}.
\end{equation}
\noindent The computation of $p$MZV$\mu_{N}$'s arises as a characterization of the coefficients of each  $\Ad_{\Phi_{p,\alpha}^{(\xi)}}(e_{\xi})$ in terms of weighted cyclotomic multiple harmonic sums, which are essentially the coefficients above (below, $m$ is any positive integer and the other indices are as above) :
\begin{equation} \label{eq:mult har sums} \text{har}_{m} \big((n_{i})_{d} ;(\xi_{i})_{d+1} \big) = m^{n_{d}+\ldots+n_{1}} \sum_{0<m_{1}<\ldots<m_{d}<m}
\frac{\big( \frac{\xi_{2}}{\xi_{1}} \big)^{m_{1}} \ldots \big(\frac{\xi_{d+1}}{\xi_{d}}\big)^{m_{d}}
\big(\frac{1}{\xi_{d+1}}\big)^{m}}{m_{1}^{n_{1}}\ldots m_{d}^{n_{d}}}.
\end{equation}
\indent The computation gives a formula for the coefficients of $\Ad_{\Phi_{p,\alpha}^{(\xi)}}(e_{\xi})$, as sums of series whose terms are linear combinations over $k_{N}$ of the numbers $\text{har}_{p^{\alpha}} \big((n_{i})_{d} ;(\xi_{i})_{d+1} \big)$ which we call the prime weighted cyclotomic multiple harmonic sums (\cite{I-1}, Definition B.0.1). It also gives a converse of that formula, which is the following ; below, the notation $f[w]$ means the coefficient of a word $w$ in a non-commutative formal power series $f$, and the sum of series converges in $K_{p}$ :

\begin{multline} \label{eq:formula for n=1}
\har_{p^{\alpha}} \big( (n_{i})_{d};(\xi_{i})_{d+1} \big) =
(-1)^{d} \sum_{\xi \in \mu_{N}(K)} \sum_{l=0}^{\infty} \xi^{-p^{\alpha}}  \Ad_{\Phi_{p,\alpha}^{(\xi)}}(e_{\xi}) \Big[ e_{0}^{l} e_{\xi_{d+1}}e_{0}^{n_{d}-1}e_{\xi_{d}}\ldots e_{0}^{n_{1}-1}e_{\xi_{1}} \Big].
\end{multline}

\indent We have a total of three ways of viewing prime weighted multiple harmonic sums : (\ref{eq:Li 0}), (\ref{eq:mult har sums}) and (\ref{eq:formula for n=1}) ; this will make three frameworks of computation, which we represent respectively by the symbols $\int$, $\Sigma$ and $\int_{1,0}$. One of the central ideas in \cite{I-2} and \cite{I-3} was to express certain byproducts of the Frobenius both in a framework of integrals and the framework of series and to use that the two expressions were equal ; this gave formulas for $p$MZV$\mu_{N}$'s in terms of series.
\newline\indent Our computation of $p$MZV$\mu_{N}$'s keeps track of the motivic nature of $\pi_{1}^{\un}(\mathbb{P}^{1} - \{0,\mu_{N},\infty\})$ : the main formulas are expressed by means of new group actions which we call pro-unipotent harmonic actions, which are certain $p$-adic byproducts of the motivic Galois action and of the Ihara action evoked above.

\subsection{Relating the explicit formulas for $p$-adic cyclotomic multiple zeta values and their algebraic theory}

The purpose of this paper and of the two subsequent ones \cite{II-2} and \cite{II-3} is to relate the explicit formulas for $p$MZV$\mu_{N}$'s obtained in \cite{I-1}, \cite{I-2} and \cite{I-3} to the algebraic theory of $p$MZV$\mu_{N}$'s.
\newline\indent The main reason for doing it is the fact that our formulas for $p$MZV$\mu_{N}$'s keep track of the motivic Galois action (via the pro-unipotent harmonic actions), and that they are explicit. Another reason for doing it is Kaneko-Zagier's notion of finite multiple zeta values \cite{Kaneko Zagier}, which are in $(\underset{p\text{ prime}}{\prod} \mathbb{Z}/p\mathbb{Z}) \big/ (\underset{p\text{ prime}}{\bigoplus} \mathbb{Z}/p\mathbb{Z})$ \cite{Kaneko Zagier}, and Rosen's lift of finite multiple zeta values, called truncated multiple zeta values, which are elements of the complete topological ring $\varprojlim_{n} (\underset{p\text{ prime}}{\prod} \mathbb{Z}/p^{n}\mathbb{Z}) \big/ (\underset{p\text{ prime}}{\bigoplus} \mathbb{Z}/p^{n}\mathbb{Z})$ \cite{Rosen}. In this paper, certain of our results will recover certain known results on these two notions, and interpret them in terms of $\pi_{1}^{\un}(\mathbb{P}^{1} - \{0,1,\infty\})$. 
\newline\indent We are going to consider intrinsically the numbers which appear in the right-hand side of (\ref{eq:formula for n=1}) (as well as their complex analogues) and the numbers $\har_{p^{\alpha}}(w)$ involved in the formulas of \cite{I-1}, \cite{I-2}, \cite{I-3}, and turn them into notions called \emph{adjoint $p$-adic cyclotomic multiple zeta values} (Ad$p$MZV$\mu_{N}$'s) (Definition \ref{def adjoint}), and \emph{cyclotomic multiple harmonic values} (MHV$\mu_{N}$'s) (Definition \ref{def harmonic}). We are going to view them as two types of ``periods'' (in a generalized sense for the second case, as will be explained in \cite{II-3}).
Since each $\har_{p^{\alpha}}(w)$ is only an algebraic number, in order to turn it into an object which can be regarded as an interesting ``period'' comparable to MZV's, we are going to consider all $p$'s at the same time (as in the work of Kaneko-Zagier and Rosen), or all $\alpha$'s at the same time, or both ; and our cyclotomic multiple harmonic values will lift both Kaneko-Zagier's finite multiple zeta values and Rosen's truncated multiple zeta values. We will establish in \S2.2 the setting for computations on these objects, in the three frameworks $\int_{1,0}$, $\int$ and $\Sigma$.
\newline\indent We adopt the following principle, for this paper and all the subsequent ones : for each question on $p$MZV$\mu_{N}$'s which we want to study by means of explicit formulas, we find its analogue for the adjoint variant of $p$MZV$\mu_{N}$'s, we solve that adjoint variant, and it remains to pass from Ad$p$MZV$\mu_{N}$'s to $p$MZV$\mu_{N}$'s ; this last step is delayed to other works. In the present sequence of papers, we only consider the Ad$p$MZV$\mu_{N}$'s which are the natural version of $p$MZV$\mu_{N}$'s adapted to dealing with explicit formulas.
\newline\indent With the above ideas, the first step of our explicit version of the algebraic theory of $p$MZV$\mu_{N}$'s is thus to develop intrinsically the properties of adjoint $p$-adic cyclotomic multiple zeta values and cyclotomic multiple harmonic values, and this is the goal of this paper. In the end, our explicit version of the algebraic theory of $p$MZV$\mu_{N}$'s will be formulated as a ``comparison'' between the properties of Ad$p$MZV$\mu_{N}$'s and those of MHV$\mu_{N}$'s, where MHV$\mu_{N}$'s are defined by explicit formulas.

\subsection{Summary of the paper}

We are going to find ``adjoint'' and ``harmonic'' analogues of the standard properties of cyclotomic multiple zeta values.

The analogues of algebraic relations for cyclotomic multiple harmonic values which we find involve infinite sums of prime weighted cyclotomic multiple harmonic sums, which will be convergent both for the weight-adic topology and for the $p$-adic topology, uniformly with respect to $(p,\alpha)$, and with certain bounds on the $p$-adic norms of the rational coefficients. At first sight, the presence of infinite sums makes us leave the framework of algebraic geometry and periods. However, we are going to prove a close relation between the ``adjoint'' and the ``harmonic'' variants of the algebraic properties of MZV$\mu_{N}$'s, namely, a torsor structure using the pro-unipotent harmonic action defined in \cite{I-3}.

The paper is organized as follows. 

In \S2 we review prerequisites and we define formally the main objects that we are going to study : adjoint cyclotomic multiple zeta values, cyclotomic multiple harmonic values, and related objects. We discuss the meaning of the definitions and we establish the setting for making computations with these objects.

In \S3 we focus on double shuffle relations and we establish some natural adjoint and harmonic variants of them. Below, $\DS_{\mu}$ is the double shuffle scheme defined by Racinet in \cite{Racinet}, $\circ^{\smallint_{1,0}}$ is the Ihara action, $\circ^{\smallint_{1,0}}_{\Ad}$ is the adjoint Ihara action introduced in \cite{I-2}, Definition 1.1.3, $\circ_{\har}^{\smallint_{1,0}}$ is the pro-unipotent harmonic action introduced in \cite{I-3}, Definition 2.1.2 ; $\mathcal{O}^{\sh}$ is the shuffle Hopf algebra over the alphabet $e_{0\cup \mu_{N}}$, graded by the number of letters of words called their weight.
The map $\comp^{\har,\Ad}$ is defined in \S2.2.1.
\newline\newline
 \textbf{Theorem 1.} \emph{
\newline\indent (i) (adjoint double shuffle) For any $\mu$, there exists an explicit affine subscheme
$\DS_{\mu,\Ad}$ of $\Spec(\mathcal{O}^{\sh})$. One has a canonical morphism $\Ad(e_{1}) : \DS_{\mu} \rightarrow \DS_{\mu,\Ad}$, which is an isomorphism on its image, and sends Racinet's torsor to a torsor for the adjoint Ihara product.
\newline\indent In particular, the non-commutative generating series of Ad$p$MZV$\mu_{N}$'s resp. AdMZV$\mu_{N}$'s is in $\DS_{0,\Ad}(K_{p})$, resp. $\DS_{2i\pi,\Ad}(\mathbb{C})$.
\newline\indent (ii) (harmonic double shuffle) There exists an explicit affine ind-scheme $\DS_{\har}$, sub-ind-scheme of $\Spf(\mathcal{O}^{\sh})$, which can be obtained by each of the three frameworks : ``$\int_{1,0}$'', ``$\int$'' and ``$\Sigma$''.
\newline\indent The non-commutative generating series MHV$\mu_{N}$'s is in $\DS_{\har}((\prod_{p}K_{p})^{\mathbb{N}})$.
\newline\indent (iii) (relation between adjoint and harmonic) The map $\comp^{\har,\Ad}$ sends $\DS_{0,\Ad} \mapsto \DS_{\har}$ and its image is a torsor under the pro-unipotent harmonic action $\circ_{\har}^{\smallint_{1,0}}$ of the group  $(\Ad_{\DS_{0}}(e_{1}),\circ^{\smallint_{1,0}}_{\Ad})$.
\newline\indent (iv) More generally, adjoint double shuffle equations are satisfied by adjoint multiple polylogarithms and harmonic double shuffle equations are satisfied by harmonic multiple polylogarithms.}
\newline\newline
We note that in (ii) above, the frameworks of computations give harmonic double shuffle equations which look different at first sight but we can prove that they are equivalent.

In \S4 we focus on associator and Kashiwara-Vergne equations. We consider an analogy between the passage from associator equations to Kashiwara-Vergne equations constructed in \cite{AT} and \cite{AET} and the passage from double shuffle relations to adjoint double shuffle relations constructed in \S3. The main result is the following.
\newline
\newline
\textbf{Theorem 2.} (rough version)
\newline\indent \emph{(i) Kashiwara-Vergne equations arise naturally as a property of ($p$-adic, complex) adjoint MZV's rather than as a property of ($p$-adic, complex) MZV's and naturally amount to certain polynomial equations on adjoint ($p$-adic, complex) MZV's. 
\newline \indent (ii) There are equations satisfied by harmonic multiple polylogarithms which are of the same source with associator equations and Kashiwara-Vergne equations (functoriality of the KZ equation).}
\newline
\newline
In \S5, as a corollary, we explain that this paper recasts the study finite and symmetric (or symmetrized) multiple zeta values (see \cite{NoteCRAS}) as a particular case and a byproduct of the study of adjoint MZV's and MHV's, which is the theme arising naturally from the study of $p$MZV$\mu_{N}$'s via explicit formulas. Indeed, finite multiple zeta values are reductions of our multiple harmonic values modulo large primes, and symmetric multiple zeta values are a particular case of our adjoint multiple zeta values ; thus they satisfy a particular case of the equations of Theorem 1 and Theorem 2. We introduce finite cyclotomic multiple zeta values and finite multiple polylogarithms, and in which we will more generally relate our study of $p$MZV$\mu_{N}$'s via explicit formulas to these objects.

In \cite{II-2} we will explain how to recover certain properties of Ad$p$MZV$\mu_{N}$'s by the explicit formulas and the properties of MHV$\mu_{N}$'s, which will mostly answer to a question of Deligne and Goncharov. In \cite{II-3} we will formalize the fact that MHV$\mu_{N}$'s can be regarded as periods in a generalized sense, using motivic multiple zeta values. In \cite{III-1} we will study adjoint and harmonic distribution relations.
 
An open question is whether $\Ad(e_{1}) : \DS_{\mu} \rightarrow  \DS_{\mu,\Ad}$ is an isomorhism. We see this question as analogous to the conjecture made in \cite{AT} of an isomorphism relating associators and solutions to the Kashiwara-Vergne problem. 
 
\emph{Acknowledgments.} I thank Masanobu Kaneko, Ivan Marin and Pierre Cartier for useful discussions, and an anonymous referee whose remarks and suggestions enabled me to improve this paper. This work has been done at Universit\'{e} Paris Diderot with support of ERC grant 257638, Universit\'{e} de Strasbourg with support of Labex IRMIA, and Universit\'{e} de Gen\`{e}ve with support of NCCR SwissMAP.

\section{Definitions and setting for computations\label{review}}

We review some material on pro-unipotent fundamental groupoids (\S2.1), we define adjoint cyclotomic multiple zeta values, in the $p$-adic case (\S2.2) and in the complex case (\S2.3), we define cyclotomic multiple harmonic values (\S2.4) and their ``overconvergent'' variants (\S2.5). At the same time we establish the setting for making computations with all these objects and we show that replacing cyclotomic multiple zeta values by their adjoint variants does not change the algebra that they generate. Finally we define more generally adjoint and harmonic analogues of multiple polylogarithms (\S2.6), which admit respectively adjoint cyclotomic multiple zeta values and cyclotomic multiple harmonic values as special values.
 
In all this text, we denote by $\mathbb{N}$ resp. $\mathbb{N}^{\ast}$ the set of nonnegative, resp. positive integers ; $d$ and $n_{i}$ ($1 \leqslant i \leqslant d$) denote positive integers, $\xi_{i}$ ($1 \leqslant i \leqslant d$ or $1 \leqslant i \leqslant d+1$ depending on the context) are $N$-th roots of unity.

\numberwithin{equation}{subsection}

\subsection{Review on pro-unipotent fundamental groupoids}

\subsubsection{Generalities on the Betti and de Rham realizations of pro-unipotent fundamental groupoids}

Let $X$ be a smooth algebraic variety over a field $K$ of characteristic zero, with $X= \overline{X} - D$ where $\overline{X}$ is proper and smooth and $D$ is a normal crossings divisor.
\newline\indent The de Rham pro-unipotent fundamental groupoid $\pi_{1}^{\un,\dR}\big( X\big)$ (\cite{Deligne}, \S10.27, \S10.30 (ii)) is a groupoid over $X$ in the category of affine schemes over $K$, whose base-points are the points of $X$, and the points of the punctured tangent spaces at points of $D$ (\cite{Deligne}, \S15). Assuming that $H^{1}(\overline{X},\mathcal{O}_{\overline{X}})=0$, which holds in the examples of this paper, one also has a canonical base-point $\omega_{\dR}$ (\cite{Deligne}, \S12.4) with, for any couple of base-points $(x,y)$, an isomorphism of schemes $\pi_{1}^{\un,\dR}(X,y,x)\simeq \pi_{1}^{\un,\dR}(X,\omega_{\dR})$, and these isomorphisms are compatible with the groupoid structure. The bundle $\pi_{1}^{\un,\dR}(X,\omega_{\dR}) \times X$ carries the universal unipotent integrable connection on $X$.
\newline\indent Let us now assume that we have an embedding $K \hookrightarrow \mathbb{C}$. Then we also have the Betti pro-unipotent fundamental groupoid $\pi_{1}^{\un,\B}(X \times_{\Spec(K)} \Spec(\mathbb{C}))$, which is another groupoid in the category of affine schemes over $X \times_{\Spec(K)} \Spec(\mathbb{C})$, defined as the Malcev completion of the topological fundamental groupoid of $X(\mathbb{C})$ \cite{Deligne}.
\newline\indent Chen's theorem \cite{Chen} or, equivalently, the Riemann-Hilbert correspondence \cite{Deligne equa diff} restricted to unipotent objects, gives a natural isomorphism
\begin{equation} \label{eq:isomorphism} \pi_{1}^{\un,\B}(X) \times_{\Spec(K)} \Spec(\mathbb{C}) \simlra \pi_{1}^{\un,\dR}(X) \times_{\Spec(K)} \Spec(\mathbb{C}).
\end{equation}
Its coefficients are iterated path integrals on $X$ in Chen's sense \cite{Chen}. If $X$ is defined over a number field, they are periods, called the Betti-de Rham periods of $\pi_{1}^{\un}(X)$.

\subsubsection{The de Rham pro-unipotent fundamental groupoid of $\mathbb{P}^{1} - \{0,\mu_{N},\infty\}$}

Let $X=(\mathbb{P}^{1} - \{0,\mu_{N},\infty\}) / K$ where $K$ contains a primitive $N$-th root of unity.
\newline\indent By \cite{Deligne}, \S12, the affine scheme $\pi_{1}^{\un,\dR}(X,\omega_{\dR})$ is canonically isomorphic to $\Spec(\mathcal{O}^{\sh})\times_{\Spec(\mathbb{Q})} \Spec(K)$, where $\mathcal{O}^{\sh}$ is the shuffle Hopf algebra on the alphabet $e_{0\cup \mu_{N}} = \{e_{x}\text{ }|\text{ }x \in \{0\} \cup \mu_{N}(K)\}$. By definition, $\mathcal{O}^{\sh}$ is the $\mathbb{Q}$-vector space $\mathbb{Q}\langle e_{0\cup \mu_{N}} \rangle$ which admits as a basis the set of words on $e_{0 \cup \mu_{N}}$ including the empty word, graded by the number of letters of words called their weight, endowed with the following operations : the shuffle product $\sh$ defined by $(e_{i_{l+l'}}\ldots e_{i_{l+1}})\text{ }\sh\text{ }(e_{i_{l}} \ldots e_{i_{1}}) =
\sum\limits_{\sigma} 
e_{i_{\sigma^{-1}(l+l')}} \ldots e_{i_{\sigma^{-1}(1)}}$ where the sum is over permutations $\sigma$ of $\{1,\ldots,l+l'\}$ such that $\sigma(1)<\ldots<\sigma(l)$ and  $\sigma(l+1)<\ldots<\sigma(l+l')$ ; the deconcatenation coproduct $\Delta_{\dec}$ defined by $\Delta_{\dec}(e_{i_{l}}\ldots e_{i_{1}}) = \sum\limits_{l'=0}^{l} e_{i_{l}}\ldots e_{i_{l'+1}} \otimes e_{i_{l'}} \ldots e_{i_{1}}$ ; the counit $\epsilon$ equal to the augmentation map ; the antipode $S$ defined by $S(e_{i_{l}}\ldots e_{i_{1}}) = (-1)^{l} e_{i_{1}}\ldots e_{i_{l}}$.

\begin{Notation} Let $K\langle\langle e_{0 \cup \mu_{N}} \rangle\rangle$ be the $K$-algebra of non-commutative formal power series over the variables $e_{x}$, $x \in \{0\} \cup \mu_{N}(K)$ and, for $f\in K\langle\langle e_{0 \cup \mu_{N}} \rangle\rangle$ and $w$ a word over the alphabet
$\{e_{x}\text{ }|\text{ }x \in \{0\} \cup \mu_{N}(K)\}$, let $f[w] \in K$ be the coefficient of $w$ in $f$.
\end{Notation}

The group scheme $\Spec(\mathcal{O}^{\sh})$ is pro-unipotent, and the completed dual of the Hopf algebra $\mathcal{O}^{\sh}$ is $K \langle \langle e_{0\cup \mu_{N}}\rangle \rangle$ viewed as the Hopf algebra obtained as the completion of the universal enveloping algebra of the free Lie algebra in the variables $e_{x}$, $x \in \{0\} \cup \mu_{N}(K)$. We will denote by $\Delta_{\sh}$ its coproduct. We have 

\begin{equation} \label{eq:shuffle equation} \begin{array}{ll} \Spec(\mathcal{O}^{\sh})(K) & = \{ f \in K\langle\langle e_{0 \cup \mu_{N}} \rangle\rangle \text{ }|\text{ }\forall w,w' \text{ words on }e_{0 \cup \mu_{N}}, f[w\text{ }\sh\text{ }w']=f[w]f[w'],\text{ and }f[\emptyset] = 1 \}
\\ & = \{ f \in K\langle\langle e_{0 \cup \mu_{N}} \rangle\rangle \text{ }|\text{ }\Delta_{\sh}(f) = f \otimes f \},
\end{array}
\end{equation}
\begin{equation} \label{eq:shuffle equation modulo products}
\begin{array}{ll}
\Lie(\Spec(\mathcal{O}^{\sh})(K)) & = \{ f \in K \langle\langle e_{0 \cup \mu_{N}} \rangle\rangle \text{ }|\text{ }\forall w,w' \text{ words on }e_{0 \cup \mu_{N}},\text{ }f[w\text{ }\sh\text{ }w']=0\}
\\ & = \{ f \in K \langle\langle e_{0 \cup \mu_{N}} \rangle\rangle \text{ }|\text{ }\Delta_{\sh}(f) = f \otimes 1 + 1 \otimes f \}.
\end{array}
\end{equation}
\newline\indent The canonical connection on $\pi_{1}^{\un,\dR}(X,\omega_{\dR})\times X$ in the sense of \cite{Deligne}, \S12 is the Knizhnik-Zamolodchikov (KZ) connection :
\begin{equation} \label{eq: nabla KZ} \nabla_{\KZ} : f \mapsto df - \bigg(  e_{0} f \frac{dz}{z} + \sum_{\xi \in \mu_{N}(K)} e_{\xi} f \frac{dz}{z-\xi} \bigg).
\end{equation}
\indent 
\begin{Notation} \label{la premiere notation}
(\cite{Deligne Goncharov}, \S5), let $\Pi = \pi_{1}^{\un,\dR}(X,\omega_{\dR})$ and $\Pi_{1,0} = \pi_{1}^{\un,\dR}(X,-\vec{1}_{1},\vec{1}_{0})$ where $\vec{v}_{x}$ means the tangent vector $\vec{v}$ at $x$.
\end{Notation}

\subsubsection{Comparison between the Betti and de Rham pro-unipotent fundamental groupoid of $\mathbb{P}^{1} - \{0,\mu_{N},\infty\}$} 

Following \cite{Deligne Goncharov}, \S5.16, let $dch \in \pi_{1}^{\un,\B}(X,\vec{1}_{0},\vec{1}_{1})(\mathbb{C})$ be the image, by the Malcev completion map  
$  \pi_{1}(X,\vec{1}_{0},\vec{1}_{1}) \rightarrow  \pi_{1}^{\un,\B}(X,\vec{1}_{0},\vec{1}_{1})(\mathbb{C})$,
of the homotopy class of $\gamma : t \in [0,1] \mapsto t \in [0,1]$ ; let
\begin{equation} \label{eq:Phi KZ}
\Phi_{\KZ} = \comp_{\B,\dR}(dch) \in \Pi_{1,0}(\mathbb{C}),
\end{equation}
which appeared first in \cite{Drinfeld}, \S2 in the $N=1$ case ; $\Phi_{\KZ}$ is the non-commutative generating series of MZV$\mu_{N}$'s : indeed, given the above definition, the formula for cyclotomic multiple zeta values as iterated integrals (\ref{eq:multizetas}) amounts to :
\begin{equation} \label{eq:coefficient}
\zeta \big( (n_{i})_{d};(\xi_{i})_{d} \big) = (-1)^{d}\Phi_{\KZ}[e_{0}^{n_{d}-1}e_{\xi_{d}}\ldots e_{0}^{n_{1}-1}e_{\xi_{1}}].
\end{equation}
\indent Multiple polylogarithms, abbreviated as MPL's \cite{Goncharov} are multivalued holomorphic functions on $X(\mathbb{C})$, defined as iterated integrals in the sense of Chen \cite{Chen}, such that the non-commutative generating series 
$\Li = 1 + \sum\limits_{\substack{n \in \mathbb{N}^{\ast} \\ \epsilon_{1},\ldots,\epsilon_{n} \in \{0,1\}^{n}}} \Li\big( e_{\epsilon_{n}} \cdots e_{\epsilon_{1}} \big)e_{\epsilon_{n}} \cdots e_{\epsilon_{1}}$
defines a solution to the $\KZ$ equation (\ref{eq: nabla KZ}). For $\gamma$ a differentiable topological path on $\mathbb{P}^{1}(\mathbb{C})$ such that $\gamma\big((0,1)\big) \subset (\mathbb{P}^{1} - \{0,\mu_{N},\infty\})(\mathbb{C})$, and $\gamma'(0)\not= 0$ and $\gamma'(1)\not= 0$,
\begin{equation} \label{eq:MPL} \displaystyle\Li\big( e_{\epsilon_{n}} \cdots e_{\epsilon_{1}} \big)(\gamma) =
\int_{t_{n}=0}^{1} \gamma^{\ast}(\frac{dz}{z-\epsilon_{n}})(t_{n}) \int_{t_{n-1}=0}^{t_{n}} \ldots   \gamma^{\ast}(\frac{dz}{z-\epsilon_{2}})(t_{2}) \int_{t_{1}=0}^{t_{2}} \gamma^{\ast}(\frac{dz}{z-\epsilon_{1}})(t_{1}).
\end{equation}
When that integral diverges, (\ref{eq:MPL}) means the regularized iterated integral defined by considering the similar integral on $\gamma([\epsilon,1-\epsilon'])$, its asymptotic expansion when $\epsilon,\epsilon' \rightarrow 0$ which is in  $\mathbb{C}[[\epsilon,\epsilon']][\log(\epsilon),\log(\epsilon')]$, and finally the coefficient of this asymptotic expansion. It depends only of the homotopy class of $\gamma$, in the extended sense which includes tangential base-points, i.e. it depends on $\gamma'(0)$ and $\gamma'(1)$.
\newline\indent If we take $\gamma(0)=0$, the formal power series expansion of the iterated integrals (\ref{eq:MPL}) at $z=0$ is convergent for $|z|<1$ :
\begin{equation} \label{eq:multiple polylogarithms power series expansion}
\Li[e_{0}^{n_{d}-1}e_{\xi_{d}}\ldots e_{0}^{n_{1}-1}e_{\xi_{1}}](z) = (-1)^{d} \sum_{0<m_{1}<\ldots <m_{d}}
\frac{\big( \frac{\xi_{2}}{\xi_{1}} \big)^{m_{1}} \ldots \big(\frac{1}{\xi_{d}}\big)^{m_{d}}}{m_{1}^{n_{1}}\ldots m_{d}^{n_{d}}}.
\end{equation}

\subsubsection{The crystalline Frobenius of the de Rham pro-unipotent fundamental groupoid of $\mathbb{P}^{1} - \{0,\mu_{N},\infty\}$ over $K_{p}$}

The differential equation $\nabla_{\KZ}$ (\ref{eq: nabla KZ}) has a crystalline Frobenius structure over $K_{p}$ (\cite{Deligne}, \S11). 

The theory of Coleman integration, which relies on this Frobenius structure, enables to define $p$-adic analogues of MPL's and MZV$\mu_{N}$'s \cite{Furusho 1} \cite{Furusho 2} \cite{Yamashita}. In particular, one has $\Li_{p,X_{K}}^{\KZ}$ resp. $\Li_{p,X_{K}^{(p^{\alpha})}}^{\KZ}$, the generating series of $p$-adic multiple polylogarithms, Coleman functions characterized as the solution to $\nabla_{\KZ}$ resp. its pull-back $\nabla_{\KZ}^{(p^{\alpha})}$ by the Frobenius $\sigma$ of $K$ iterated $\alpha$ times, equivalent to $e^{\log_{p}(z)e_{0}}$ at $\vec{1}_{0}$ (defined in \cite{Furusho 1}, \cite{Furusho 2} for $N=1$ and \cite{Yamashita} for any $N$). 

We consider the Frobenius iterated $\alpha$ times $(\alpha \in \mathbb{N}^{\ast})$ in the sense of \cite{I-1}, \S1. 
Then, another type of $p$-adic analogue of MPL's and MZV$\mu_{N}$'s can be defined in a more ad hoc way, using canonical de Rham paths \cite{Deligne Goncharov} \cite{Unver MZV} \cite{U2} \cite{I-1} \cite{I-3}. In the end, with \cite{I-1} and \cite{I-3}, we have for each $\alpha \in \mathbb{Z}\cup \{\pm \infty\} - \{0\}$, an element $\Phi_{p,\alpha} \in \Pi_{1,0}(K)$, which characterizes the Frobenius at base-points $(\vec{1}_{1},\vec{1}_{0})$ iterated $\alpha$ times, where $\alpha=-\infty$ corresponds to Coleman integration. $p$-adic cyclotomic multiple zeta values are defined as its coefficients, namely :
$$ \zeta_{p,\alpha}\big((n_{i})_{d};(\xi_{i})_{d}\big) = (-1)^{d} \Phi_{p,\alpha}[e_{0}^{n_{d}-1}e_{\xi_{d}} \ldots e_{0}^{n_{1}-1}e_{\xi_{1}}]. $$

One also has $\Li_{p,\alpha}^{\dagger}$, the non-commutative generating series of overconvergent $p$-adic multiple polylgarithms (\cite{I-1}, \S1), which are overconvergent analytic functions on the affinoid analytic space $U_{0\infty}^{\an}=\mathbb{P}^{1,\an} - \underset{\xi \in \mu_{N}(K)}{\cup} \B(\xi,1)$, where $\B(x,r)$ means the open ball of center $x$ and radius $r$. It is characterized by $\Li_{p,\alpha}^{\dagger}(0)=1$ and the following differential equation (\cite{I-1}, Proposition 2.1)
\begin{equation}
\label{eq:horizontality equation}
d\Li_{p,\alpha}^{\dagger} = \bigg( p^{\alpha}e_{0}\omega_{0}(z) + \sum_{\xi \in \mu_{N}(K)} p^{\alpha} \omega_{\xi}(z) e_{\xi} \bigg) \Li_{p,\alpha}^{\dagger} - \Li_{p,\alpha}^{\dagger} \bigg( \omega_{0}(z^{p^{\alpha}})e_{0} + \sum_{\xi \in \mu_{N}(K)} \omega_{z_{0}^{p^{\alpha}}}(z^{p^{\alpha}}) \Ad_{\Phi^{(\xi)}_{p,\alpha}} (e_{\xi}) \bigg) 
\end{equation}
\noindent where $\Phi_{p,\alpha}^{(\xi)} = (x \mapsto \xi x)_{\ast}(\Phi_{p,\alpha})$. Equivalently, it is characterized by
\begin{multline} \label{eq:horizontality1}
\Li_{p,\alpha}^{\dagger}(z)(e_{0},(e_{\xi})_{\xi \in \mu_{N}(K)})
\times\Li_{p,X_{K}^{(p^{\alpha})}}^{\KZ}(z^{p^{\alpha}})\big(e_{0},(\Ad_{\Phi^{(\xi)}_{p,\alpha}}(e_{\xi}))_{\xi \in \mu_{N}(K)} \big) 
\\ = \Li_{p,X_{K}}^{\KZ}(z) \big(p^{\alpha}e_{0},(p^{\alpha}e_{\xi})_{\xi \in \mu_{N}(K)} \big)
\end{multline}
We note that $\Li_{p,X_{K}}^{\KZ}$ and $\Li_{p,X_{K}^{(p^{\alpha})}}^{\KZ}$ depend on the choice of a branch of the $p$-adic logarithm, but not $\Li_{p,\alpha}^{\dagger}$, nor $p$MZV$\mu_{N}$'s.

\subsubsection{The motivic pro-unipotent fundamental groupoid of $\mathbb{P}^{1} - \{0,\mu_{N},\infty\}$}

The motivic pro-unipotent fundamental groupoid $\pi_{1}^{\un,\text{mot}}(\mathbb{P}^{1} - \{0,\mu_{N},\infty\})$ is defined and studied in \cite{Deligne Goncharov}. Let $G_{\omega}$ be the fundamental group associated with the Tannakian category of mixed Tate motives over $k_{N}$ which are unramified at primes $p$ prime to $N$, and the canonical fiber functor $\omega$. It is a semi-direct product $G_{\omega}= \mathbb{G}_{m} \ltimes U_{\omega}$ where $U_{\omega}$ is pro-unipotent. It acts on $\Pi_{1,0}$, and this action encodes the algebraic theory of MZV$\mu_{N}$'s and $p$MZV$\mu_{N}$'s according to the conjecture of periods (\cite{Deligne Goncharov}, \S5). The action of $\mathbb{G}_{m}$ encodes the weight grading and is 
\begin{equation} \label{eq:tau} \tau : \begin{array}{cc} \mathbb{G}_{m} \times \Pi_{1,0} \rightarrow \Pi_{1,0}
\\ \big( \lambda,f(e_{0},(e_{\xi})_{\xi \in \mu_{N}(K)})\big) \mapsto f(\lambda e_{0},(\lambda e_{\xi})_{\xi \in \mu_{N}(K)})
\end{array}.
\end{equation}
The action of $U_{\omega}$ has been computed by Goncharov \cite{Goncharov}. The image of this action by a certain morphism is isomorphic to the Ihara product (our notation below is not standard)
\begin{equation} \label{eq:Ihara} \circ^{\smallint_{1,0}} :
\begin{array}{cc} \Pi_{1,0} \times \Pi_{1,0} \rightarrow \Pi_{1,0}
\\
(g,f) \mapsto g \circ^{\smallint_{1,0}} f = g(e_{0},(e_{\xi})_{\xi \in \mu_{N}(K)}) \times f\big(e_{0},(\Ad_{g^{(\xi)}}(e_{\xi}))_{\xi \in \mu_{N}(K)}\big)
\end{array} .
\end{equation}

\subsubsection{The Betti and de Rham pro-unipotent fundamental groupoid of a more general $\mathbb{P}^{1} - D$}

By \cite{Deligne}, the description of the Betti and de Rham realizations of $\pi_{1}^{\un}(\mathbb{P}^{1} - \{0,\mu_{N},\infty\})$ of \S2.1.2 remains true in the case of an arbitrary punctured projective line $\mathbb{P}^{1} - D$ over a subfield of $\mathbb{C}$, provided that we replace the alphabet $e_{0 \cup \mu_{N}}$ by the alphabet $\{e_{x}\text{ }|\text{ }x \in D - \{\infty\}\}$. It can also be deduced from the description of $\pi_{1}^{\un}(\mathcal{M}_{0,n})$ of \S2.1.3.
We write the words on that alphabet on the form $e_{0}^{n_{d}-1}e_{z_{d}}\ldots e_{0}^{n_{1}-1}e_{z_{1}}e_{0}^{n_{0}-1}$ with $d$ and $n_{i}$ positive integers and $z_{i} \in D - \{0,\infty\}$. For most computations, we can restrict to words such that $n_{d}\geqslant 2$ and $n_{0}=1$. The word above is then denoted also by $((n_{i})_{d},(z_{i})_{d})$.
\newline\indent The multiple polylogarithms on $\mathcal{M}_{0,n}$ (\S2.1.3) induce multiple polylogarithms on $\mathbb{P}^{1} - \{0,x_{1},\ldots,x_{r},\infty\}$, whose power series expansion are given by (\ref{eq:Li series bis}) as functions of $z$, which are solution to the KZ equation, generalization of (\ref{eq: nabla KZ}) :
\begin{equation} \label{eq: nabla KZ prime} \nabla_{\KZ} : f \mapsto df - \bigg(  e_{0} f \frac{dz}{z} + \sum_{\xi \in \mu_{N}(K)} e_{\xi} f \frac{dz}{z-\xi} \bigg) .
\end{equation}
The weighted multiple harmonic sums, obtained by the coefficients of the power series expansion (\ref{eq:Li series bis}), are now the following numbers (where the $z_{i}$'s are in $D - \{0,\infty\}$) :
$$ \har_{m}((n_{i})_{d},(z_{i})_{d+1}) = 
\sum_{0<m_{1} <\ldots < m_{d}<m} 
\frac{\big( \frac{z_{2}}{z_{1}} \big)^{m_{1}} \ldots \big(\frac{z_{d+1}}{z_{d}}\big)^{m_{d}} \big(\frac{1}{z_{d+1}}\big)^{m}}{m_{1}^{n_{1}} \ldots m_{d}^{n_{d}}} . $$

\subsubsection{The Betti and de Rham pro-unipotent fundamental groupoid of the moduli spaces $\mathcal{M}_{0,n}$}

Let, for $n \in \mathbb{N}^{\ast}$, the scheme
$\mathcal{M}_{0,n+3} = \{ (x_{1},\ldots,x_{n+3}) \in (\mathbb{P}^{1})^{n+3} \text{ } |\text{ } x_{i}\not= x_{j} \text{ for all i} \not= \text{j} \} / \PGL_{2}$ over $\mathbb{Q}$, and let $\overline{\mathcal{M}}_{0,n+3}$ be its Deligne-Mumford compactification, which is a smooth projective variety such that $\overline{\mathcal{M}}_{0,n+3} - \mathcal{M}_{0,n+3}$ is a normal crossings divisor. The $x_{i}$'s are called the canonical coordinates. The homography of $\mathbb{P}^{1}$ which sends $(x_{n+1},x_{n+2},x_{n+3})$ to $(0,1,\infty)$ induces an isomorphism between $\mathcal{M}_{0,n+3}$ and the affine variety
$\{(y_{1},y_{2},\ldots,y_{n}) \in (\mathbb{P}^{1} - \{0,1,\infty\})^{n} \text{ }|\text{ for all i,j}, y_{i} \not= y_{j} \}$, defined by $\displaystyle y_{i} = \frac{x_{i} - x_{n+1}}{x_{i}-x_{n+3}} \frac{x_{n+2} - x_{n+3}}{x_{n+2}-x_{n+1}}$, $(1 \leqslant i \leqslant n)$ ; the $y_{i}$'s are called simplicial coordinates. In particular, this gives $\mathcal{M}_{0,4} \simeq \mathbb{P}^{1} - \{0,1,\infty\}$, $\overline{\mathcal{M}_{0,4}} =\mathbb{P}^{1}$ ; $\overline{\mathcal{M}_{0,5}}$ is obtained blowing up $(\mathbb{P}^{1})^{2} \supset \mathcal{M}_{0,5}$ at the three points where $(\mathbb{P}^{1})^{2} - \mathcal{M}_{0,5}$ is not normal crossings, namely, in simplicial coordinates, $(0,0)$, $(1,1)$ and $(\infty,\infty)$.
\newline\indent By \cite{Deligne} \S12, $\Lie \big(\pi_{1}^{\un,\dR}(\mathcal{M}_{0,n+3},\omega_{\dR}) \big)$ is the pro-nilpotent Lie algebra with generators $e_{ij}$, $1 \leq i \not= j \leqslant n+3$,
\noindent and relations $e_{ij} = e_{ji}$, for all $i,j$, 
$\sum\limits_{j=1}^{n} e_{ij} = 0$ for all $i$, and $[e_{ij},e_{kl}] = 0$ for all $i,j,k,l$ pairwise distinct. This determines the pro-unipotent affine group scheme $\pi_{1}^{\un,\dR}(\mathcal{M}_{0,n},\omega_{\dR})$.
\newline\indent The canonical connection on $\pi_{1}^{\un,\dR}(\mathcal{M}_{0,n},\omega_{\dR}) \times \mathcal{M}_{0,n}$ in the sense of \cite{Deligne}, \S12, is the KZ connection in several variables : 
\begin{equation} \label{eq:nablaKZ M0,n} \nabla_{\KZ} : f \mapsto df - \sum_{1 \leqslant i<j \leqslant n+3} e_{ij} d\log(x_{i}-x_{j}) f ,
\end{equation}
i.e., in the cubic coordinates $c_{i}$ defined by 
$y_{i} = c_{i} \ldots c_{n}$, $(1 \leqslant i \leqslant n)$,
\begin{equation} \nabla_{\KZ} : f \mapsto df
- \bigg( \sum_{u=1}^{r} \frac{dc_{u}}{c_{u}} \sum_{u\leqslant i<j \leqslant n} e_{ij}
- \sum_{\substack{1\leqslant v \leqslant v' \leqslant n \\  2 \leqslant i \leqslant j }} \frac{d(c_{v} \ldots c_{v'})}{c_{v} \ldots c_{v'} -1} e_{v-1,v'}  
-  \sum_{\substack{1\leqslant v \leqslant v' \leqslant n \\ 1 \leqslant i \leqslant j}} \frac{d(c_{v} \ldots c_{v'})}{c_{v} \ldots c_{v'} -1} e_{v',n-1} \bigg) f.
\end{equation}
\indent 
The groupoid $\pi_{1}^{\un,\B}(\mathcal{M}_{0,n})$ can be computed as follows by induction on $n$ \cite{FR}. Each forgetful map $\mathcal{M}_{0,n+1} \rightarrow \mathcal{M}_{0,n}$ is a fibration and induces a long exact sequence in homotopy ; given that the $\mathcal{M}_{0,n}$'s are $K(\pi,1)$ spaces, that long exact sequence amounts to a short exact sequence :
$1 \rightarrow \pi_{1}^{\text{top}}\big(F(\mathbb{C}) \big)
\rightarrow \pi_{1}^{\text{top}}(\mathcal{M}_{0,n+1}(\mathbb{C}))
\rightarrow \pi_{1}^{\text{top}}(\mathcal{M}_{0,n}(\mathbb{C})) \rightarrow 1$ where $F$ is the fiber of the forgetful map ; moreover, that short exact sequence is split. By the exactness of the functor of Malcev completion, it gives rise to a short split exact sequence $1 \rightarrow \pi_{1}^{\un,\B}\big(F(\mathbb{C}) \big)
\rightarrow \pi_{1}^{\un,\B}(\mathcal{M}_{0,n+1}(\mathbb{C}))
\rightarrow \pi_{1}^{\un,\B}(\mathcal{M}_{0,n}(\mathbb{C})) \rightarrow 1$.
\newline\indent Multiple polylogarithms in several variables are defined as follows \cite{Goncharov}. One consider first the following family of iterated integrals : for $a_{0},\ldots,a_{n+1} \in \mathbb{C}$, and $\gamma$ a path in $\mathbb{C} - \{a_{1},\ldots,a_{n}\}$ such that $\gamma(0)=a_{0}$ and $\gamma(1)=a_{n+1}$ :
\begin{equation} \label{eq:Li} I\big( a_{n+1} ; a_{n},\ldots,a_{1} ; a_{0} \big)(\gamma) = \int_{t_{n}=0}^{1}\gamma^{\ast}\bigg(\frac{dz}{z-a_{n}}\bigg)(t_{n}) \int_{t_{n-1}=0}^{t_{n}} \ldots \gamma^{\ast}\bigg(\frac{dz}{z-a_{2}}\bigg)(t_{2}) \int_{t_{1}=0}^{t_{2}} \gamma^{\ast}\big( \frac{dz}{z-a_{1}}\big)(t_{1}) .
\end{equation}
One usually writes $(a_{n},\ldots,a_{1})$ as 
$(\overbrace{0,\ldots,0}^{n_{d}-1},z_{d},\ldots,\overbrace{0,\ldots,0}^{n_{1}-1},z_{1},\overbrace{0,\ldots,0}^{n_{0}-1})$ and $a_{n+1}=z$, and an affine change of variable allows to assume that $a_{0}=0$. Finally one writes $(x_{1},\ldots,x_{d})= (\frac{z_{2}}{z_{1}},\ldots,\frac{z}{z_{d}})$. The functions obtained after these transformations are multiple polylogarithms in several variables. They define a solution to the KZ equation (\ref{eq:nablaKZ M0,n}). If $n_{0}=1$, one writes
$(\underbrace{0,\ldots,0}_{n_{d}-1},z_{d},\ldots,\underbrace{0,\ldots,0}_{n_{1}-1},z_{1})=\big( (n_{i})_{d}; (z_{i})_{d} \big)$. MPL's have the following power series expansion :
\begin{equation} \label{eq:Li series bis} \Li\big( (n_{i})_{d};(z_{i})_{d} \big)(z) = \sum_{0<m_{1} <\ldots < m_{d}}
\frac{\big( \frac{z_{2}}{z_{1}} \big)^{m_{1}} \ldots \big(\frac{z}{z_{d}}\big)^{m_{d}}}{m_{1}^{n_{1}}\ldots m_{d}^{n_{d}}} .
\end{equation}

\subsection{Adjoint $p$-adic cyclotomic multiple zeta values}

For each prime number $p$ which does not divide $N$, $K_{p}$ is the extension of $\mathbb{Q}_{p}$ generated by the $N$-th roots of unity in $\overline{\mathbb{Q}_{p}}$.
\newline\indent Let $\mathcal{O}^{\ast}$ be the $\mathbb{Q}$-vector space generated by the empty word and the words of the form
$\big((n_{i})_{d};(\xi_{i})_{d} \big)$, identified to elements of $\mathcal{O}^{\sh}$ (defined in \S2.1.2) by 
$\big((n_{i})_{d};(\xi_{i})_{d} \big) = e_{0}^{n_{d}-1}e_{\xi_{d}} \ldots e_{0}^{n_{1}-1}e_{\xi_{1}}$.

\subsubsection{Definition}

In addition with the notations above, $b$ is a non-negative integer, and $\Lambda$ is a formal variable.

\begin{Definition} \label{def adjoint}(i) Let adjoint $p$-adic cyclotomic multiple zeta values (Ad$p$MZV$\mu_{N}$'s) be the numbers
\begin{multline} \label{eq:zeta adjoint}
\zeta^{\Ad}_{p,\alpha} \big( (n_{i})_{d};b;(\xi_{i})_{d+1}\big) = (-1)^{d} \sum_{\xi \in \mu_{N}(K)} \xi^{-p^{\alpha}} \Ad_{\Phi^{(\xi)}_{p,\alpha}}(e_{\xi}) \big[e_{0}^{b}e_{\xi_{d+1}}e_{0}^{n_{d}-1}e_{\xi_{d}}\ldots e_{0}^{n_{1}-1}e_{\xi_{1}} \big]
\\ = (-1)^{d} \sum_{d'=0}^{d} \bigg( \prod_{i=d'+1}^{d} {-n_{i} \choose l_{i}} \bigg)\text{ }\xi_{d'}^{-p^{\alpha}}\text{ }
\zeta_{p,\alpha}^{(\xi_{d'})}\big( (n_{d-i}+l_{d-i});(\xi_{d+1-i}) \big)_{0 \leq i \leq d-d'}\text{ }
\zeta_{p,\alpha}^{(\xi_{d'})}\big((n_{i});(\xi_{i})\big)_{1 \leq i \leq d'-1} .
\end{multline}
\noindent (ii) Let the $\Lambda$-adjoint $p$-adic cyclotomic multiple zeta values ($\Lambda$Ad$p$MZV$\mu_{N}$'s) be the following power series :
\begin{multline} \zeta^{\Lambda \Ad}_{p,\alpha} \big( (n_{i})_{d};(\xi_{i})_{d+1} \big) = \Lambda^{n_{1}+\ldots+n_{d}} \sum_{b=0}^{\infty} \Lambda^{b} \zeta_{p,\alpha}^{\Ad} \big((n_{i})_{d};b;(\xi_{i})_{d+1} \big) 
\\
=  (-1)^{d} \sum_{\xi \in \mu_{N}(K)} \xi^{-p^{\alpha}} \Ad_{\Phi^{(\xi)}_{p,\alpha}}(e_{\xi})
\bigg[ \frac{\Lambda^{n_{1}+\ldots+n_{d}}}{1-\Lambda e_{0}} e_{\xi_{d+1}}e_{0}^{n_{d}-1}e_{\xi_{d}} \ldots e_{0}^{n_{1}-1}e_{\xi_{1}} \bigg] .
\end{multline}
\end{Definition}

In the case of $\mathbb{P}^{1} - \{0,1,\infty\}$, adjoint $p$-adic multiple zeta values are the numbers
$\zeta^{\Ad}_{p,\alpha}((n_{i})_{d};b) = (\Phi_{p,\alpha}^{-1}e_{1}\Phi_{p,\alpha})[e_{0}^{b}e_{1}e_{0}^{n_{d}-1}e_{1}\ldots e_{0}^{n_{1}-1}e_{1}] \in \mathbb{Q}_{p}$. 
\newline\indent Let $\mathcal{O}^{\ast}_{\Ad}$ be the $\mathbb{Q}$-vector space generated by the words of the form $ \big( \big((n_{i})_{d};b;(\xi_{i})_{d+1} \big)\big)$ : $b$ is written separately because it will play a particular role. Moreover, we will see that the $b=0$ and $b>0$ cases will sometimes have different properties (the former corresponds to finite cyclotomic multiple zeta values, see \S6), as already suggested by the formulas for $p$MZV$\mu_{N}$'s in \cite{I-2} in which this distinction exists.

\begin{Definition} Let $K\langle\langle Y_{N}^{\Ad} \rangle\rangle$ be the set of non-commutative formal power series of the following type : $f=\sum\limits_{w=\big( (n_{i})_{d};(\xi_{i})_{d+1}\big)}\sum\limits_{b\in \mathbb{N}} f[w;b](w;b)$ with $f[w;b] \in K$.
\end{Definition}

\subsubsection{Relation with $p$-adic cyclotomic multiple zeta values}

We now show that $p$MZV$\mu_{N}$'s and Ad$p$MZV$\mu_{N}$'s generate the same $k_{N}$-algebra, compatibly with the weight and depth.

Let $K$ be a field of characteristic $0$. Let $\Theta(K)$ be the group of characters of $\mu_{N}$ with values in $K^{\ast}$, i.e. the set of multiplicative morphisms $\mu_{N}(K) \rightarrow (K^{\ast},\ast)$. Let $\mathcal{O}^{\sh}_{n,\leq d}$
be the $k_{N}$-vector space generated by the elements $\sh_{i=1}^{r}w_{i}$ where $r\geq 1$ and $w_{i}$ are words such that $\sum\limits_{i=1}^{r}\weight(w_{i})=n$ and $\sum\limits_{i=1}^{r}\depth(w_{i})\leq d$. Below we implicitly view $\Phi$ and $\Phi_{\Ad}$ as functions $\mathcal{O}^{\sh} \rightarrow K$.

\begin{Definition} (i) For any $\chi \in \Theta$, let $\Moy_{\chi} : K \langle\langle e_{0\cup\mu_{N}} \rangle\rangle \rightarrow K \langle\langle e_{0\cup\mu_{N}} \rangle\rangle$ the map which sends $f \mapsto \sum\limits_{\xi \in \mu_{N}(K)} \chi(\xi) f^{(\xi)}$.
	
	(ii) For any $\Phi \in \tilde{\Pi}_{1,0}(K)$, 
	let $\Phi_{\Ad,\chi} = \Moy_{\chi} \Ad_{\Phi}(e_{1}) = \sum\limits_{\xi \in \mu_{N}(K)} \chi(\xi) {\Phi^{(\xi)}}^{-1}e_{\xi}\Phi^{(\xi)}$.
\end{Definition}

\begin{Proposition} The maps $\tilde{\Pi}_{1,0}(K) \rightarrow K \langle \langle e_{0 \cup \mu_{N}} \rangle\rangle$, $\Phi \mapsto \Ad_{e_{1}}(\Phi)$ and $\tilde{\Pi}_{1,0}(K) \times \Theta(K) \rightarrow \Theta \langle \langle e_{0 \cup \mu_{N}} \rangle\rangle$, $(\Phi,\chi) \mapsto \Phi_{\Ad,\chi}$, are injective. Moreover, for all $n \geq d\geq 0$, we have $\Phi(\mathcal{O}^{\sh}_{n,\leq d}) = \Phi_{\Ad,\chi}(\mathcal{O}^{\sh}_{n+1,\leq d+1})$.
\end{Proposition}

\begin{proof} The injectivity of $\Phi \mapsto \Phi^{-1}e_{1}\Phi$ follows from $K \langle \langle e_{1}\rangle\rangle \cap \tilde{\Pi}_{1,0}(K)= \{1\}$ and the implication $f e_{1} = e_{1}f \Rightarrow f \in K \langle \langle e_{1}\rangle\rangle$. That implication is proved as follows : if $w$ is a word which contains at least a letter different from $e_{1}$, one shows that $f[w]=0$, by writing $w=e_{1}^{n}e_{x}z$ with $x \not=1$ and by an induction on $n$.
	
Let us prove the rest of the statement. We start with a few preliminary properties :
\noindent\newline\indent (a) for all $x_{1},\ldots,x_{n} \in \{0\} \cup \mu_{N}$, $\Phi^{(\xi)}[e_{x_{n}} \ldots e_{x_{1}}] = \Phi [e_{\xi^{-1}x_{n}} \ldots e_{\xi^{-1}x_{1}}] $. Thus, for all $n\geq d \geq 1$, the $k_{N}$-vector space generated by the $\Phi^{(\xi)}[w]$'s with $w$ a word of weight $n$ and depth $d$ is independent of $\xi$.
\newline\indent (b) For $f \in \tilde{\Pi}_{1,0}(K)$, writing $f^{-1}f=1$ we get by induction that, for all $n\geq d \geq 1$,
we have $f(\mathcal{O}^{\sh}_{n,\leq d}) = f^{-1}(\mathcal{O}^{\sh}_{n,\leq d})$ and we also have $f^{-1}[w] \equiv - f[w] \mod f(\mathcal{O}^{\sh}_{n,\leq d-1})$ for all words $w$.
\newline\indent (c) For $f$ a solution to the shuffle equation or to the shuffle equation modulo products and such that $f[e_{0}] = 0$, any $f[w]$ with $w$ a word of weight $n$ and depth $d$ is a $\mathbb{Z}$-linear combination of the $f[w']$ with $w'$ a word of weight $n$ and depth $d$ of the form $w'=e_{0}w''e_{\xi}$ with $\xi \in \mu_{N}(K)$.
\newline\indent (d) By (a) and the shuffle equation, for $\Phi \in \tilde{\Pi}_{1,0}(K)$, we have $\Phi^{(\xi)}[e_{0}^{l}] = \Phi^{(\xi)}[e_{\xi}] = 0$ for all $\xi \in \mu_{N}(K)$ and $l\geq 1$.
\newline\indent The inclusion $\Phi(\mathcal{O}^{\sh}_{n,\leq d}) \supset \Phi_{\Ad}(\mathcal{O}^{\sh}_{n+1,\leq d+1})$ is clear by the shuffle equation for $\Phi$, and we want to prove the converse inclusion and the injectivity, by an induction on $d$.
\newline\indent The only non-zero coefficient in $\Phi$ in depth $0$ is the coefficient of the empty word (weight 0), equal to 1. The only non-zero coefficients in $\Phi_{\Ad,\chi}$ in depth $\leq 1$ are the coefficients $\Phi_{\Ad,\chi}[e_{\xi}]$ (weight 1), equal to $\chi(\xi)$. They are in $k_{N} \subset K$. This determines $\chi$ in terms of $\Phi_{\Ad,\chi}$. We have $\Phi(\mathcal{O}^{\sh}_{n,\leq 0}) = \Phi_{\Ad,\chi}(\mathcal{O}^{\sh}_{n+1,\leq 1})$ for all $n$. This proves the result for $d=0$.
\newline\indent Now for any $d>0$, we consider a word of the form 
$w=e_{0}^{l}e_{\xi_{d+1}}e_{0}^{n_{d}-1} \ldots e_{0}^{n_{2}-1}e_{\xi_{2}}e_{\xi_{1}}$ with $l>0$, and the $\xi_{i}$'s in $\mu_{N}(K)$ (i.e. the special case $n_{1}=1$ and $l>0$ in usual words). We have by (d)
$$ \Phi_{\Ad,\chi}[e_{0}^{l}e_{\xi_{d+1}}e_{0}^{n_{d}-1} \ldots e_{0}^{n_{2}-1}e_{\xi_{2}}e_{\xi_{1}}] \equiv \chi(\xi_{0}) {\Phi^{(\xi_{0})}}^{-1}[e_{0}^{l}e_{\xi_{d+1}}e_{0}^{n_{d}-1} \ldots e_{0}^{n_{2}-1}e_{\xi_{2}}] \mod \Phi(\mathcal{O}^{\sh}_{n,\leq d-1}) $$
By induction on $d$ and (a), (b), (c) above, this determines $\Phi^{(\xi)}$, thus $\Phi$, in depth $\leq d$, in terms of $\Phi_{\Ad}$, thus it proves the injectivity in depth $\leq d$, and this also proves the inclusion $\Phi(\mathcal{O}^{\sh}_{n,\leq d}) \subset \Phi_{\Ad,\chi}(\mathcal{O}^{\sh}_{n,\leq d+1})$.
\end{proof}

\begin{Corollary} The $p$MZV$\mu_{N}$'s and Ad$p$MZV$\mu_{N}$'s generate the same $k_{N}$-algebra. 
	
More precisely,	for any $n \geq d\geq 1$, the two following vector spaces are equal : the $k_{N}$-vector space of generated by products 
$\prod_{i=1}^{r}\zeta^{\Ad}_{p,\alpha}(w)$ with $r \geq 1$, 
$\sum\limits_{i=1}^{r} \weight(w_{i}) = n$ and $\sum\limits_{i=1}^{r} \depth(w_{i}) \leq d$, and The $k_{N}$-vector space generated by products  $\prod_{i=1}^{r}\zeta^{\Ad}_{p,\alpha}(w)$ with $r \geq 1$, 
$\sum\limits_{i=1}^{r} (\weight(w_{i})-1) = n$ and $\sum\limits_{i=1}^{r} (\depth(w_{i})-1) \leq d$.
\end{Corollary}

\begin{proof} We apply Proposition 2.2.4 to $\Phi = \Phi_{p,\alpha}$ and $\chi : \xi \mapsto \xi^{-p^{\alpha}}$.
\end{proof}

\subsection{Adjoint complex cyclotomic multiple zeta values}

\subsubsection{Definition}
 
The most direct complex analogue of Ad$p$MZV$\mu_{N}$'s (Definition \ref{def adjoint}) is the following.

Let a character $\chi : \mu_{N}(\mathbb{C}) \rightarrow \mathbb{C}^{\ast}$ and let $\Phi_{\KZ,\Ad,\chi} = \sum\limits_{\xi \in \mu_{N}(K)} \chi(\xi) {\Phi_{\KZ}^{-1}}^{(\xi)}e_{\xi}{\Phi_{\KZ}}^{(\xi)}$.

\begin{Definition} \label{def adjoint complex}(i) Let the adjoint cyclotomic multiple zeta values be the numbers : 
	\begin{equation} \label{eq:complex01} \zeta^{\Ad}\big( (n_{i})_{d};(\xi_{i})_{d+1};l;\chi) = \Phi_{\KZ,\Ad,\chi} [e_{0}^{l}e_{\xi_{d+1}}e_{0}^{n_{d}-1}e_{\xi_{d}} \ldots e_{0}^{n_{1}-1}e_{\xi_{1}}] .
	\end{equation}
	(ii) Let the $\Lambda$-adjoint cyclotomic multiple zeta values be the numbers : 
	\begin{equation} \label{eq:complex02} \zeta^{\Lambda \Ad}((n_{i})_{d};(\xi_{i})_{d+1};\chi) = \Phi_{\KZ,\Ad,\chi} \big[ \frac{1}{1-\Lambda e_{0}}e_{\xi_{d+1}}e_{0}^{n_{d}-1}e_{\xi_{d}} \ldots e_{0}^{n_{1}-1}e_{\xi_{1}} \big] .
	\end{equation}
\end{Definition}

Whereas in the crystalline setting, one has the Frobenius automorphism, whose restriction to $\Lie(\Pi_{0,0})$ sends $(e_{0},e_{1}) \mapsto (\frac{e_{0}}{p},\Phi_{p,1}^{-1}\frac{e_{1}}{p}\Phi_{p,1})$, in the Betti-de Rham setting, one has the automorphism of $\Pi_{0,0}$ which expresses the monodromy of the KZ connection (\ref{eq: nabla KZ}), which sends $(e^{e_{0}},e^{e_{1}}) \mapsto (e^{2i\pi e_{0}}, \Phi_{\KZ}^{-1}e^{2i\pi e_{1}}\Phi_{\KZ})$. Let $\comp^{B,DR}$ be the Betti-de Rham comparison isomorphism of $\Pi_{0,0}$ ; let $\gamma$ be the straight path $[0,1] \rightarrow \mathbb{C}$, $t\mapsto t$ ; let $c_{0}$, resp. $c_{1}$ be a simple loop around $0$, resp. $1$ positively oriented ; we have 
$$ (e^{2i\pi e_{0}},\Phi_{\KZ}^{-1}e^{2i\pi e_{1}}\Phi_{\KZ}) = \big(\comp^{\B,\dR}(c_{0}), \comp^{\B,\dR} (\gamma^{-1}c_{1}\gamma) \big) . $$ 
It is thus also natural to consider the following :

\begin{Definition} Let 
\begin{equation} \label{eq:complex01 bis} \zeta_{\exp}^{\Ad}\big( (n_{i})_{d};(\xi_{i})_{d+1};l;\chi) = \sum\limits_{\xi \in \mu_{N}(\mathbb{C})} \chi(\xi) ({\Phi_{\KZ}^{(\xi)}}^{-1}e^{2i\pi e_{\xi}}\Phi^{(\xi)}_{\KZ}) [e_{0}^{l}e_{\xi_{d+1}}e_{0}^{n_{d}-1}e_{\xi_{d}} \ldots e_{0}^{n_{1}-1}e_{\xi_{1}}] .
\end{equation}
\end{Definition}

We define similarly $\zeta_{\exp}^{\Lambda,\Ad}$. Finally, since unlike in the $p$-adic case there is no privileged choice of $\chi$, it is also natural to consider the following :

\begin{Definition} Let 
	\begin{equation} \label{eq:complex01 bis} \zeta_{\exp,1}^{\Ad}\big( (n_{i})_{d};(\xi_{i})_{d+1};l) = ({\Phi_{\KZ}^{(\xi)}}^{-1}e^{2i\pi e_{\xi}}\Phi^{(\xi)}_{\KZ}) [e_{0}^{l}e_{\xi_{d+1}}e_{0}^{n_{d}-1}e_{\xi_{d}} \ldots e_{0}^{n_{1}-1}e_{\xi_{1}}] .
	\end{equation}
\end{Definition}

We define similarly $\zeta_{\exp,1}^{\Lambda,\Ad}$

We have ${\Phi^{(\xi)}_{\KZ}}^{-1}e^{2i\pi e_{\xi}}\Phi^{(\xi)}_{\KZ} = 1 + 2i\pi {\Phi^{(\xi)}_{\KZ}}^{-1}e_{\xi}\Phi^{(\xi)}_{\KZ} + \zeta(2)
{\Phi_{\KZ}^{(\xi)}}^{-1}\sum\limits_{n\geq 2} \frac{(2i\pi)^{n-2}e_{1}^{n}}{n!} {\Phi_{\KZ}^{(\xi)}}^{-1}$. 
Thus, the numbers (\ref{eq:complex01 bis}) divided by $2 \pi i$ are congruent to the numbers (\ref{eq:complex01}) modulo $\zeta(2)$ and they satisfy the results of \S3 and \S4.

One has a variant of cyclotomic multiple zeta values defined as follows. 
Let $\phi_{\infty}$ be the Frobenius at infinity of $\pi_{1}^{\un,\dR}(\mathbb{P}^{1} - \{0,\mu_{N},\infty\})$ ;
it induces the automorphism of $\Lie(\Pi_{0,0})$ which sends $(e_{0},e_{1}) \mapsto (-e_{0},-\Phi_{-}^{-1}e_{1}\Phi_{-})$, where $\Phi_{-} = \phi_{\infty}({}_{\vec{1}_{1}} 1 _{\vec{1}_{0}}) \in \Pi_{1,0}(\mathbb{R})$. $\Phi_{-}$ is $h_{\KZ}$ in \cite{Enriquez}, \S11.1 whereas the $N=1$ case of $\Phi_{-}$ is $g_{\KZ}$ in \cite{Enriquez}, \S11.1. The coefficients of $\Phi_{-}$ are denoted by 
$\zeta_{-}\big((n_{i})_{d};(\xi_{i})_{d}\big) = \Phi_{-}[e_{0}^{n_{d}-1}e_{\xi_{d}} \ldots e_{0}^{n_{1}-1}e_{\xi_{1}}]$. They generate a $\mathbb{Q}$-vector subspace of the one of MZV$\mu_{N}$'s. (For details in the $N=1$ case see \cite{Furusho 2}, \S2.2).

\begin{Definition} Let the adjoint variants of the numbers $\zeta_{-}$ be :
\begin{equation} \zeta_{-}^{\Ad} \big( (n_{i})_{d};(\xi_{i})_{d+1};b \big) = (-1)^{d}
\sum\limits_{\xi \in \mu_{N}(\mathbb{C})} \xi 
({\Phi^{(\xi)}_{-}}^{-1}e_{\xi}\Phi^{(\xi)}_{-})
\big[ e_{0}^{b} e_{\xi_{d+1}}e_{0}^{n_{d}-1}e_{\xi_{d}} \ldots e_{0}^{n_{1}-1}e_{\xi_{1}} \big] .
\end{equation}
\end{Definition}

In the $N=1$ case, $\Phi_{-}$ is in $\GRT_{1}(\mathbb{R})$ and in $\DS_{0}(\mathbb{R})$, \cite{Furusho 2} (in particular $\Phi_{-}[e_{0}e_{1}]=0$). In the general case, we have a similar fact by \cite{Enriquez}, \S11.1 and \cite{Racinet}.

\subsubsection{Relation with cyclotomic multiple zeta values} 

\begin{Proposition} (i) The respective $k_{N}[2\pi i]$-algebras generated by the numbers $\zeta^{\Ad}(.)$, $\zeta_{\exp}^{\Ad}(.)$, $\zeta_{\exp,1}^{\Ad}(.)$ and $\zeta(.)$ are equal. 
\newline (ii) The respective $k_{N}$-algebras generated by the numbers $\zeta_{-}^{\Ad}(.)$ and $\zeta^{-}(.)$ are equal. 	
\end{Proposition}

\begin{proof} Similar to the proof of Proposition 2.2.4 plus the inversion of a Vandermonde linear system.
	
By inverting a Vandermonde linear system, the data of (\ref{eq:complex01}) resp. (\ref{eq:complex01 bis}) is equivalent respectively to the data of the sequences
$\displaystyle \bigg( ({\Phi_{\KZ}^{(\xi)}}^{-1} e_{\xi}\Phi^{(\xi)}_{\KZ}) [e_{0}^{l}e_{\xi_{d+1}}e_{0}^{n_{d}-1}e_{\xi_{d}} \ldots e_{0}^{n_{1}-1}e_{\xi_{1}}] \bigg)_{\xi \in \mu_{N}(K)}$, and 
\newline 
$\bigg( ({\Phi_{\KZ}^{(\xi)}}^{-1}e^{2i\pi e_{\xi}}\Phi^{(\xi)}_{\KZ}) [e_{0}^{l}e_{\xi_{d+1}}e_{0}^{n_{d}-1}e_{\xi_{d}} \ldots e_{0}^{n_{1}-1}e_{\xi_{1}}]\bigg)_{\xi \in \mu_{N}(K)}$. Moreover, we have $\exp \big( {\Phi_{\KZ}^{(\xi)}}^{-1} e_{\xi}\Phi^{(\xi)}_{\KZ}) \big) = {\Phi_{\KZ}^{(\xi)}}^{-1} \exp(2\pi i e_{\xi}) \Phi^{(\xi)}_{\KZ}$. We also note that in the two above sequences, all the terms can be obtained one from another by applying an automorphism $(z \mapsto \xi z)_{\ast}$ of 
$\pi_{1}^{\un,\dR}(\mathbb{P}^{1} - \{0,\mu_{N},\infty\})$.
\end{proof}

Note that, as in the $p$-adic case, the proof gives more precise information on the fact that the above equalities are compatible with depth filtration. They are also of course compatible with the weight.

In the end, in the complex case there are several natural objects which can be called adjoint variants of cyclotomic multiple zeta values, and they can be considered as equivalent. In function of the property that we are looking at, it can be more convenient to consider a version or another of this object.

\subsection{Cyclotomic multiple harmonic values}

\subsubsection{Definition\label{mhv definition}}

\begin{Definition} \label{def harmonic}For any index $w=((n_{i})_{d};(\xi_{i})_{d+1})$, let
\newline\indent (i) for each $p \in \mathcal{P}_{N}$, $\har_{p^{\mathbb{N}}}(w) 
= \big( \har_{p^{\alpha}} (w)\big)_{\alpha \in \mathbb{N}} \in K_{p}^{\mathbb{N}}$, which we call a $p$-adic cyclotomic multiple harmonic value.
\newline\indent (ii) for each $\alpha \in \mathbb{N}^{\ast}$, $\har_{\mathcal{P}^{\alpha}}(w) = \big( \har_{p^{\alpha}} (w)\big)_{p  \in  \mathcal{P}_{N}} \in \prod_{p\in\mathcal{P}_{N}} K_{p}$, which we call an adelic cyclotomic multiple harmonic value.
\newline\indent (iii) $\har_{\mathcal{P}_{N}^{\mathbb{N}}}(w) =
\big( \har_{p^{\alpha}} (w)\big)_{(p,\alpha)  \in  \mathcal{P}_{N} \times \mathbb{N}} \in \big( \prod_{p\in\mathcal{P}_{N}} K_{p} \big)^{\mathbb{N}}$, which we call a ($p$-adic $\times$ adelic) cyclotomic multiple harmonic value.
\end{Definition}

We will view (i) as the natural explicit $p$-adic substitute to $p$MZV$\mu_{N}$'s ; (ii) as the natural lift of the cyclotomic finite MZV's, defined in \S6 ; (iii) as the natural way to formulate the algebraic properties of (i) and (ii) in a unified way. We will refer most of the time to (iii) and omit the adjective ``$p$-adic $\times$ adelic''. We will abbreviate cyclotomic multiple harmonic values as MHV$\mu_{N}$'s.
\newline\indent Let $\mathcal{O}^{\ast}_{\har}$ be the $\mathbb{Q}$-vector space generated by the empty word and the words of the form $\big((n_{i})_{d};(\xi_{i})_{d+1} \big)$.

\subsubsection{Setting for computations, integrals at (1,0)\label{setting infinite sums}}

In the framework of integrals at $(1,0)$, represented by the symbol $\int_{1,0}$, we will transfer algebraic properties from MZV$\mu_{N}$'s to $\Lambda$-adic AdMZV$\mu_{N}$'s, and the $\Lambda=1$ case will give properties of MHV$\mu_{N}$'s, via equation (\ref{eq:formula for n=1}).
By that equation, the index $\big((n_{i})_{d};(\xi_{i})_{d+1} \big)$ of MHV$\mu_{N}$'s is identified to $\frac{1}{1-e_{0}}e_{\xi_{d+1}}e_{0}^{n_{d}-1}e_{\xi_{d}}\ldots e_{0}^{n_{1}-1}e_{\xi_{1}}$.
The map $\tau(\Lambda) : \mathcal{O}^{\sh} \rightarrow \mathcal{O}^{\sh}[\Lambda] ; w \mapsto \Lambda^{\weight(w)}w$ sends an element to its orbit under the motivic Galois action $\tau$ of $\mathbb{G}_{m} (\ref{eq:tau})$.

\begin{Notation} \label{definition sigma inv DR}
(i) Let $\comp^{\text{har},\Ad} : \mathcal{O}^{\sh} \rightarrow \widehat{\mathcal{O}^{\sh}}$, $w \mapsto \frac{1}{1-e_{0}}w$.
\newline (ii) Let $\comp^{\Lambda \Ad,\Ad} : \mathcal{O}^{\sh} \rightarrow \widehat{\tau(\Lambda)\mathcal{O}^{\sh}}$, $w \mapsto \frac{1}{1-\Lambda e_{0}} w$.
\newline (iii) Let $\comp_{\Lambda}^{\Lambda \Ad,\har} : \widehat{\mathcal{O}^{\sh}} \mapsto  \widehat{\tau(\Lambda)\mathcal{O}^{\sh}}$, $\frac{1}{1-e_{0}}w \mapsto \frac{1}{1-\Lambda e_{0}}w$.
\end{Notation}

We also denote in the same way the duals of these maps. The map $\comp^{\text{har},\Ad}$ was called $\comp^{\Sigma \smallint}$, ``comparison from integrals to series'', in \cite{I-2}. In this paper, it seems more relevant to call it the ``comparison from adjoint (analogues of MZV$\mu_{N}$'s) to harmonic (analogues of MZV$\mu_{N}$'s)''. We note that $\comp^{\Lambda \Ad, \Ad} =  \widehat{\tau(\Lambda)}^{-1} \circ \big(\frac{1}{\Lambda}\comp^{\har,\Ad}\big) \circ \widehat{\tau(\Lambda)}$.
\newline\indent The strategy of transferring algebraic properties stated above is motivated by the notions of adjoint Ihara action $\circ_{\Ad}^{\smallint_{1,0}}$ (\cite{I-2} Definition 1.1.3) and pro-unipotent harmonic action $\circ_{\har}^{\smallint_{1,0}}$ (\cite{I-3} Definition 2.1.2), which are pushforwards of the usual Ihara action (\ref{eq:Ihara}) by $\Ad(e_{1})$ and $\comp^{\text{har},\Ad}$ which appears implicitly in equation (\ref{eq:formula for n=1}) ; and with the fact that the Ihara action is compatible with the algebraic relations of MZV$\mu_{N}$'s.
\newline\indent $\Ad(e_{1})$ restricted to $ \tilde{\Pi}_{1,0}(K) = \{f \in \Pi_{1,0}(K) \text{ }|\text{ }f[e_{0}] = f[e_{1}] = 0\}$ is injective. The map $\comp^{\Sigma \smallint}$, viewed here as  $\comp^{\text{har},\Ad}$, is not injective, but taking into account all the relations of iteration of the Frobenius enables to replace it by an injective map (\cite{I-3}, \S4.2).
Thus we consider that applying the push-forward by $\Ad(e_{1})$ and $\comp^{\text{har},\Ad}$ should not lose information.

\subsubsection{Setting for computations, power series expansions of integrals at $0$}

In that framework, we will view prime weighted cyclotomic multiple harmonic sums via equation (\ref{eq:Li 0}). Let us write more precisely how they are connected to multiple polylogarithms.

\begin{Notation} \label{notation coefficient of power series}For $S = \sum\limits_{(n_{1},\ldots,n_{r}) \in \mathbb{N}^{r}} c_{n_{1},\ldots,n_{r}} x_{1}^{n_{1}} \ldots x_{r}^{n_{r}} \in R[[x_{1},\ldots,x_{r}]]$ a formal power series, where $R$ is a ring, for all $(n_{1},\ldots,n_{r}) \in \mathbb{N}^{r}$, we write $c_{n_{1},\ldots,n_{r}}= S[x_{1}^{n_{1}} \ldots x_{r}^{n_{r}}]$.
\end{Notation}

By the power series expansion of multiple polylogarithms (\ref{eq:Li series bis}), we have two slightly different formulas :
\begin{equation} \har_{m}\big( (n_{i})_{d};(\xi_{i})_{d+1} \big) =  m^{n_{1}+\ldots+n_{d}+l}\Li[e_{0}^{l-1}e_{\xi_{d+1}}e_{0}^{n_{d}-1}e_{\xi_{d}} \ldots e_{0}^{n_{1}-1}e_{\xi_{1}}][z^{m}];
\end{equation}
\begin{equation} \label{eq:har et Li 2} \har_{m}((n_{i})_{d};(\xi_{i})_{d+1}) = m^{n_{1}+\ldots+n_{d}} \sum_{0<m'<m} \xi_{d+1}^{m-m'}  \Li[ e_{\xi_{d+1}}e_{0}^{n_{d}-1}e_{\xi_{d}} \ldots e_{0}^{n_{1}-1}e_{\xi_{1}}][z^{m'}].
\end{equation}
\indent We consider multiple polylogarithms restricted to a function on a $(\mathbb{P}^{1} - Z)^{2} - \text{diagonal}$ where $Z$ is a finite subset of $\mathbb{P}^{1}$, and we choose $O = (\vec{1}_{0},\vec{1}_{0})$ in cubic coordinates as the origin of paths of integration. 
\newline\indent We use the notation $\omega_{z}=\frac{dy}{y-z}$ for $z \in Z- \{0,\infty\}$. By (\ref{eq:Li series bis}), for $(y_{1},y_{2})$ on a neighbourhood of the chosen reference base-point $O$, for $\gamma$ the straight path from $O$ to $(y_{1},y_{2})$, for any two words
$w = \big( (n_{i})_{d}:(z_{i})_{d} \big)$, 
$\tilde{w} = \big( (\tilde{n}_{i})_{d}:(\tilde{z}_{i})_{d} \big)$, we have :
\begin{multline}
\int_{\gamma} \big( \omega_{0}^{t_{d'}-1}\omega_{z_{d'}} \ldots \omega_{0}^{t_{1}-1}\omega_{z_{1}}\big)(y_{2}) \big(\omega_{0}^{n_{d}-1}\omega_{z_{d}} \ldots \omega_{0}^{n_{1}-1}\omega_{z_{1}}\big)(y_{1}y_{2})
\\ =
\sum_{0<m_{1}<\ldots<m_{d}<m'_{1}<\ldots<m'_{d'}=m}
\frac{\big( \frac{z_{2}}{z_{1}} \big)^{m_{1}} \ldots \big(\frac{y_{1}}{z_{d}}\big)^{m} 
\big( \frac{\tilde{z}_{2}}{\tilde{z}_{1}} \big)^{m_{1}} \ldots \big(\frac{y_{2}}{\tilde{z}_{d'}}\big)^{m}} {m_{1}^{n_{1}}\ldots m_{d}^{n_{d}} m_{1}^{\tilde{n}_{1}}\ldots m_{d'-1}^{\tilde{n}_{d'-1}} m_{d}^{\tilde{n}_{d'}}}.
\end{multline}

\subsubsection{Setting for computations, series}

In this framework, we will view prime weighted multiple harmonic sums via equation (\ref{eq:mult har sums}). We consider them as functions of the upper bound $m$ of their domain of summation, $\mathbb{N} \rightarrow \mathbb{C}$, resp. $\mathbb{N} \rightarrow K_{p}$. We will study them by using only the structure of topological field of $\mathbb{C}$ and $\mathbb{C}_{p}$, both at the source and at the target.
\newline\indent If $c : \mathbb{N} \rightarrow \mathbb{C}$ is a function such that $c(0)=0$, the Newton series of $c$ is the function : $\displaystyle z \mapsto \sum\limits_{m \in\mathbb{N}} (\nabla c)_{m} {z \choose m}$, where 
$\nabla c = \Big(\sum\limits_{m'=0}^{m} (-1)^{m'}{m \choose m'} c(m')\Big)_{m \in \mathbb{N}}$ and $\displaystyle {z \choose m} = \frac{z(z-1)\ldots(z-m+1)}{m!}$. It is usually defined on a half-plane of the type $\{ \text{Re}(z) > \rho \}$ of $\mathbb{C}$. This notion appears in certain proofs of results on multiple harmonic sums by several other authors, which we will generalize and interpret in terms of cyclotomic multiple harmonic values.

\subsection{``Overconvergent'' variants of cyclotomic multiple harmonic values}

\subsubsection{Definition}

For any word $w$ on $e_{0 \cup \mu_{N}}$, the power series expansion of the functions $\Li_{p,X_{K}}^{\KZ}[w]$ at $0$ in $k_{N}[[z]][\log(z)]$ in is identical to the one of the complex multiple polylogarithm $\Li[w]$ of (\ref{eq:multiple polylogarithms power series expansion}). Thus, the coefficients of the power series expansions of $\Li_{p,\alpha}^{\dagger}$ can be written in terms of multiple harmonic sums and cyclotomic $p$-adic multiple zeta values. In the next statement, we use Notation \ref{notation coefficient of power series}.

\begin{Definition} 
	\noindent (i) Let the overconvergent variants of cyclotomic multiple harmonic sums be the numbers $\frak{h}_{m}^{\dagger_{p,\alpha}}(w) = \Li_{p,\alpha}^{\dagger}[w][z^{m}] \in K$, where $w$ is a word on $e_{0 \cup \mu_{N}}$ and $m \in \mathbb{N}$.
	\newline (ii) Let the overconvergent variants of weighted cyclotomic multiple harmonic sums be the numbers $\har_{m}^{\dagger_{p,\alpha}}(w) = m^{\weight(w)}\Li_{p,\alpha}^{\dagger}[w][z^{m}]$, where $w$ is a word on $e_{0 \cup \mu_{N}}$ and $m \in \mathbb{N}$.
	\newline (iii) Let overconvergent prime weighted cyclotomic multiple harmonic sums be the numbers
	$\har_{p^{\alpha}}^{\dagger_{p,\alpha}}(w)$, where $w$ is a word on $e_{0 \cup \mu_{N}}$.
\end{Definition}

The following is a variant of Definition \ref{def harmonic}, from which we take the same notations. We now consider all $p$'s resp. all $\alpha$'s at the same time.

\begin{Definition} \label{variant harmonic}For any index $w= \big( (n_{i})_{d};(\xi_{i})_{d+1} \big)$, we call respectively
	\newline (i) for each $p \in \mathcal{P}_{N}$, $\har^{\dagger}_{p^{\mathbb{N}}}(w) 
	= \big( \har^{\dagger_{p,\alpha}}_{p^{\alpha}} (w)\big)_{\alpha \in \mathbb{N}} \in K_{p}^{\mathbb{N}}$, overconvergent variants of adelic cyclotomic multiple harmonic values.
	\newline (ii) for each $\alpha \in \mathbb{N}^{\ast}$, $\har^{\dagger}_{\mathcal{P}^{\alpha}}(w) = \big( \har^{\dagger_{p,\alpha}}_{p^{\alpha}} (w)\big)_{p  \in \mathcal{P}_{N}} \in \underset{p\in\mathcal{P}_{N}}{\prod} K_{p}$, overconvergent variants of $p$-adic cyclotomic multiple harmonic values.
	\newline (iii) $\har^{\dagger}_{\mathcal{P}_{N}^{\mathbb{N}}}(w)=\big( \har^{\dagger_{p,\alpha}}_{p^{\alpha}} (w)\big)_{(p,\alpha)  \in  \mathcal{P}_{N} \times \mathbb{N}} \in \big( \underset{p\in\mathcal{P}_{N}}{\prod} K_{p} \big)^{\mathbb{N}}$, overconvergent variants of ($p$-adic $\times$ adelic) cyclotomic multiple harmonic values.
\end{Definition}

As in the previous sections, we will usually refer to (iii) and omit the adjective ``$p$-adic $\times$ adelic''. We will use the abbreviation MHV$\mu_{N}^{\dagger}$.

\subsubsection{Setting for computations}

In the previous parts, we have used three settings for computations (\S2.2.2, \S2.2.3, \S2.2.4) for studying prime weighted multiple harmonic values, corresponding to the frameworks $\int_{1,0}$, $\int$ and $\Sigma$ respectively. We now add to them a complement :

\begin{Proposition} \label{prop formula dagger} (i) ($\int_{1,0}$) We have : 
	\begin{equation} \label{eq:overconv expression}
	\har_{p^{\alpha}}^{\dagger_{p,\alpha}}[e_{0}^{l}e_{\xi_{d+1}}e_{0}^{n_{d}-1}e_{\xi_{d}}\ldots e_{0}^{n_{1}-1}e_{\xi_{1}}]
	= \sum_{b=l}^{\infty} \sum_{\xi \in \mu_{N}(K)} \xi^{-p^{\alpha}}
	\Ad_{{\Phi_{p,\alpha}^{(\xi)}}}(e_{\xi}) [e_{0}^{b}e_{\xi_{d+1}}e_{0}^{n_{d}-1}e_{\xi_{d}}\ldots e_{0}^{n_{1}-1}e_{\xi_{1}}] .
	\end{equation}
	\noindent (ii) ($\int$) \label{lemma power series expansion Li dagger} We have $ \Li_{p,\alpha}^{\dagger}[w][z^{p^{\alpha}}] = \har^{\dagger}_{p,\alpha}$ and, for $m \in \{1,\ldots,p^{\alpha}-1\}$,  $\Li_{p,\alpha}^{\dagger}[w][z^{m}] = (p^{\alpha})^{\weight(w)}\Li(z)[z^{m}] = \frac{1}{m^{l}}\frak{h}_{m} \big((n_{i})_{d};(\xi_{i})_{d+1} \big)$.
	\newline\noindent (iii) ($\Sigma$) We can write a formula as sums of series for each $\har_{p^{\alpha}}^{\dagger_{p,\alpha}}(w)$ by composing (\ref{eq:overconv expression}) with the formula for $\Ad_{\Phi_{p,\alpha}^{(\xi)}}(e_{\xi})$ as sums of series obtained in the main theorem of \cite{I-2}.
\end{Proposition}

\begin{proof} (i) For any word $w=e_{0}^{l}e_{\xi_{d+1}}e_{0}^{n_{d}-1}e_{\xi_{d}}\ldots e_{0}^{n_{1}-1}e_{\xi_{1}}$ and let us consider the coefficient $[w][z^{p^{\alpha}}]$ of equation (\ref{eq:horizontality1}) ; we obtain : 
	\begin{multline} \label{eq:remainder} \har_{p^{\alpha}}^{\dagger_{p,\alpha}}
	(e_{0}^{l}e_{\xi_{d+1}}e_{0}^{n_{d}-1}e_{\xi_{d}}\ldots e_{0}^{n_{1}-1}e_{\xi_{1}})
	\\ = \har_{p^{\alpha}} \big( (n_{i})_{d};(\xi_{i})_{d+1} \big) -
	\sum_{b=0}^{l} \sum_{\xi \in \mu_{N}(K)} \xi^{-1}
	\Ad_{{\Phi_{p,\alpha}^{(\xi)}}}(e_{\xi}) [e_{0}^{b}e_{\xi_{d+1}}e_{0}^{n_{d}-1}e_{\xi_{d}}\ldots e_{0}^{n_{1}-1}e_{\xi_{1}}] .
	\end{multline}
	If we combine this and equation (\ref{eq:formula for n=1}) proven in \cite{I-2}, we obtain the result. (In \cite{I-2}, we actually obtained equation (\ref{eq:formula for n=1}) by proving that $\har_{p^{\alpha}}^{\dagger_{p,\alpha}}(e_{0}^{l}e_{z_{i_{d+1}}}e_{0}^{n_{d}-1}e_{\xi_{d}}\ldots e_{0}^{n_{1}-1}e_{\xi_{1}}) \displaystyle \underset{l \rightarrow \infty}{\rightarrow} 0$, which was an immediate consequence of the main result of \cite{I-1}.)
	\newline (ii) Follows from equation (\ref{eq:horizontality1}) and from that we have $\Li(z^{p^{\alpha}})[z^{m}]=0$ for all $m \in \{1,\ldots,p^{\alpha}-1\}$. 
	\newline (iii) Immediate.
\end{proof}

In particular, (\ref{eq:overconv expression}) means that the numbers $\har_{p^{\alpha}}^{\dagger_{p,\alpha}}(w)$ are the remainders of the sums of series of (\ref{eq:formula for n=1}) which express $\har_{p^{\alpha}}$ in terms of cyclotomic $p$-adic multiple zeta values. This leads us to a variant of Definition \ref{def adjoint} (ii) :

\begin{Definition} \label{def over adjoint}
	Let, for all words, the overconvergent variants of $\Lambda$-adjoint $p$-adic cyclotomic multiple zeta values  ($\Lambda$Ad$p$MZV$\mu_{N}^{\dagger}$'s) be 
	\begin{multline} \zeta^{\Lambda,\Ad,\dagger}_{p^{\alpha}} (e_{0}^{l}e_{\xi_{d+1}}e_{0}^{n_{d}-1}e_{\xi_{d}}\ldots e_{0}^{n_{1}-1}e_{\xi_{i_{1}}})
	\\ = \sum_{b \geqslant l} \sum_{\xi \in \mu_{N}(K)} \xi^{-p^{\alpha}} \Lambda^{b+n_{d}+\ldots+n_{1}}
	\Ad_{\Phi_{p,\alpha}^{(\xi)}}(e_{\xi})
	[e_{0}^{b}e_{\xi_{d+1}}e_{0}^{n_{d}-1}e_{\xi_{d}}\ldots e_{0}^{n_{1}-1}e_{\xi_{1}}] .
	\end{multline}
\end{Definition}

In particular, the $l=0$ case of $\Lambda$Ad$p$MZV$\mu_{N}^{\dagger}$'s are just $\Lambda$Ad$p$MZV$\mu_{N}$'s. For $l\geqslant 0$, they are remainders of $\Lambda$Ad$p$MZV$\mu_{N}$'s viewed as power series in $\Lambda$.

\subsection{Adjoint multiple polylogarithms and harmonic multiple polylogarithms}

Adjoint cyclotomic multiple zeta values and cyclotomic multiple harmonic values are particular values of functions which retain a lot of their properties, and which we introduce now.

\subsubsection{Adjoint multiple polylogarithms}

In this context, the generalization of cyclotomic multiple zeta values are multiple polylogarithms, with the only restriction that we assume them to be iterated integrals on the straight path from $0$ to $1$.
\newline\indent This suggests the following generalization of Definition \ref{def adjoint} :

\begin{Definition} Let adjoint multiple polylogarithms be 
	\begin{multline*} \Li^{\Ad}\big( (n_{i})_{d},l;(z_{i})_{d+1} \big) = (-1)^{d} \sum_{z \in Z} z^{-1} ({\Li^{(z)}}^{-1}e_{z}\Li^{(z)})\big[e_{0}^{l}e_{z_{d+1}}e_{0}^{n_{d}-1}e_{z_{d}}\ldots e_{0}^{n_{1}-1}e_{z_{1}} \big] \\ = \sum_{d'=1}^{d+1} \prod_{i=d'}^{d} {-n_{i} \choose l_{i}}
	z_{d'}^{-1} (-1)^{n_{d'}+\ldots+n_{d}}
	\Li^{(z_{d'})} \big( (n_{d-i}+l_{d-i})_{d-d'},(z_{d-i})_{d-d'} \big)
	\Li^{(z_{d'})} \big( (n_{i})_{d'-1},(z_{i})_{d'-1} \big) 
	\end{multline*}
	where $\Li^{(z)}$ denotes an iterated integral on a straight path from $0$ to $z$.
\end{Definition}

Of course, as in \S7.1 we also have a variant where $e_{z}$ above is replaced by $e^{2i\pi e_{z}}$. We note that we have the following adjoint variant of the KZ equation (\ref{eq: nabla KZ}) : for all $u \in \Lie(\Pi_{0,0})$,
$$ d\Ad_{\Li}(u) = 
\Ad_{\Li_{p,X_{K}}^{\KZ}} \bigg( \ad_{u} \Big( e_{0} \omega_{0} + \sum_{z_{0} \in D- \{0,\infty\}} e_{z_{0}}\omega_{z_{0}} \Big) \bigg) . $$
Moreover, in the case of $\mathbb{P}^{1} - \{0,\mu_{N},\infty\}$, we have an adjoint variant of the differential equation satisfied by the overconvergent $p$-adic multiple polylogarithms (\ref{eq:horizontality equation}) : for all $u \in \Lie(\Pi_{0,0})$,
\begin{multline*}
d\Ad_{\Li_{p,\alpha}^{\dagger}}(u) = 
\\
\Ad_{\Li_{p,\alpha}^{\dagger}} \bigg( \ad_{u} \Big( \omega_{0}(z^{p^{\alpha}})e_{0} +  \sum_{\xi \in \mu_{N}(K)}  \omega_{\xi^{p^{\alpha}}}(z^{p^{\alpha}})e_{\xi^{p^{\alpha}}}  \Big) \bigg) -
\ad_{\Ad_{\Li_{p,\alpha}^{\dagger}}(u)} \Big(
\omega_{0}(z) e_{0} + \sum_{\xi \in \mu_{N}(K)} \omega_{\xi}(z)\Ad_{\Phi^{(\xi^{p^{\alpha}})}_{p,\alpha}}(e_{\xi^{p^{\alpha}}}) \Big) .
\end{multline*}

\subsubsection{A generalization of cyclotomic multiple harmonic values and finite cyclotomic multiple zeta values}

We now assume that $\mathbb{P}^{1} - D$ is defined over a number field, which is embedded in $\mathbb{C}_{p}$ for all primes $p$.

\begin{Definition} Let multiple harmonic polylogarithms be 
	$$ \har_{\mathcal{P}^{\mathbb{N}}} = 
	\bigg(  \sum_{0<m_{1} <\ldots < m_{d}<p^{\alpha}}
	\frac{\big( \frac{z_{2}}{z_{1}} \big)^{m_{1}} \ldots \big(\frac{z_{d+1}}{z_{d}}\big)^{m_{d}}\big(\frac{1}{z_{d+1}}\big)^{p^{\alpha}}}{m_{1}^{n_{1}}\ldots m_{d}^{n_{d}}}  \bigg) \in \prod_{p \in \mathcal{P}} \mathbb{C}_{p}^{\mathbb{N}^{\ast}} . $$	
\end{Definition}

Unlike in \S2.4 and \S2.5, we do not claim a crystalline meaning for this notion.

\section{Around double shuffle equations\label{double shuffle}}

We review the regularized double shuffle equations (\S3.1), we construct adjoint double shuffle equations (\S3.2) and harmonic double shuffle relations (\S3.3), in the three frameworks $\int_{1,0}$, $\int$, $\Sigma$, which proves the theorem stated in \S1.3. We discuss a particular consequence of the double shuffle equation, which we call the ``reversal'' equation, and we find its adjoint and harmonic counterparts (\S3.4).

\subsection{Review on the double shuffle equations}

Let $K$ be a field of characteristic zero which contains a primitive $N$-th root of unity. Let the alphabet $Y_{N} =\{ y_{n}^{(\xi)}\text{ }|\text{ }n\geqslant 1, \xi \in \mu_{N}(K)\}$. A word on $e_{0\cup \mu_{N}}$ of the form $e_{0}^{n_{d}-1}e_{\xi_{d}} \ldots e_{0}^{n_{1}-1}e_{\xi_{1}}$ can be viewed as the word $y_{n_{1}}^{(\xi_{1})} \ldots y_{n_{d}}^{(\xi_{d})}$ on $Y_{N}$.
\newline\indent Racinet has defined \cite{Racinet}, for each $\lambda\in K$ the set $\DS_{\lambda}(K)$ : this is the set of couples $(\psi_{\sh},\psi_{\ast}) \in K \langle\langle e_{0\cup \mu_{N}} \rangle\rangle \times K \langle\langle Y_{N} \rangle\rangle$ such that $\psi_{\sh}$ satisfies the shuffle equation (equation (\ref{eq:shuffle eq}) below), $\psi_{\ast}$ satisfies the quasi-shuffle equation (equation (\ref{eq:stuffle eq}) below) ; a certain relation between  $\psi_{\sh}$ and $\psi_{\ast}$ (equation (\ref{eq:rel between reg}) below) ; as well as $\psi_{\sh}[e_{0}] = \psi_{\sh}[e_{1}]=0$, and $\psi_{\sh}[e_{0}e_{1}]=-\lambda^{2}/24$. This defines $\DS_{\lambda}$ as a subscheme of $\Pi_{1,0}$. Racinet proves that the Ihara product (\ref{eq:Ihara}) restricts to a group law on $\DS_{0}$ and to an action of $\DS_{0}$ on $\DS_{\lambda}$ which makes it a torsor.
\newline\indent We have $\Phi_{\KZ} \in \DS_{2\pi i}(\mathbb{C})$, and $\Phi_{p,\alpha} \in \DS_{0}(K_{p})$ : for $\alpha=1$ and $\alpha=-\infty$, this follows from \cite{Besser Furusho} ($N=1$), \cite{Yamashita} (any $N$); by the relations of iteration of the Frobenius (\cite{I-3}, equations (1.11), (1.12), (1.13) and Proposition 1.5.2), and by Racinet's theorem, it follows that this is true for any $\alpha \in \mathbb{Z} \cup \{\pm \infty\} - \{0\}$.

\subsubsection{The shuffle relation of iterated integrals}

The shuffle equation appears in equation (\ref{eq:shuffle equation}) : an element $f \in k \langle \langle e_{0 \cup \mu_{N}} \rangle\rangle$ satisfies the shuffle equation if :
\begin{equation} \label{eq:shuffle eq} \forall w,w' \text{ words on }e_{0\cup \mu_{N}}, \text{ we have } f[w]f[w'] = f[w\text{ }\sh\text{ }w'] 
\end{equation}
where $\sh$ is the shuffle product reviewed in \S2.1.1, characterized by induction on the weight by $w\text{ }\sh\text{ }1 = 1\text{ }\sh\text{ }w = w$ ($1$ is the empty word) and $e_{x}w \text{ } \sh \text{ } e_{x'}w'  
= e_{x}(w \text{ }\sh\text{ } e_{x'}w') + e_{x'}(e_{x}w \text{ } \sh\text{ } w')$ or equivalently
$w e_{x} \text{ } \sh\text{ }  w'e_{x}
= (w\text{ } \sh\text{ } w'e_{x'})e_{x} + (w\text{ } \sh\text{ } w'e_{x})e_{x'}$ for all $x,x' \in \{0\} \cup \mu_{N}(K)$. The shuffle equation amounts to 
\begin{equation} \label{eq:shuffle eq bis} \hat{\Delta}_{\sh}(f) = f \hat{\otimes} f 
\end{equation}
where $\hat{\Delta}_{\sh}$ is the coproduct in $K \langle \langle e_{0 \cup \mu_{N}} \rangle\rangle$ viewed as the completed dual of the Hopf algebra $\mathcal{O}^{\sh}$.
\newline\indent The fact that the generating series $\Phi_{\KZ}$, resp. $\Phi_{p,\alpha}$ satisfies the shuffle equation follows directly from their definition as points of $\Pi_{1,0} \simeq \Spec(\mathcal{O}^{\sh,e_{0\cup \mu_{N}}})$ and from equation (\ref{eq:shuffle equation}), and amounts to a family of relations on MZV$\mu_{N}$'s, resp. $p$MZV$\mu_{N}$'s. The coefficient of $\Phi_{\KZ}$ at a word $e_{0}^{n_{d}-1}e_{\xi_{d}} \ldots e_{0}^{n_{1}-1}e_{\xi_{1}}e_{0}^{n_{0}-1}=e_{x_{n}}\ldots e_{x_{1}}$ which does not necessarily satisfy the hypothesis $n_{d}\geqslant 2$ and $n_{0}=1$, is the regularized iterated integral
$\displaystyle \underset{\epsilon,\epsilon' \rightarrow 0}{\Reg\lim} \int_{\epsilon}^{1-\epsilon'} \frac{dt_{n}}{t_{n}-x_{n}} \ldots \frac{dt_{1}}{t_{1}-x_{1}}$
defined as the constant term in the asymptotic expansion of $\int_{\epsilon}^{1-\epsilon'} \frac{dt_{n}}{t_{n}-x_{n}} \ldots \frac{dt_{1}}{t_{1}-x_{1}}$, which is in $\mathbb{C}[[\epsilon,\epsilon']][\log(\epsilon),\log(\epsilon')]$. That coefficient is denoted by $\zeta_{\sh}(e_{x_{n}}\ldots e_{x_{1}})$ and is called a (shuffle-)regularized MZV$\mu_{N}$. Similarly, the coefficient of $\Phi_{p,\alpha}$ at such a word is denoted by $(\zeta_{p,\alpha})_{\sh}(e_{x_{n}}\ldots e_{x_{1}})$ and called a (shuffle-)regularized $p$MZV$\mu_{N}$. The shuffle equation for $\Phi_{\KZ}$ also follows from the identity
\begin{equation} \label{eq:eq shuffle source} \displaystyle \int_{\epsilon<t_{1}< \ldots < t_{n}<1-\epsilon'} \times \int_{\epsilon<t_{n+1}<\ldots t_{n+n'}<1-\epsilon'} = \sum_{\substack{\sigma \text{ permutation of }\{1,\ldots,l+l'\} \\\text{s.t. } \sigma(1)<\ldots<\sigma(n)\\ \text{and } \sigma(n+1)<\ldots<\sigma(n+n')}} \int_{\epsilon<t_{\sigma^{-1}(1)} < \ldots<t_{\sigma^{-1}(n+n')}<1-\epsilon'},
\end{equation}
which follows from an equality between domains of integration. 
The definition of multiple polylogarithms as coefficients of a point of $\pi_{1}^{\un,\dR}(\mathbb{P}^{1} - \{0,\mu_{N},\infty\})$ and equation (\ref{eq:shuffle equation}) imply that they satisfy the shuffle equation ; this also follows from their definition as iterated integrals (\ref{eq:MPL}) and equation (\ref{eq:eq shuffle source}).

\subsubsection{The quasi-shuffle relation of iterated series}

An element $f \in K \langle\langle Y_{N} \rangle\rangle$ satisfies the quasi-shuffle equation if :
\begin{equation} \label{eq:stuffle eq} 
\forall w,w' \text{ words on }Y_{N}, \text{ we have }
f[w]f[w'] = f[w \ast w']
\end{equation}
where $\ast$, the quasi-shuffle product on $\mathbb{Q} \langle Y_{N} \rangle$, is defined as follows, by induction on the depth : $1 \ast w = w \ast 1 = w$ ($1$ is the empty word) and
$w y^{(\xi)}_{n} \ast w'y^{(\xi')}_{n'} = (w \ast w'y^{(\xi')}_{n'})y^{(\xi)}_{n} + (wy^{(\xi)}_{n} \ast w') y^{(\xi')}_{n'} + (w \ast w')y^{(\xi\xi')}_{n+n'}$.
The quasi-shuffle product makes $K \langle\langle Y_{N} \rangle\rangle$ into a commutative algebra. The quasi-shuffle equation amounts to 
\begin{equation} \label{eq:stuffle eq bis}
\hat{\Delta}_{\ast}(f) = f \hat{\otimes} f 
\end{equation}
where $\hat{\Delta}_{\ast}$ is the coproduct in $K \langle \langle Y_{N} \rangle\rangle$ viewed as the completed dual coalgebra of $\mathbb{Q} \langle Y_{N} \rangle$. Actually, one has actually a natural structure of Hopf algebra $\mathcal{O}^{\ast,e_{0\cup \mu_{N}}}$ on $\mathbb{Q} \langle Y_{N} \rangle$ in which the product is $\ast$ (\cite{Hoffman} for $N=1$).
\newline\indent For a word $w=\big((n_{i})_{d};(\xi_{i})_{d}\big)$ such that we do not necessarily have the hypothesis $(n_{d},\xi_{d})\not=(1,1)$ of (\ref{eq:multizetas}), let $\Reg \displaystyle \sum\limits_{0<m_{1}<\ldots<m_{d}} \frac{\big( \frac{\xi_{2}}{\xi_{1}} \big)^{m_{1}} \ldots \big(\frac{1}{\xi_{d}}\big)^{m_{d}}}{m_{1}^{n_{1}}\ldots m_{d}^{n_{d}}}$ be the constant term in the asymptotic expansion of $\displaystyle \sum\limits_{0<m_{1}<\ldots<m_{d}<m} \frac{\big( \frac{\xi_{2}}{\xi_{1}} \big)^{m_{1}} \ldots \big(\frac{1}{\xi_{d}}\big)^{m_{d}}}{m_{1}^{n_{1}}\ldots m_{d}^{n_{d}}}$ when $m \rightarrow \infty$, which is in $\mathbb{C}[\log(m)][[\frac{1}{m}]]$, in which the terms involving the Euler-Mascheroni constant are ``removed'' in a canonical way (see \cite{Racinet}). We denote it by $\zeta_{\ast}(w)$ and we call it a regularized MZV$\mu_{N}$. Let $(\Phi_{\KZ})_{\ast}= 1+ \sum\limits_{w\text{ word on }Y_{N}} \zeta_{\ast}(w)w$ ; $(\Phi_{\KZ})_{\ast}$ satisfies the quasi-shuffle equation ; this follows from the identities
\begin{equation} \label{eq:quasi shuffle eq source} \sum_{0<m_{1}<\ldots<m_{d}<m} \times \sum_{0<m'_{1}<\ldots<m'_{d'}<m} = \sum_{\text{quasi-shuffle elements}}\text{ }\sum_{0<m''_{1}<\ldots<m''_{r}<m}
\end{equation}
where a quasi-shuffle element is a way to order $m_{1},\ldots,m_{d}$ and $m'_{1},\ldots,m'_{d'}$, which determines an integer $r \in \{\max(d,d'),\ldots,d+d'\}$ and variables $m''_{1},\ldots,m''_{r}$ : for example, for $d=d'=1$, there are three quasi-shuffle elements : $\{m_{1}<m'_{1}\}$, $\{m_{1}=m'_{1}\}$ and $\{m_{1}>m'_{1}\}$.
\newline\indent In the $p$-adic case, the fact that $p$MZV$\mu_{N}$'s satisfy the quasi-shuffle equations \cite{Besser Furusho} is proved by using formal properties of Coleman functions. Moreover, an analogue of the quasi-shuffle regularization which still satisfies the quasi-shuffle equation has been defined in \cite{Furusho Jafari} in the case of $N=1$. We can deduce easily a similar definition for any $N$.
\newline\indent One has a variant $\ast_{\har}$ on $\mathcal{O}_{\har}^{\ast}$ (defined in \S2.2.1) adapted to cyclotomic multiple harmonic sums. For any words $w = \big((n_{i})_{d};(\xi_{i})_{d+1})$ and
$w' = \big((\tilde{n}_{i})_{d};(\tilde{\xi}_{i})_{d+1})$, $w \ast_{\har} w'$ is the sum, indexed by the set of quasi-shuffle elements $(u_{i})_{d''}$ of 
$\big((n_{i})_{d},(\tilde{n}_{i})_{d}\big)$, of the sequences 
$\big((u_{i})_{d''}, (\xi_{a_{i}} \tilde{\xi}_{b_{i}})_{d''+1} \big)$ defined as follows : $a_{1} = b_{1} = 1$ and, 
for $2 \leqslant i \leqslant d''$,
$ (a_{i},b_{i}) = \left\{ \begin{array}{ll} (a_{i-1}+1 ,b_{i-1})&
\text{ if } \exists l\text{ }|\text{ } u_{i-1}=n_{l},
\\ (a_{i-1}, b_{i-1}+1) &\text{ if }\exists l'\text{ }|\text{ } u_{i-1}=\tilde{n}_{l'},\text{ } 
\\ (a_{i-1}+1,b_{i-1}+1)&\text{ if }\exists l,l' \text{ }|\text{ } u_{i-1}=n_{l}+\tilde{n}_{l'},\text{ } \end{array}\right.$. This makes $\mathcal{O}_{\har}^{\ast}$ into a a commutative algebra. We have natural morphisms of algebras $i : \mathcal{O}^{\ast} \hookrightarrow \mathcal{O}_{\har}^{\ast}$, $\big( (n_{i})_{d};(\xi_{i})_{d} \big) \mapsto \big( (n_{i})_{d};((\xi_{i})_{d},1) \big)$ and $r : \mathcal{O}_{\har}^{\ast} \twoheadrightarrow \mathcal{O}^{\ast}$,
$\big( (n_{i})_{d};(\xi_{i})_{d+1} \big) \mapsto \big( (n_{i})_{d}; (\frac{\xi_{i}}{\xi_{d+1}})\big) $, with $r \circ i=\id$.
\newline\indent By equation (\ref{eq:quasi shuffle eq source}), the power series expansion of multiple polylogarithms (\ref{eq:multiple polylogarithms power series expansion}) satisfy a version of the quasi-shuffle relation. It can be encoded by means of $\ast_{\har}$ ; we will use it in \S3.3.2.

\subsubsection{The relation between the two regularizations of MZV$\mu_{N}$'s}

The two regularizations of MZV$\mu_{N}$'s are related as follows. Let $\pr : K\langle \langle e_{0\cup \mu_{N}} \rangle\rangle \rightarrow K\langle \langle Y_{N} \rangle\rangle$ be the unique continuous (for the weight-adic topology) and linear map which sends a word $w$ to itself if its rightmost letter is not $e_{0}$ and to $0$ if it is $e_{0}$, where we view words on $Y_{N}$ as words on $e_{0\cup \mu_{N}}$ as usual.
Let the maps $\textbf{p},\textbf{q} : \mathcal{O}^{\sh} \rightarrow \mathcal{O}^{\sh}$ defined as follows  (\cite{Racinet}, \S2.2.3) :
$$ \textbf{p}(e_{0}^{n_{d}-1}e_{\sigma_{d}} \cdots e_{0}^{n_{1}-1}e_{\sigma_{1}}e_{0}^{n_{0}-1} ) = e_{0}^{n_{d}-1}e_{\sigma_{d}^{-1}} \cdots
e_{0}^{n_{2}-1}e_{(\sigma_{d}\cdots \sigma_{2})^{-1}} e_{0}^{n_{1}-1}e_{(\sigma_{d}\cdots \sigma_{1})^{-1}}e_{0}^{n_{0}-1} $$
$$ \textbf{q} (e_{0}^{n_{d}-1}e_{\xi_{d}} \cdots e_{0}^{n_{1}-1}e_{\xi_{1}}e_{0}^{n_{0}-1} ) = e_{0}^{n_{d}-1}e_{\xi_{d}^{-1}} \cdots e_{0}^{n_{1}-1}e_{\xi_{3}^{-1}\xi_{2}} e_{0}^{n_{1}-1}e_{\xi_{2}^{-1}\xi_{1}}e_{0}^{n_{0}-1} $$
By duality they define maps $K \langle\langle e_{0\cup \mu_{N}} \rangle\rangle \rightarrow K \langle\langle e_{0\cup \mu_{N}} \rangle\rangle$, which we will also denote by $\textbf{p}$ and $\textbf{q}$ and, by restriction, they define maps $K \langle \langle Y_{N} \rangle\rangle  K \langle \langle Y_{N} \rangle\rangle$. We have (\cite{Racinet}, Corollaire 2.24, D\'{e}finition 3.1) :
\begin{equation} \label{eq:rel between reg} \textbf{q}\pr(\Phi_{\KZ}) = \exp \bigg( \sum\limits_{n=2}^{\infty} \frac{(-1)^{n}}{n}\zeta(n)y_{1}^{n}
\bigg)(\Phi_{\KZ})_{\ast} .
\end{equation}
This formula can be compared with equations (\ref{eq:multizetas}) and (\ref{eq:multizetas integral}).
\newline\indent The $p$-adic analogue of this formula has been proved  in \cite{Furusho Jafari}, Theorem 0.1, (iii), in the $N=1$ case, and this can be easily adapted to any $N$.

\subsection{Adjoint $p$-adic double shuffle equations}

For $\Phi \in \DS_{0}(K)$, we are going to find some ``adjoint double shuffle equations'' for $\Phi_{\Ad,\chi}$, as defined in Definition 2.2.3. In particular, we will obtain adjoint double shuffle equations for adjoint $p$MZV$\mu_{N}$'s.

\subsubsection{Expression of $\Phi_{\Ad,\chi}$ in terms of $K \langle \langle Y_{N}\rangle\rangle$}

The first step is to rewrite $\pr(f_{\Ad})$ in terms of operations on $K\langle\langle Y_{N}\rangle\rangle$.
\newline\indent Below the hat denotes the completion with respect to the weight grading. The notation below refers to the word ``shifting'', this terminology will be explained in a subsequent paper.

\begin{Definition} \label{def shft} (i) Let 
$\shft_{\ast} : \mathcal{O}^{\ast} \rightarrow  \widehat{\mathcal{O}^{\ast}}$, 
$w(e_{0},(e_{\xi})_{\xi \in \mu_{N}(K)}) \mapsto w (\frac{1}{1+e_{0}}e_{0},(\frac{1}{1+e_{0}}e_{\xi})_{\xi \in \mu_{N}(K)})$.
\newline (ii) For any $l \in \mathbb{N}$, let
$\shft_{l} : \mathcal{O}^{\ast} \mapsto \mathcal{O}^{\ast}$, 
$\displaystyle e_{0}^{n_{d}-1}e_{\xi_{d}} \ldots e_{0}^{n_{1}-1}e_{\xi_{1}} \mapsto \sum\limits_{\substack{l_{1}+\ldots+l_{d}=l \\ l_{1},\ldots,l_{d} \geqslant  0}} \prod_{i=1}^{d}  {-n_{i} \choose l_{i}} e_{0}^{n_{d}+l_{d}-1}e_{\xi_{d}} \ldots e_{0}^{n_{1}+l_{1}-1}e_{\xi_{1}}$. 
\newline Let $\shft^{\vee}_{l} :  K \langle \langle Y_{N}\rangle\rangle \rightarrow  K \langle \langle Y_{N}\rangle\rangle$, defined by for all $w$,
$(\shft_{l}^{\vee}G)[w] = G[\shft_{l}(w)]$.
\newline (iii) \label{def SY} For any $G \in K \langle \langle Y_{N}\rangle\rangle$ and $\xi \in \mu_{N}(K)$, let $G^{\inv,\xi} \in K \langle\langle Y_{N}\rangle\rangle$ be defined by, for all $d\geq 0$,
$G^{\inv,\xi}[y_{l+1}^{(\xi_{d+1})}y_{n_{d}}^{(\xi_{d})} \ldots y_{n_{1}}^{(\xi_{1})}] = \left\{ \begin{array}{ll} (-1)^{n_{1}+\ldots+n_{d}}G[y_{n_{1}}^{(\xi_{2})} \ldots y_{n_{d}}^{(\xi_{d+1})}] & \xi_{1}=\xi
\\ 0 & \xi_{1}\not=\xi \end{array} \right.$.
\end{Definition}

We note that $\shft_{l}$ is the coefficient of $\Lambda^{l}$ in $\widehat{\tau(\Lambda)}^{-1} \shft_{\ast} \tau(\Lambda)$.

\begin{Proposition} \label{prop adjoint star} For any $f\in  \tilde{\Pi}_{1,0}(K)$, then we have 
$$ \pr(f_{\Ad,\chi}) = \sum_{l\geq 0} \sum_{w_{l}} \sum_{\xi \in \mu_{N}(K)} \bigg( \chi(\xi)(-1)^{l} \big(  \shft^{\vee}_{\ast,l} \pr(f^{(\xi)})\big)^{\inv,\xi}  \pr(f^{(\xi)}) \bigg) [w_{l}] w_{l} $$
\end{Proposition}

\begin{proof} Consider a word whose rightmost letter is not $e_{0}$, and write it as $w=e_{0}^{l}e_{\xi_{d+1}}e_{0}^{n_{d}-1}e_{\xi_{d}} \ldots e_{0}^{n_{1}-1}e_{\xi_{1}}$ ($\xi_{i}$'s are roots of unity, $l\geq 0$, $n_{i}\geq 1$). Then, for $f$ a solution to the shuffle equation such that $f[e_{0}]=0$, we have
$$ \sum_{\xi \in \mu_{N}(K)} \chi(\xi)({f^{(\xi)}}^{-1}e_{\xi}f^{(\xi)})[w] = 
\sum_{d'=1}^{d+1} \chi(\xi_{d'})
{f^{(\xi_{d'})}}^{-1}[e_{0}^{l}e_{\xi_{d+1}}e_{0}^{n_{d}-1} \ldots e_{0}^{n_{d'}-1}]
f^{(\xi_{d'})}[e_{0}^{n_{d'-1}-1}e_{\xi_{d'-1}} \ldots e_{0}^{n_{1}-1}e_{\xi_{1}}]
$$
By using the formula for the antipode of the shuffle Hopf algebra we have
$$ {f^{(\xi_{d'})}}^{-1}[e_{0}^{l}e_{\xi_{d+1}}e_{0}^{n_{d}-1} \ldots e_{0}^{n_{d'}-1}] = (-1)^{n_{d'}+\ldots+n_{d}+l}
f^{(\xi_{d'})}[e_{0}^{n_{d'}-1}e_{\xi_{d'+1}} \ldots e_{0}^{n_{d}-1}e_{\xi_{d+1}}e_{0}^{l}]
$$
By using a well-known consequence of the shuffle equation and of $f[e_{0}]=0$, we have
$$ 
f^{(\xi_{d'})}[e_{0}^{n_{d'}-1}e_{\xi_{d'+1}} \ldots e_{0}^{n_{d}-1}e_{\xi_{d+1}}e_{0}^{l}] = 
\sum_{l_{d'}+\ldots+l_{d}=l} 
\prod_{i=d'}^{d} {-n_{i} \choose l_{i}}
f^{(\xi_{d'})}[e_{0}^{n_{d'}+l_{d'}-1}e_{\xi_{d'+1}} \ldots e_{0}^{n_{d}+l_{d}-1}e_{\xi_{d+1}}] $$
By the definition of $\shft_{l}$ (Definition 3.2.1 (iii)), we have
$$ 
\sum_{l_{d'}+\ldots+l_{d}=l} 
\prod_{i=d'}^{d} {-n_{i} \choose l_{i}}
f^{(\xi_{d'})}[e_{0}^{n_{d'}+l_{d'}-1}e_{\xi_{d'+1}} \ldots e_{0}^{n_{d}+l_{d}-1}e_{\xi_{d+1}}] = 
(\shft_{l}^{\vee}f^{(\xi_{d'})})[(e_{0}^{n_{d'}-1}e_{\xi_{d'+1}} \ldots e_{0}^{n_{d}-1}e_{\xi_{d+1}})] $$
Moreover,
$$ (-1)^{l+n_{d'}+\ldots+n_{d}}
(pr f^{(\xi_{d'})})[\shft_{l}(e_{0}^{n_{d'}-1}e_{\xi_{d'+1}} \ldots e_{0}^{n_{d}-1}e_{\xi_{d+1}})] = (-1)^{l} \big( \shft^{\vee}_{\ast,l} \pr(f^{(\xi)})\big) ^{\inv,\xi_{d'}}[e_{0}^{l}e_{\xi_{d+1}}\ldots e_{0}^{n_{d'}-1}e_{\xi_{d'}}] $$
and, finally, since all the words involved have their rightmost letter not equal to $e_{0}$ we can replace everywhere $f^{(\xi)}$ by $\pr(f^{(\xi)})$. In the end, we have 
\begin{multline} \sum_{\xi \in \mu_{N}(K)} \chi(\xi)({f^{(\xi)}}^{-1}e_{\xi}f^{(\xi)})[w]
\\
\begin{array}{l} = \displaystyle \sum_{d'=1}^{d+1} \chi(\xi_{d'}) (-1)^{l} \big( \shft^{\vee}_{l} \pr(f^{(\xi)}) \big)^{\inv,\xi_{d'}}[e_{0}^{l}e_{\xi_{d+1}}\ldots e_{0}^{n_{d'}-1}e_{\xi_{d'}}] \pr(f^{(\xi_{d'})})[e_{0}^{n_{d'-1}-1}e_{\xi_{d'-1}} \ldots e_{0}^{n_{1}-1}e_{\xi_{1}}] 
\\ \displaystyle  = \Big( \sum_{d'=1}^{d+1} \chi(\xi_{d'}) (-1)^{l} (\shft^{\vee}_{l} \pr(f^{(\xi_{d'})}))^{\inv,\xi_{d'}}  \pr(f^{(\xi_{d'})}) \Big) [w] 
\\ \displaystyle  = \Big( \sum_{\xi \in \mu_{N}(K)} \chi(\xi) (-1)^{l} (\shft^{\vee}_{\ast,l} \pr(f^{(\xi)}))^{\inv,\xi}  \pr(f^{(\xi)}) \Big) [w]
\end{array} 
\end{multline}
\end{proof}

\subsubsection{Relation between the two regularizations}

The second step is to observe that the relation between the two regularizations becomes trivial in the adjoint setting. This is an aspect of the proximity between adjoint cyclotomic multiple zeta values and cyclotomic multiple harmonic sums which we are using in this work.
\newline\indent In view of the Lemma \ref{prop adjoint star}, we now define an analogue of $\pr(f_{\Ad})$ in which $f$ is replaced by $f_{\ast}$ defined by the equalities :
\begin{equation} \label{eq:passage a star 1} E_{f} = \exp(\sum_{n=2}^{\infty} \frac{(-1)^{n}}{n} f[e_{0}^{n-1}e_{1}] e_{1}^{n}) 
\end{equation}
\begin{equation} \label{eq:passage a star 2} \textbf{q} \pr(f) = E_{f} f_{\ast} 
\end{equation}

\begin{Definition} For any $f \in \tilde{\Pi}_{1,0}(K)$, let 
$\displaystyle f_{\Ad,\chi,\ast} = \sum_{l\geq 0} \sum_{w_{l}}  \sum_{\xi \in \mu_{N}(K)} \chi(\xi) \big((-1)^{l} \big( \shft^{\vee}_{\ast,l}(f_{\ast}^{(\xi)})\big)^{\inv,\xi}f_{\ast}^{(\xi)} \big)[w_{l}]w_{l}$
where the sum is over $w_{l}$ of the form $e_{0}^{l}e_{\xi_{d+1}}e_{0}^{n_{d}-1}e_{\xi_{d}} \ldots e_{0}^{n_{1}-1}e_{\xi_{1}}=y_{l+1}^{(\xi_{d+1})}y_{n_{d}}^{(\xi_{d})} \ldots y_{n_{1}}^{(\xi_{1})}$, $d\geq 0$.
\end{Definition}

In \cite{I-2}, Definition 1.1.3 and Proposition 1.1.4, we have defined the adjoint Ihara product on $\Ad_{\Pi_{1,0}}(e_{1})$ by the formula $(g,f) \mapsto g \circ_{\Ad}^{\smallint_{1,0}} f = f (e_{0},(g^{(\xi)})_{\xi \in \mu_{N}(K)})$, and we have proved that it is a group law on $\Ad_{\Pi_{1,0}}(e_{1})$.
\newline\indent In the next statement, (i) says that the two regularizations (integrals and series) give the same result for adjoint $p$-adic cyclotomic multiple zeta values : this is coherent with the fact that their adjoint quasi-shuffle relation ((iii) below) can be understooed via cyclotomic multiple harmonic sums, as we are going to see in the subsequent paper.

\begin{Proposition} \label{comparaison reg adjoint}We have 
\begin{equation} \label{eq:comparison of regularisations} \textbf{q}\pr(\Phi_{\Ad,\chi}) = \Phi_{\Ad,\chi,\ast}
\end{equation}
\end{Proposition}

\begin{proof} (i) We are going to simplify the expression of $pr(\Phi)$ given by Lemma \ref{prop adjoint star}. By equation (\ref{eq:passage a star 2}) we have, for each $\xi \in \mu_{N}(K)$ and $l\geq 0$,
$$ (-1)^{l}\shft^{\vee}_{l} (pr(\Phi^{(\xi)})^{\inv,\xi} ) \pr(\Phi^{(\xi)}) = (-1)^{l}\shft^{\vee}_{l}\big( E^{(\xi)}(\Phi_{\ast}^{(\xi)})^{\inv,\xi}) \big) E^{(\xi)} \Phi_{\ast}^{(\xi)} $$
Let now a word $w=e_{0}^{l}e_{\xi_{d+1}}e_{0}^{n_{d}-1}e_{\xi_{d}} \ldots e_{0}^{n_{1}-1}e_{\xi_{1}}$. By the definition of the map $G \mapsto G^{\inv,\xi}$ (Definition \ref{def shft} (iii)), for all $d'$ ($1 \leq d' \leq d+1$), assuming $\xi=\xi_{d'}$, we have
\begin{multline*} (-1)^{l}\big( \shft_{l}^{\vee}(E_{\Phi}^{(\xi)}\Phi_{\ast}^{(\xi)}) \big)^{\inv,\xi}[e_{0}^{l}e_{\xi_{d+1}}e_{0}^{n_{d}-1}e_{\xi_{d}} \ldots e_{0}^{n_{d'}-1}e_{\xi_{d'}}] 
\\ =
(-1)^{l+n_{d}+\ldots+n_{d'}}
\big( \shft_{l}^{\vee}(E^{(\xi)}\Phi_{\ast}^{(\xi)})\big) [e_{0}^{n_{d'}-1}e_{\xi_{d'+1}} \ldots e_{0}^{n_{d}-1}e_{\xi_{d+1}}] 
\end{multline*}
by the definition of $\shft_{l}$ (Definition 3.2.1, (iii)), and by the fact that $E_{\Phi}^{(\xi)}$ has non-zero coefficients only at certain words which are of the form $e_{\xi}^{n}$,  we have
\begin{multline} \label{eq:3} \big(\shft^{\vee}_{\ast,l}\big( E_{\Phi}^{(\xi)}(\Phi_{\ast}^{(\xi)}) \big)[e_{0}^{n_{d'}-1}e_{\xi_{d'+1}} \ldots e_{0}^{n_{d}-1}e_{\xi_{d+1}}] 
\\
\begin{array}{l} \displaystyle
= \big(E_{\Phi}^{(\xi)}
\shft^{\vee}_{\ast,l}(\Phi_{\ast}^{(\xi)}) \big) [e_{0}^{n_{d'}-1}e_{\xi_{d'+1}} \ldots e_{0}^{n_{d}-1}e_{\xi_{d+1}}]
\\ \displaystyle
= \sum_{r=d'}^{d+1} E_{\Phi}^{(\xi)}[e_{0}^{n_{d'}-1}e_{\xi_{d'+1}} \ldots e_{0}^{n_{r}-1}e_{\xi_{r+1}}]
\shft^{\vee}_{\ast,l}(\Phi_{\ast}^{(\xi)})  [e_{0}^{n_{r+1}-1}e_{\xi_{r+2}} \ldots e_{0}^{n_{d}-1}e_{\xi_{d+1}}] 
\end{array}
\end{multline}
By the definition of the map $G \mapsto G^{\inv,\xi}$ (Definition \ref{def shft} (iii)), we have 
$$ (-1)^{n_{r+1}+\ldots+n_{d}} \shft^{\vee}_{l}(\Phi_{\ast}^{(\xi)})  [e_{0}^{n_{r+1}-1}e_{\xi_{r+2}} \ldots e_{0}^{n_{d}-1}e_{\xi_{d+1}}] = \big(
\shft^{\vee}_{l}(\Phi_{\ast}^{(\xi)}) \big)^{(\inv,\xi_{r+1})} [e_{0}^{l}e_{\xi_{d+1}} \ldots e_{0}^{n_{r+2}-1}e_{\xi_{r+1}}] $$
Since $E_{\Phi}^{(\xi)}$ has non-zero coefficients only at words of the form $e_{\xi}^{n}$, in the sum over $r$ in (\ref{eq:3}), a non-zero term can appear only if $\xi_{r+1}=\xi_{d'}$, and we can write
$$ E_{\Phi}^{(\xi)}[e_{0}^{n_{d'}-1}e_{\xi_{d'+1}} \ldots e_{0}^{n_{r}-1}e_{\xi_{r+1}}] = E_{\Phi}^{(\xi)}[e_{\xi_{r+1}}e_{0}^{n_{r}-1} \ldots e_{\xi_{d'+1}}e_{0}^{n_{d'}-1}]
= E_{\Phi}^{(\xi)}[e_{0}^{n_{r}-1} \ldots e_{0}^{n_{d'}-1}e_{\xi_{d'}}]  $$
Below we use the notation $\tau$ introduced in (\ref{eq:tau}). By the two previous equalities, (\ref{eq:3}) multiplied by $(-1)^{l+n_{d}+\ldots+n_{d'}}$ is equal to
\begin{multline*} \sum_{r=d'}^{d+1} (-1)^{l}
\big(\shft^{\vee}_{\ast,l}(\Phi_{\ast}^{(\xi_{d})}) \big)^{(\inv,\xi_{d'})} [e_{0}^{l}e_{\xi_{d+1}} \ldots e_{0}^{n_{r+1}-1}e_{\xi_{r+1}}] \big( \tau(-1)E_{\Phi}^{(\xi)}[e_{0}^{n_{r}-1} \ldots e_{0}^{n_{d'}-1}e_{\xi_{d'}}] 
\\ =
(-1)^{l}\bigg( \big(\shft^{\vee}_{\ast,l}(\Phi_{\ast}^{(\xi_{d})}) \big)^{(\inv,\xi_{d'})} (\tau(-1)E^{(\xi)}) \bigg)[e_{0}^{l}e_{\xi_{d+1}} \ldots e_{\xi_{d'+1}} e_{0}^{n_{d'}-1}e_{\xi_{d'}}] 
\end{multline*}
This formula combined to Lemma \ref{prop adjoint star} shows that 
\begin{multline} \label{eq:avant derniere} \pr(\Phi)[w] 
= \sum_{d'=1}^{d+1} \bigg(  (-1)^{l} \xi_{d'}^{-1} 
\big(\shft^{\vee}_{\ast,l}(\Phi_{\ast}^{(\xi_{d'})}) \big)^{(\inv,\xi_{d'})} (\tau(-1)E^{(\xi_{d'})}) E^{(\xi_{d'})} \Phi_{\ast}^{(\xi_{d'})}  \bigg)[w] 
\\ =  \sum_{\xi \in \mu_{N}(K)} \bigg(  (-1)^{l} \xi^{-1}
\big(\shft^{\vee}_{\ast,l}(\Phi_{\ast}^{(\xi)}) \big)^{(\inv,\xi)} (\tau(-1)E^{(\xi)}) E^{(\xi)} \Phi_{\ast}^{(\xi)}  \bigg)[w]
\end{multline}
For any $S \in K\langle \langle e_{0\cup \mu_{N}}\rangle\rangle$, and $n \geq 0$, we have $(\tau(-1)S)^{n} = \tau(-1)(S^{n})$ because for any words $w_{i}$, $\displaystyle \sum_{i=1}^{n}\weight(w_{i}) = \weight(w_{1}\ldots w_{n})$. Thus, if the coefficient of the empty word in $S$ is $0$, we can write
$\displaystyle \tau(-1)\exp(S) = \sum_{n=0}^{\infty} \tau(-1)\frac{S^{n}}{n!} = \sum_{n=0}^{\infty} \frac{(\tau(-1)S)^{n}}{n!} = \exp (\tau(-1)S)$. In particular, by (\ref{eq:passage a star 1}) $\displaystyle \tau(-1)E_{\Phi} = 
\exp \big(\sum_{n=2}^{\infty} \frac{1}{n}\Phi[e_{0}^{n-1}e_{1}]e_{1}^{n} \big)$. Thus
$$(\tau(-1)E_{\Phi}) E_{\Phi}  = \exp \big(\sum_{n=2}^{\infty} \frac{1}{n}\Phi[e_{0}^{n-1}e_{1}]e_{1}^{n} \big)  \exp \big(\sum_{n=2}^{\infty} \frac{(-1)^{n}}{n}\Phi[e_{0}^{n-1}e_{1}]e_{1}^{n} \big)  = 
\exp \big(\sum_{n=2}^{\infty} \frac{(1+(-1)^{n})}{n}\Phi[e_{0}^{n-1}e_{1}]e_{1}^{n} \big) $$
For $n$ odd, we have $1+(-1)^{n}=0$ ; for $n$ even, since $\Phi \in \DS_{0}(K)$ we have $\Phi[e_{0}^{n-1}e_{1}]=0$. This proves
\begin{equation} \label{eq:derniere} \tau(-1)(E_{\Phi})E_{\Phi} = 1
\end{equation} 
The result follows from (\ref{eq:avant derniere}) and (\ref{eq:derniere}).
\end{proof}

\subsubsection{Definition of the adjoint double shuffle relations}

\begin{Proposition} \label{221}Let $\Phi \in \DS_{0}(K)$ and $\chi \in \Theta$.
\newline (i) We have
\begin{equation} \label{eq:ds adjoint 1}
\text{for all non-empty words }w,w',\text{ } \Phi_{\Ad,\chi}[w \text{ }\sh \text{ }w'] = 0 .
\end{equation}
\noindent (ii) We have 
\begin{multline} \label{eq:ds adjoint 2}
\text{for all words }w,w'\text{ and }L \in \mathbb{N},
 \text{ }\sum_{\substack{l,l'\geqslant 0 \\ l+l'=L}} \Phi_{\Ad,\chi,\ast}[w;l] \Phi_{\Ad,\chi,\ast}[w';l'] = \Phi_{\Ad,\chi,\ast}[ w \ast_{\har} w';L] .
\end{multline}
\end{Proposition}

\begin{proof} (i) Let $\Delta_{\sh}$ be the shuffle coproduct. For all $\xi \in \mu_{N}(K)$, we have 
$\Delta_{\sh}(e_{\xi})= e_{\xi} \otimes 1 + 1 \otimes e_{\xi}$, and $\Delta_{\sh}(\Phi^{(\xi)})=\Phi^{(\xi)} \otimes \Phi^{(\xi)}$, whence $\Delta( {\Phi^{(\xi)}}^{-1} e_{\xi} \Phi^{(\xi)}) = {\Phi^{(\xi)}}^{-1} e_{\xi} \Phi^{(\xi)} \otimes 1 + 1 \otimes {\Phi^{(\xi)}}^{-1} e_{\xi} \Phi^{(\xi)}$.
This implies $\Delta_{\sh}(\Phi_{\Ad}) = \Phi_{\Ad} \otimes 1 + 1 \otimes \Phi_{\Ad}$.
\newline\indent (ii) Let us prove that $\shft_{\ast}$ is a morphism of quasi-shuffle algebras. For simplicity, we do the proof for $\mathbb{P}^{1} - \{0,1,\infty\}$. The dual of $\shft_{\ast}$ is the concatenation algebra morphism $\imath_{\ast}^{\vee}$ defined by 
$$ y_{n} \mapsto \Lambda^{n} \sum_{l=0}^{n-1} \Lambda^{l} {n-1 \choose l} (-1)^{n-l} y_{n-l} 
= \Lambda^{n} \sum_{l=1}^{n} \Lambda^{n-l} (-1)^{l} y_{l} {n-1 \choose n-l} . $$
\noindent We have 
$(\imath^{\vee} \otimes \imath^{\vee})\Delta_{\ast}(y_{n}) 
= 
1 \otimes \imath^{\vee}(y_{n}) 
+ 
\imath^{\vee}(y_{n}) \otimes 1 
+
\sum_{k=1}^{n-1} \imath^{\vee}(y_{k}) \otimes \imath^{\vee}(y_{n-k})$ ; the third term of this sum is 
$$ \Lambda^{n} \sum_{k=1}^{n-1} \big( 
\sum_{l=1}^{k} \Lambda^{k-l} {k-1 \choose k-l} (-1)^{l} y_{l}
\big)
\otimes 
\big( 
\sum_{l=1}^{n-k} \Lambda^{n-k-l} {n-k-1 \choose n-k-l'} (-1)^{l'} y_{l'}
\big) $$
$$ = \Lambda^{n}
\sum_{L=2}^{n} \Lambda^{n-L}(-1)^{n-L} 
\big( 
   \sum_{\substack{l+l'=L\\ l,l'\geqslant  1}} y_{l} \otimes y_{l'}
   \big)
\sum_{\substack{l\leqslant k \leqslant n-L+l}} {k-1 \choose k-l} {n-k-1 \choose n-k-l'}   
$$    
\noindent and for all $l,l'$ such that $l+l'=L$, we have
$$
\sum_{\substack{l\leqslant k \leqslant n-L+l}} {k-1 \choose k-l} {n-k-1 \choose n-k-l'}    
= 
\sum_{k'=0}^{n-L} {k' + l-1 \choose k'} {n-L-k'+l'-1 \choose n-L-k'}
= {n-L+L-1 \choose n-L} .
$$
\noindent On the other hand, $S_{Y}$ is an anti-morphism of quasi-shuffle algebras. This gives the result.
\end{proof}

\begin{Definition} We call (\ref{eq:ds adjoint 1}) the \emph{adjoint shuffle equation} and (\ref{eq:ds adjoint 2}) the \emph{adjoint quasi-shuffle equation} ; and we call the collection of (\ref{eq:ds adjoint 1}) and (\ref{eq:ds adjoint 2}) the \emph{adjoint double shuffle equations}.

We denote by $\DS_{0,\Ad}(K)$ the set of $\psi \in K \langle \langle e_{0\cup\mu_{N}} \rangle\rangle$ which satisfy the adjoint shuffle equation and such that for all $\chi \in \Theta$, $\textbf{q} \pr \Moy_{\chi}(\psi)$ satisfies the adjoint quasi-shuffle equation.

These equations define an affine scheme $\DS_{0,\Ad}$ over $k_{N}$, which we call the \emph{adjoint double shuffle scheme}.
\end{Definition}

By the previous propositions, we have proved that if $\Phi \in \DS_{0}(K)$ then $\Ad_{\Phi}(e_{1})$ is in $\DS_{0,\Ad}(K)$. 
We note that the adjoint shuffle equation is also known as the linearlized shuffle equation, or the shuffle equation modulo products, and amounts to say that $\Ad_{\Phi}(e_{1})$ is primitive for the shuffle coproduct $\Delta_{\sh}$. We also note that the adjoint quasi-shuffle equation can be reformulated as saying that $\Ad_{\Phi}(e_{1})$ is a "grouplike" element for a certain adjoint variant $\Delta_{\ast}^{\Ad}$ of $\Delta_{\ast}$.

\subsubsection{Stability by the Ihara product}

We now deduce from Racinet's theorem that $\DS_{0}$ is a group for the Ihara product \cite{Racinet} an adjoint variant of this statement.

\begin{Proposition} The image of the map $\Ad(e_{1}) : \DS_{0} \mapsto \DS_{0,\Ad}$ is an algebraic group with the adjoint Ihara product $\circ_{\Ad}^{\smallint_{1,0}}$.
\end{Proposition}

\begin{proof} By Racinet's theorem \cite{Racinet}, $(\DS_{0},\circ^{\smallint_{1,0}})$ is an algebraic group. By \cite{I-2}, Proposition 1.1.4, $\Ad(e_{1})$ is a morphism of algebraic groups $(\tilde{\Pi}_{1,0},\circ^{\smallint_{1,0}}) \simlra (\Ad_{\tilde{\Pi}_{1,0}}(e_{1}),\circ_{\Ad}^{\smallint_{1,0}})$.
\end{proof}

The following question seems interesting, as we will explain in the next section :

\begin{Question} \label{question} Is $\Ad(e_{1}) : \DS_{0} \mapsto \DS_{0,\Ad}$ an isomorphism ?
\end{Question}

\subsection{Adjoint complex double shuffle equations}

\subsubsection{Adjoint double shuffle equation\label{paragraph lifts}}

Let $K$ be a field of characteristic $0$. For $\mu \in K - \{0\}$, let $\DS_{\mu}$ be the scheme of regularized double shuffle relations with parameter $\mu$ defined by Racinet \cite{Racinet}.

In order to write the relation between the two regularizations of adjoint cyclotomic multiple zeta values, we have to consider a slightly more general object.

Let $\rho_{\Ad}$ be the linear map $K[T_{1},T_{2}] \rightarrow K[T_{1},T_{2}]$ defined by the formula 

$$ \rho (e^{T_{2}e_{1}} e^{T_{1}e_{1}}) \mapsto (-1)^{\weight}(e^{T_{1}e_{1}})(-1)^{weight}(E)E e^{T_{2}e_{1}} $$ 

We now define the regularized versions of adjoint complex cyclotomic multiple zeta values. Let the two regularizations $\Phi_{\KZ}(T)$ and $\Phi_{\KZ,\ast}(T)$ of $\Phi_{\KZ}$ resp. $\Phi_{\KZ,\ast}$ which can be found in \cite{Racinet} ($T$ is a formal variable). We define their adjoint analogues, the two regularizations of $\Phi_{\KZ,\Ad,\chi}$ :

(a) The regularization in the sense of integrals :
$\displaystyle \Phi_{\KZ,\Ad,\chi}(T) = \sum\limits_{\xi \in \mu_{N}(K)} \chi(\xi) {\Phi(T)^{-1}}^{(\xi)}e_{\xi}\Phi(T)^{(\xi)} $

(b) The regularization in the sense of series (in view of Proposition \ref{prop adjoint star}) :
\newline 
$\displaystyle \Phi_{\KZ,\Ad,\chi,\ast}(T) = \sum_{l\geq 0} \sum_{w_{l}}  \sum_{\xi \in \mu_{N}(K)} \chi(\xi) \big((-1)^{l} \big( \shft^{\vee}_{\ast,l}(\Phi_{\ast}(T)^{(\xi)})\big)^{\inv,\xi} \Phi_{\ast}^{(\xi)}(T) \big)[w_{l}]w_{l} \in \mathbb{C} \langle\langle Y_{N} \rangle\rangle$

By considering their coefficients, we deduce the regularized version of adjoint MZV$\mu_{N}$'s :

\begin{Definition} (a) The regularized (in the sense of integrals) AdMZV$\mu_{N}$'s are 
	$$ \zeta^{\Ad}\big( (n_{i})_{d};(\xi_{i})_{d+1};l;\chi)(T) = 
	\Phi_{\KZ,\Ad,\chi}(T) [e_{0}^{l}e_{\xi_{d+1}}e_{0}^{n_{d}-1}e_{\xi_{d}} \ldots e_{0}^{n_{1}-1}e_{\xi_{1}}]. $$
	(b) The regularized (in the sense of series) AdMZV$\mu_{N}$'s are 
	$$ \zeta^{\Ad}\big( (n_{i})_{d};(\xi_{i})_{d+1};l;\chi)(T) = 
	\Phi_{\KZ,\Ad,\chi,\ast}(T) [e_{0}^{l}e_{\xi_{d+1}}e_{0}^{n_{d}-1}e_{\xi_{d}} \ldots e_{0}^{n_{1}-1}e_{\xi_{1}}]. $$
	And similarly for the $\Lambda$-adic AdMZV$\mu_{N}$'s.
\end{Definition}

We generalize the previous definitions as follows :
$$ \Phi_{\KZ,\chi,\Ad}(T_{1},T_{2}) = \sum_{\xi \in \mu_{N}(K)} \chi(\xi) {\Phi(T_{2})^{-1}}^{(\xi)}e_{\xi}\Phi(T_{1})^{(\xi)} $$
$$ \Phi_{\KZ,\ast,\chi,\Ad}(T_{1},T_{2}) = \sum_{l\geq 0} \sum_{w_{l}}  \sum_{\xi \in \mu_{N}(K)} \chi(\xi) \big((-1)^{l} \big( \shft^{\vee}_{\ast,l}(\Phi_{\ast}(T_{2})^{(\xi)})\big)^{\inv,\xi} \Phi_{\ast}^{(\xi)}(T_{1}) \big)[w_{l}]w_{l} $$

\begin{Proposition} If $\Phi \in \DS_{\mu}(K)$, we have 
	\newline (i) $\Phi_{\Ad,\chi}(T)$ satisfies the adjoint shuffle equation
	\newline (ii) $\Phi_{\Ad,\chi,\ast}(T)$ satisfies the adjoint quasi-shuffle equation
	\newline (iii) We have $\textbf{q} \pr \rho (\Phi_{\Ad,\chi}(T_{1},T_{2})) = \Phi_{\Ad,\chi,\ast}(T_{1},T_{2})$, 
	$\Phi_{\Ad,\chi}(T) = \Phi_{\Ad,\chi}(T,T)$, $ \Phi_{\Ad,\chi,\ast}(T) = \Phi_{\Ad,\chi,\ast}(T,T)$
	\newline (iv) We have $\Phi_{\Ad,\chi}[e_{1}e_{0}e_{1}] = (\Phi^{-1}e_{1}\Phi)[e_{1}e_{0}e_{1}] = 2\Phi[e_{0}e_{1}]$.
\end{Proposition}

\begin{proof} Similar to the proofs in the $p$-adic case (\S3.2). (iv) follows from $\Phi^{-1}[e_{1}e_{0}] = (-1)^{2}\Phi[e_{0}e_{1}]$ and $({\Phi^{-1}}^{(\xi)}e_{\xi}\Phi^{(\xi)})[e_{1}e_{0}e_{1}]=0$ for $\xi\not=1$.
\end{proof}

\begin{Definition} Let $\DS_{\Ad,\eta}(\mathbb{C})$ be the set of 
$(\Phi_{\Ad,\chi}(T_{1},T_{2}),\Phi_{\Ad,\chi,\ast}(T_{1},T_{2})) \in \mathbb{C}(T)\langle \langle e_{0\cup \mu_{N}} \rangle\rangle^{2}$ which satisfy the above equations and 
$\Phi_{\Ad,\chi}[e_{1}e_{0}e_{1}]=2 (-\mu^{2}/24)$.
\newline The above equations define an affine scheme which we call the \emph{scheme of adjoint regularized double shuffle relations}.
\end{Definition}

We deduce :

\begin{Corollary} The map $\Ad(e_{1})$ defines a morphism $\DS_{\mu} \rightarrow \DS_{\mu,\Ad}$.
	
	By the adjoint Ihara action, the image of the map $\DS_{\mu} \rightarrow \DS_{\mu,\Ad}$ is a torsor under the group defined as the image of the map $\DS_{0} \rightarrow \DS_{0,\Ad}$.
\end{Corollary}

\begin{proof} This follows from Racinet's theorem that $\DS_{\mu}$ is a torsor under the group $\DS_{0}$ for the Ihara action, and that $\Ad(e_{1})$ is a morphism of algebraic groups $(\tilde{\Pi}_{1,0},\circ^{\smallint_{1,0}}) \simlra (\Ad_{\tilde{\Pi}_{1,0}}(e_{1}),\circ_{\Ad}^{\smallint_{1,0}}))$ (\cite{I-2}, Proposition 1.1.4).
\end{proof}

Particular cases of our adjoint regularized double shuffle relations can be found in \cite{Hi} and \cite{HMS} which appeared recently on the arXiv.

\begin{Remark} We note that $$ (-1)^{\weight}(E)E= \exp(\sum_{n\geq 2,\text{even}} \frac{2}{n} \Phi[e_{0}^{n-1}e_{1}]e_{1}^{n} ) = \exp(\sum_{n\geq 1} \frac{1}{n} \Phi[e_{0}^{2n-1}e_{1}]e_{1}^{2n} ) $$
	In the context of \S3, this was equal to 1. Here, in the case where $\Phi=\Phi_{\KZ}$, we have 
	$\Phi[e_{0}^{2n-1}e_{1}] = \frac{(-1)^{n-1} B_{2n}}{2.(2n)!} (24\Phi[e_{0}e_{1}])^{n}$. It follows from \cite{Ihara Kaneko Zagier}, Theorem 7 that this relation holds for all solutions to the double shuffle equations. Thus 
	$$ (-1)^{\weight}(E)E =  \exp(\sum_{n\geq 1} \frac{1}{n} \frac{(-1)^{n-1} B_{2n}}{2.(2n)!} (24\Phi[e_{0}e_{1}])^{n}e_{1}^{2n} )=  \exp(\sum_{n\geq 1}  \frac{- B_{2n}}{2n.(2n)!} (-24\Phi[e_{0}e_{1}]e_{1}^{2})^{n}) $$
	This term appears implicitly in the relation between the two regularizations of adjoint cyclotomic multiple zeta values.
\end{Remark}

\subsubsection{Other aspects of the adjoint double shuffle relations\label{paragraph lifts}}

We give two $\Lambda$-adjoint formulations of the above adjoint shuffle equation.

\begin{Proposition} \label{4.8} For all $w,w' \in \mathcal{O}^{\ast}$, and for $\Lambda, \Lambda'$ formal variables, we have, for all $w,w' \in \mathcal{O}^{\ast}$, and for any $\xi, \xi' \in \mu_{N}(\mathbb{C})$
	\begin{multline}
	(\Phi_{\KZ}^{-1}e^{2\pi i e_{1}} \Phi_{\KZ}) (\frac{1}{1-\Lambda e_{0}}e_{\xi}w)\text{ }
	(\Phi_{\KZ}^{-1}e^{2\pi i e_{1}} \Phi_{\KZ})(\frac{1}{1-\Lambda' e_{0}}e_{\xi'}w')
	\\ = 
	(\Phi_{\KZ}^{-1}e^{2\pi i \xi} \Phi_{\KZ}) \bigg( \frac{1}{1-(\Lambda+\Lambda') e_{0}} e_{\xi} \big( w  \text{ }\sh  \text{ } (\frac{1}{1-\Lambda' e_{0}}e_{\xi'}w') \big)
	+ e_{\xi'} \big( (\frac{1}{1-\Lambda e_{0}}e_{\xi}w) \text{ }\sh \text{ }w' \big) \bigg),
	\end{multline}
	\begin{equation}
	(\Phi_{\KZ}^{-1}e_{1} \Phi_{\KZ}) \bigg( \frac{1}{1-(\Lambda+\Lambda') e_{0}} e_{\xi}( w  \text{ }\sh  \text{ } (\frac{1}{1-\Lambda' e_{0}}e_{\xi'}w'))
	+ e_{\xi'}((\frac{1}{1-\Lambda e_{0}}e_{\xi}w) \text{ }\sh \text{ }w')\bigg) = 0.
	\end{equation}
\end{Proposition}

\begin{proof} Below we use using Notation \ref{definition sigma inv DR} for $\comp^{\Lambda \Ad,\Ad}$. (a) Let us show that 
	\begin{multline} \label{eq:thing to prove} \comp^{\Lambda \Ad,\Ad}e_{\xi}w
	\text{ }\sh\text{ }\comp^{\Lambda'\Ad,\Ad}e_{\xi'}w' =
	\\ 
	\comp^{(\Lambda+\Lambda')\Ad, \Ad}e_{\xi}( w  \text{ }\sh  \text{ } \comp^{\Lambda'\Ad,\Ad}e_{\xi'}w')
	+ 
	\comp^{(\Lambda+\Lambda')\Ad,\Ad}
	e_{\xi'}(\comp^{\Lambda\Ad,\Ad}e_{\xi}w \text{ }\sh \text{ }w') .
	\end{multline}
	\noindent Let $L$ be the left-hand side in (\ref{eq:thing to prove}). We compute the image of $\partial_{e_{x}}(L)$ for all $x \in \{0\} \cup \mu_{N}(K)$ (for the definition of $\partial_{e_{x}}$, see the proof of Proposition \ref{prop 2.3.1}). Given that $\partial_{e_{\xi}}$ is a derivation for $\sh$, we obtain (where $\delta$ is Kronecker's symbol) :
\begin{equation*}
\begin{array}{l}
\partial_{e_{0}}L = 
(\Lambda+\Lambda')L,
\\
\partial_{e_{\xi}}L = (w\text{ }\sh\text{ } \comp^{\Lambda'\Ad,\Ad}e_{\xi'}w') + \delta_{\xi,\xi'}( \comp^{\Lambda\Ad,\Ad}e_{\xi}w \text{ }\sh\text{ }w'),
\\
\partial_{e_{\xi'}}L = \delta_{\xi,\xi'} (w\text{ } \sh\text{ } \comp^{\Lambda'\Ad,\Ad}e_{\xi'}w') + ( \comp^{\Lambda\Ad,\Ad}e_{\xi}w\text{ }\sh\text{ } w'),
\\ 
\partial_{e_{x}}L = 0 \text{ }\text{if}\text{ }x \not\in \{0,\xi,\xi'\}.
\end{array}
\end{equation*}
We deduce (\ref{eq:thing to prove}) by $L=\sum\limits_{x\in \{0\} \cup \mu_{N}(K)} e_{x}\partial_{e_{x}}(L)$.
\newline\indent (b) On the other hand, $f^{-1}e^{2i \pi e_{1}} f$ resp. $f^{-1}e_{1} f$ satisfies the shuffle equation, resp. the shuffle equation modulo products. We deduce the result.
\end{proof}

\subsubsection{Adjoint and harmonic double shuffle equations for multiple polylogarithms}

The shuffle relation of multiple polylogarithms is true by their definition in terms of the iterated integrals (\ref{eq:Li}).
\newline\indent A quasi-shuffle relation for multiple polylogarithms is still true by their power series expansion (\ref{eq:Li series bis}). However, it involves several curves at the same time : if $D,D'$ are two finite subsets of $\mathbb{P}^{1}(K)$, both containing $0$ and $\infty$, the quasi-shuffle equation expressing product of multiple polylogarithms associated respectively with $\mathbb{P}^{1} - D$ and $\mathbb{P}^{1} - D'$ will involve $\mathbb{P}^{1} - DD'$ where $DD' = \{0,\infty\} \cup \{xx'\text{ | }\text{ } x \in D - \{0,\infty\} , x' \in D' - \{0,\infty\} \}$. We leave the exact definitions to the reader.
\newline\indent In the end, by imitating the proofs of \S3, we obtain :

\begin{Proposition} (rough version)
	The adjoint multiple polylogarithms and multiple harmonic polylogarithms resp. finite multiple polylogarithms satisfy a generalization of the double shuffle equations of \S3 resp. \S6.
\end{Proposition}

\subsection{Harmonic double shuffle equations}

We define harmonic double shuffle equations for cyclotomic multiple harmonic values of Definition \ref{def harmonic}, in the frameworks $\smallint_{1,0}$ (\S3.3.1), $\int$ (\S3.3.2), $\Sigma$ (\S3.3.3), and we prove that they are equivalent (\S3.3.4).

\subsubsection{In the framework $\smallint_{1,0}$}

We use the adjoint double shuffle equations of \S3.2 to define the harmonic double shuffle equations in the framework of $\smallint_{1,0}$. In \cite{I-3}, Definition 2.1.2, we have introduced $\circ_{\har}^{\smallint_{1,0}}$, the pro-unipotent harmonic action of integrals at (1,0).

\begin{Proposition} \label{prop 2.3.1}Let $\psi \in K \langle \langle e_{0 \cup \mu_{N}} \rangle\rangle$ such that $\psi[e_{0}]=0$. Let $h=\comp^{\Lambda \Ad,\Ad}\psi$, i.e. $h= \sum\limits_{w} \psi[\frac{1}{1-\Lambda e_{0}}w]w$ where the sum is over words $w$ of the form $e_{\xi_{d+1}} e_{0}^{n_{d}-1}e_{\xi_{d}} \ldots e_{0}^{n_{1}-1} e_{\xi_{1}}$.
\newline (i) $\psi$ satisfies the adjoint quasi-shuffle equation if and only if $h$ satisfies
\begin{equation} \label{eq: DS har int1,0 1}
\text{for all words }w,w',\text{ }h(w \ast w') = h(w) \text{ }h(w') .
\end{equation}
\noindent (ii) $\psi$ satisfies the adjoint shuffle equation if and only if $h$ satisfies
\begin{equation} \label{eq: DS har int1,0 2} 
\text{for all words }w,w'\text{ and for all }n \in \mathbb{N}^{\ast},\text{ }
h (e_{\xi'} (e_{0}^{n-1}e_{\xi} w \text{ } \sh \text{ } w' )) =
h ( e_{\xi} ( w \text{ }\sh \text{ }\frac{1}{1 - \Lambda e_{0}} e_{0}^{n-1}e_{\xi'}w') ) .
\end{equation}
\noindent (iii) The set 
$\{ \comp^{\Lambda \Ad,\Ad}\psi \text{ }|\text{ }\psi \in \Ad_{\DS_{0}(K)}(e_{1})\}$ is a torsor under the group $(\Ad_{\DS_{0}(K)}(e_{1}),\circ_{\Ad}^{\smallint_{1,0}})$ for the pro-unipotent harmonic action $\circ_{\har}^{\smallint_{1,0}}$.
\end{Proposition}

\begin{proof} (i) The adjoint quasi shuffle equation amounts to
$$ \big(\sum_{L \geqslant 0} \psi_{\ast}[w;l]\Lambda^{l} \big) \big(\sum_{L\geqslant 0} \psi_{\ast}[w';l']\Lambda^{l'} \big) = \sum_{L \geqslant 0} \Lambda^{L} \sum_{l+l'=L} \psi[w;l] \psi[w';l'] = \sum_{L \geqslant 0}\Lambda^{L} \psi[ w \ast w';L], $$
i.e. $h(w)h(w') = h(w \ast w')$.
\newline (ii) Let us prove that, for all $w,w'$ words, and $n \in \mathbb{N}^{\ast}$, we have :
\begin{multline} \label{eq:equation 4.10}
-\comp^{\Lambda \Ad,\Ad} \bigg( (e_{0}^{n-1}e_{\xi}w) \text{ }\sh\text{ } w' - w \text{ }\sh\text{ } \shft_{\ast}(e_{0}^{n-1}e_{\xi'})w') \bigg)
= 
\\
\sum_{t=0}^{n-1} \big( e_{0}^{t}e_{\xi} w \big) \text{ }\sh\text{ } \big( (-1)^{n-t}(\shft_{\ast}(e_{0}^{n-1-t}e_{\xi'})w' \big).
\end{multline}
For any $x \in \{0\} \cup \mu_{N}(K)$, let $\partial_{e_{x}}$, resp. $\tilde{\partial}_{e_{x}}$
be the unique linear maps $\mathcal{O}^{\sh} \rightarrow \mathcal{O}^{\sh}$ that send the empty word to $0$ and defined on the other words by 
$\partial_{e_{x}}(e_{x'}w) =\tilde{\partial}_{e_{x}}(we_{x'}) = \left\{ \begin{array}{l} w  \text{ if }x = x'
\\  0 \text{ if }x \not= x'
\end{array} \right.$. For all $w \in \mathcal{O}^{\sh}$, we have $w = \sum\limits_{x \in \{0\} \cup \mu_{N}(K)} e_{x}\partial_{e_{x}}(w) = \sum\limits_{x \in \{0\} \cup \mu_{N}(K)} \tilde{\partial}_{e_{x}}(w)e_{x}$, and  $\partial_{e_{x}}$ and $\tilde{\partial}_{e_{x}}$ are derivations for the shuffle product.
\newline\indent Let $R$ be the right-hand side in (\ref{eq:equation 4.10}). We have $R= e_{0} \partial_{e_{0}}R + \sum\limits_{\xi \in \mu_{N}(K)} e_{\xi}\partial_{e_{\xi}}(R)$. In order to prove (\ref{eq:equation 4.10}) it suffices to show the equalities :
if $\xi \not= \xi'$,
$\left\{\begin{array}{l}\partial_{e_{\xi}}(R) =  w \text{ }\sh\text{ } (-1)^{n}\shft_{\ast}(e_{0}^{n-1}e_{\xi'}w') 
\\ \partial_{e_{\xi'}}(R) =  - (e_{0}^{n-1}e_{\xi}w)\text{ }\sh\text{ } w'\end{array} \right.$ ; if $\xi=\xi'$, $\partial_{e_{\xi}}(R) = w \text{ }\sh\text{ } (-1)^{n}\shft_{\ast}(e_{0}^{n-1}e_{\xi}w')  - (e_{0}^{n-1}e_{\xi}w) \text{ }\sh\text{ } w'$, and $\partial_{e_{0}}(R) = \Lambda R$.
\newline The two first ones are clear ; let us show the last one. One has :
\begin{multline} \label{eq:e0x}
\partial_{e_{0}}(R) = 
\sum_{t=1}^{n-1} (e_{0}^{t-1}e_{\xi}w) 
\text{ }\sh\text{ }
\frac{1}{(\Lambda e_{0}-1)^{n-t}} e_{0}^{n-1-t} e_{\xi'}w'
\\ 
+ \sum_{t=0}^{n-2} (e_{0}^{t}e_{\xi}w) 
\text{ }\sh\text{ }
\frac{1}{(\Lambda e_{0}-1)^{n-t}} e_{0}^{n-2-t} e_{\xi'}w'
+ (e_{0}^{n-1}e_{\xi}w)
\sh
\frac{\Lambda}{\Lambda e_{0}-1} e_{\xi'}w' .
\end{multline}
The sum of the two first terms of the right hand side of (\ref{eq:e0x}) equals
$$ \sum_{t=0}^{n-2} (e_{0}^{t}e_{\xi}w)
\text{ }\sh\text{ }
\bigg[1 + \frac{1}{\Lambda e_{0}-1}\bigg]
\frac{1}{(\Lambda e_{0}-1)^{n-1-t}}
e_{0}^{n-2-t}e_{\xi'}w' = \Lambda \sum_{t=0}^{n-2} (e_{0}^{t}e_{\xi}w) 
\text{ }\sh\text{ }
\frac{1}{(\Lambda e_{0}-1)^{n-t}}
e_{0}^{n-1-t}e_{\xi'}w' $$
This and the third term of (\ref{eq:e0x}) are, respectively, the $0 \leqslant t\leqslant n-2$ terms and the $t=n-1$ term of $\Lambda R$. This proves (\ref{eq:equation 4.10}) which implies the result.
\newline (iii) By \cite{I-3}, equation (2.1), $\circ_{\har}^{\smallint_{1,0}}$ is characterized by an equation which we can rewrite with the notation of this paper as
$\comp^{\har,\Ad}(g \circ_{\Ad}^{\smallint_{1,0}} f) = g \circ_{\har}^{\smallint_{1,0}} \comp^{\har,\Ad}f$. This gives the result.
\end{proof}

\begin{Definition}
We call (\ref{eq: DS har int1,0 1}) the \emph{harmonic quasi-shuffle equation} and (\ref{eq: DS har int1,0 2}) the \emph{harmonic shuffle equation of the framework $\int_{1,0}$} ; and we call the collection of (\ref{eq: DS har int1,0 1}) and (\ref{eq: DS har int1,0 2}) the \emph{harmonic double shuffle equations} of the framework $\int_{1,0}$.
\newline Let $\DS_{\har}^{\smallint_{1,0}}$ be the affine ind-scheme defined by the equations (\ref{eq: DS har int1,0 1}) and (\ref{eq: DS har int1,0 2}), which we call the \emph{harmonic double shuffle equations} of the framework $\int_{1,0}$.
\end{Definition}

\subsubsection{In the framework $\smallint$}

We construct the harmonic double shuffle equations in the framework of $\smallint$, i.e. by considering power series expansions of multiple polylogarithms.

\begin{Proposition} \label{lemma shuffle rien} (i) For any words $w,\tilde{w}$, we have 
\begin{equation} \label{eq: DS har int 1} \har_{\mathcal{P}^{\mathbb{N}}}(w)\har_{\mathcal{P}^{\mathbb{N}}}(\tilde{w}) = \har_{\mathcal{P}^{\mathbb{N}}}(w \ast_{\har} \tilde{w}) .
\end{equation}
\noindent (ii) For any words $w= \big( (n_{i})_{d};(\xi_{i})_{d} \big)$, $\tilde{w} =\big( (\tilde{n}_{i})_{d};(\tilde{\xi}_{i})_{d}\big)$,
we have
\begin{multline} \label{eq: DS har int 2}
\har_{\mathcal{P}^{\mathbb{N}}}(w\text{ }\sh\text{ }\tilde{w}) = \har_{\mathcal{P}^{\mathbb{N}}}\big( (\shft_{\ast}S_{Y})(\tilde{w}) w\big)
\\ = \sum_{l_{1},\ldots,l_{d} \in \mathbb{N}} 
\prod_{i=1}^{d'} {-\tilde{n}_{i} \choose l_{i}} (-1)^{\tilde{n}_{i}}
\har_{\mathcal{P}^{\mathbb{N}}}\big(
n_{1},\ldots,n_{d-1},n_{d}+\tilde{n}_{d}+l_{d'},\tilde{n}_{d'-1}+l_{d'-1},\ldots,\tilde{n}_{1}+l_{1};
\xi_{1},\ldots,\xi_{d},\tilde{\xi}_{d'},\ldots,\tilde{\xi}_{1} \big) .
\end{multline}
\end{Proposition}

\begin{proof} (i) This amounts to the quasi-shuffle relation for prime weighted cyclotomic multiple harmonic sums, $\har_{p^{\alpha}}(w) \har_{p^{\alpha}}(w') \har_{p^{\alpha}}(w \ast w')$, which is a consequence of the known double shuffle relations for multiple polylogarithms in two variables.
\newline (ii) The shuffle equation for multiple polylogarithms in one variable gives $\Li[w\text{ }\sh\text{ }\tilde{w}] = \Li[w] \Li[\tilde{w}]$. Let $m \in \mathbb{N}^{\ast}$ ; by equation (\ref{eq:multiple polylogarithms power series expansion}) for all $m' \in \{1,\ldots,m-1\}$,
$$ \Li[\tilde{w}][z^{m-m'}] = \sum_{0<m_{1}<\ldots <m_{d} = m-m'} \frac{(\frac{\xi_{2}}{\xi_{1}})^{m_{1}} \ldots (\frac{1}{\xi_{d}})^{m_{d}}}{m_{1}^{n_{1}} \ldots m_{d}^{n_{d}}} = \sum_{m'=m'_{d}<\ldots<m'_{1}<m} 
\frac{(\frac{\xi_{2}}{\xi_{1}})^{m'-m'_{1}} \ldots (\frac{1}{\xi_{d}})^{m'-m'_{d}}}{(m'-m'_{d})^{n_{d}} \ldots (m'-m'_{1})^{n_{1}}} $$
and
\begin{equation*} \label{eq:intermediate}
(\Li[w]\Li[\tilde{w}])[z^{m}] = \sum_{0<m_{1}<\ldots<m_{d} < m' < m'_{d'}<\ldots<m'_{1}<m} 
\frac{\big( \frac{\xi_{2}}{\xi_{1}} \big)^{m_{1}} \ldots \big( \frac{\xi_{d+1}}{\xi_{d}} \big)^{m_{d}}
\big( \frac{\xi_{d'+1}}{\xi_{d+1}} \big)^{m'}
\big( \frac{\xi_{d'}}{\xi_{d+1}} \big)^{m'_{d'}} \ldots \big(\frac{\xi_{2}}{\xi_{1}} \big)^{m_{1}} }
{ m_{1}^{n_{1}}\ldots {m'}_{d}^{n_{d}} (m-{m'}_{d'})^{\tilde{n}_{d'}} \ldots (m-{m'}_{1})^{\tilde{n}_{1}}} .
\end{equation*}
Now assume that $m = p^{\alpha}$. For all $m'_{i} \in \{1,\ldots,p^{\alpha}-1\}$. We have $v_{p}(m'_{i})< v_{p}(p^{\alpha})$, thus  $\displaystyle\frac{1}{(m'_{i}-p^{\alpha})^{\tilde{n}_{i}}}=  m^{-\tilde{n}_{i}} \sum\limits_{l\geqslant  0} {-\tilde{n}_{i} \choose l} \big( \frac{p^{\alpha}}{m'_{i}} \big)^{l} \in \mathbb{Z}_{p}$, whence :
\begin{multline*}
(\Li[w]\Li[\tilde{w}])[z^{p^{\alpha}}] =
\\ 
\sum_{l_{1},\ldots,l_{d'} \in \mathbb{N}} \prod_{i=1}^{d} {-{n'}_{i} \choose l_{i}} (-1)^{t'_{i}} \frac{1}{z_{j_{1}}^{p^{\alpha}}} \sum_{\substack{0<m_{1} < \ldots < m_{d-1} < m_{d} = m \\ = m_{d} < m_{d-1} < \ldots m_{1} < p^{\alpha}}}
\frac{ (\frac{z_{i_{2}}}{z_{i_{1}}})^{n_{1}} \ldots (\frac{1}{z_{i_{d}}})^{n} y_{d'}^{n} (\frac{z_{j_{d'-1}}}{z_{j_{d'}}})^{m_{d-1}} \ldots (\frac{z_{j_{2}}}{z_{i_{1}}})^{n_{1}}}
{m_{1}^{n_{1}} \ldots l_{d}^{n_{d}+t_{d'}+l_{d'}} l_{d-1}^{n_{d-1}+l_{d-1}} \ldots l_{1}^{n_{1}+l_{1}}} .
\end{multline*}
On the other hand, $w\text{ }\sh\text{ } \tilde{w}$ is a linear combination of words whose rightmost letter is not $e_{0}$ and, by equation (\ref{eq:har et Li 2}), we have
$$ \sum_{0<m'<p^{\alpha}} \Li[w\text{ }\sh\text{ }\tilde{w}][m'] = \har_{p^{\alpha}}[w\text{ }\sh\text{ }\tilde{w}] . $$
\noindent We deduce the result, given the definition of $\shft_{\ast}$ and $S_{Y}$ (Definition \ref{def SY}).
\end{proof}

\begin{Definition}
We call (\ref{eq: DS har int 1}) the
\emph{harmonic quasi-shuffle equation}
and (\ref{eq: DS har int 2}) the  \emph{harmonic shuffle equation of the framework $\int$} ; and we call the collection of (\ref{eq: DS har int 1}) and (\ref{eq: DS har int 2}) the \emph{harmonic double shuffle equations} of the framework $\int$.
\newline Let $\DS_{\har}^{\smallint}$ be the affine ind-scheme defined by the equations (\ref{eq: DS har int 1}) and (\ref{eq: DS har int 2}), which we call the \emph{harmonic double shuffle equations} of the framework $\int$.
\end{Definition}

\subsubsection{In the framework $\Sigma$}

Given that the power series expansion of multiple polylogarithms is written in terms of sums of series, we get directly the same harmonic double shuffle equations in the framework $\Sigma$. Indeed, it is known that the integral shuffle relation of cyclotomic multiple zeta values can be understood purely in terms of their formula as iterated series : it is a generalization of a proof which goes back to Euler, who defined multiple zeta values in depth $1$ and $2$ by the series formula in (\ref{eq:multizetas}), and found the double shuffle relations in this particular case. The same proof gives the shuffle relation for multiple harmonic sums : for $w=\big((n_{i})_{d};(\xi_{i})_{d+1}\big)$, $w'=\big(({n'}_{i})_{d};({\xi'}_{i})_{d+1}\big)$,
$$ \har_{m}(w\text{ }\sh\text{ }w') = \sum_{0<m_{1}<\ldots<m_{d} < m' < m'_{d'}<\ldots<m'_{1}<m} 
\frac{\big( \frac{\xi_{2}}{\xi_{1}} \big)^{m_{1}} \ldots \big( \frac{\xi_{d+1}}{\xi_{d}} \big)^{m_{d}}
\big( \frac{\xi_{d'+1}}{\xi_{d+1}} \big)^{m'}
\big( \frac{\xi_{d'}}{\xi_{d+1}} \big)^{m'_{d'}} \ldots \big(\frac{\xi_{2}}{\xi_{1}} \big)^{m_{1}} }
{ m_{1}^{n_{1}}\ldots {m'}_{d}^{n_{d}} (m-{m'}_{d'})^{\tilde{n}_{d'}} \ldots (m-{m'}_{1})^{\tilde{n}_{1}}} , $$
from which one can deduce the harmonic shuffle relation defined in \S3.3.2.
As concerns the quasi-shuffle relation for cyclotomic multiple harmonic values, it is of course immediate in terms of series. Thus we can directly define : 

\begin{Definition} Let $\DS_{\har}^{\Sigma} = \DS_{\har}^{\smallint}$.
\end{Definition}

\subsubsection{Comparison between the results of the three frameworks}

We now show that the definitions of \S3.3.1 and \S3.3.2 are equivalent.

\begin{Lemma} \label{lemma for comparison of ds} Let $F$ a function $\mathcal{O}^{\ast} \rightarrow K$ and $\tilde{\imath}$, a function $\mathcal{O}^{\ast} \rightarrow \mathcal{O}^{\ast}[[\Lambda]]$ satisfying, for all $a,b$ words in $\mathcal{O}^{\ast}$, $\tilde{\imath}(ab)= \tilde{\imath}(b)\tilde{\imath}(a)$. We have an equivalence between :
\newline (i) $\forall n \in \mathbb{N}^{\ast}$, $\forall \xi \in \mu_{N}(K)$, $\forall w,w'$ words in $\mathcal{O}^{\ast}$, $F\big((e_{0}^{n-1}e_{\xi}w) \sh w'\big) = F\big(w \sh (\tilde{\imath}(e_{0}^{n-1}e_{\xi})w')\big)$
\newline (ii) $\forall u,w,w'$ words in $\mathcal{O}^{\ast}$, $F\big((uw) \sh w'\big) = F\big(w \sh (\tilde{\imath}(u)w')\big)$
\newline (iii) $\forall w,w'$ words in $\mathcal{O}^{\ast}$, $F(w\text{ }\sh\text{ }w') = F(\tilde{\imath}(w) w')$.
\end{Lemma}

\begin{proof} (i) $\Rightarrow$ (ii) : we write $u$ as a concatenation of words of the form $e_{0}^{n_{i}-1}e_{1}$ and we iterate (iii).
\newline (ii) $\Rightarrow$ (iii) : we take $w = \emptyset$. 
\newline (iii) $\Rightarrow$ (i) : we apply (iii) to each member of (i).
\end{proof}

\begin{Proposition} We have $\DS_{\har}^{\smallint_{1,0}} = \DS_{\har}^{\smallint}$.
\end{Proposition} 

\begin{proof} The harmonic quasi-shuffle relations in $ \DS_{\har_{\mathcal{P}^{\mathbb{N}^{\ast}}}}^{\smallint}$ and $ \DS_{\har_{\mathcal{P}^{\mathbb{N}^{\ast}}}}^{\smallint_{1,0}}$, equations (\ref{eq: DS har int1,0 1}) and (\ref{eq: DS har int 1}), are identical. The equivalence between the harmonic shuffle equations (\ref{eq: DS har int1,0 2}) and (\ref{eq: DS har int 2}) follows from Lemma \ref{lemma for comparison of ds}.
\end{proof}

\begin{Definition} Let us call denote by 
	$\DS_{\har}$ the common value of $\DS_{\har}^{\smallint_{1,0}}$, $\DS_{\har}^{\smallint}$ and $\DS_{\har}^{\Sigma}$ and call its equations the prime harmonic double shuffle equations.
\end{Definition}

\begin{Remark} We see that comparing the results in the frameworks $\int_{1,0}$ and $\int$ requires a proof, whereas comparing the results in the frameworks $\int$ and $\Sigma$ is trivial, and this is because the power series expansion of multiple polylogarithms is expressed in terms of series. By contrast, in \cite{I-2} and \cite{I-3}, comparing the pro-unipotent harmonic actions $\circ_{\har}^{\smallint}$ and $\circ_{\har}^{\Sigma}$ required a proof, whereas comparing the pro-unipotent harmonic actions $\circ_{\har}^{\smallint_{1,0}}$ and $\circ_{\har}^{\smallint}$ was simple.
\end{Remark}

We now write the "overconvergent" variants of the previous results.

\begin{Proposition} \label{prop other view on }
	(i) ($\int_{1,0}$) $\Lambda$Ad$p$MZV$\mu_{N}^{\dagger}$'s satisfy equations obtained as remainders (in the sense of power series in $\Lambda$) of the equations of $\Lambda$Ad$p$MZV$\mu_{N}$'s. By taking $\Lambda=1$, we deduce equations satisfied by MHV$\mu_{N}^{\dagger}$'s.
	\newline (ii) ($\int$) (a) For all words $w,w'$ on $e_{0\cup \mu_{N}}$ we have :
	$$ \har_{p^{\alpha}}^{\dagger_{p,\alpha}}(w \text{ }\sh\text{ }w') = \har_{p^{\alpha}}^{\dagger_{p,\alpha}}(w) +
	\har_{p^{\alpha}}^{\dagger_{p,\alpha}}(w') + \har_{p^{\alpha}}(\shft_{\ast}(S_{Y}(w))w') . $$
	(b) For any automorphism $\sigma$ of $\mathbb{P}^{1} - \{0,\mu_{N},\infty\}$ which fixes $0$, we have an equation satisfied by $\har_{p^{\alpha}}^{\dagger_{p,\alpha}}$ obtained by 
	$$ \Li_{p,\alpha}^{\dagger}(\sigma(z)) = \sigma_{\ast} \har_{p^{\alpha}}^{\dagger_{p,\alpha}}(z) . $$  
\end{Proposition}

\begin{proof} (i) This is immediate by Definition \ref{def over adjoint}. 
	\newline (ii) (a) This follows from the shuffle relation $\Li_{p,\alpha}^{\dagger}[w\text{ }\sh\text{ }w'] = \Li_{p,\alpha}^{\dagger}[w]\Li_{p,\alpha}^{\dagger}[w']$ specialized to the coefficient $[z^{p^{\alpha}}]$, combined with the differential equation (\ref{eq:horizontality equation}) characterizing $\Li_{p,\alpha}^{\dagger}$, Lemma \ref{lemma power series expansion Li dagger}, and the proof of Lemma \ref{lemma shuffle rien}.
	\newline (b) This follows from the definition of $\Li_{p,\alpha}^{\dagger}$ in terms of the value of the Frobenius on the canonical de Rham path (\cite{I-1}, \S1) and the functoriality of the Frobenius.
\end{proof}

We see that, as in \S3.3, the formulas for the harmonic shuffle relation in the framework $\int_{1,0}$ and in the framework $\int$ are different. We leave to the reader to check that they are equivalent, following \S3.3.3. Using the framework $\Sigma$ here is beyond the scope of this paper : this requires to uses the formulas for Ad$p$MZV$\mu_{N}$'s found in \cite{I-2}, and this will be the subject of \cite{II-2}.

\begin{Definition} Let $\DS_{\har}^{\dagger}$, resp. $\M_{\har}^{\dagger}$ be the affine ind-scheme defined by the equations obtained in Proposition \ref{prop other view on } as variants of those of $\DS_{\har}$, resp. $\M_{\har}$ defined in \S3 and \S4 respectively.
\end{Definition}

We have canonical isomorphisms
$\DS_{\har}^{\dagger} \simeq \DS_{\har}$, resp. $\M_{\har}^{\dagger} \simeq  \M_{\har}$.

\begin{Remark} An extension of Proposition \ref{prop other view on } (ii), which would include a quasi-shuffle equation for $\har_{p^{\alpha}}^{\dagger_{p,\alpha}}$ and the equations obtained by the functoriality with respect to automorphisms of 
	$\mathcal{M}_{0,5}^{(N)}$, could be obtained by using the variant of $\Li_{p,\alpha}^{\dagger}$ on  $\mathcal{M}_{0,5}^{(N)}$, defined by the Frobenius of  $\pi_{1}^{\un,\crys}(\mathcal{M}_{0,5}^{(N)})$. Here, in simplicial coordinates, $\mathcal{M}_{0,5}^{(N)}$ is 
	$\{ (y_{1},y_{2},\ldots,y_{n}) \in (\mathbb{P}^{1} - \{0,1,\infty\})^{n} \text{ }|\text{ }\forall i,j, \forall \xi \text{ s.t. }\xi^{N}=1, y_{i} \not= \xi y_{j} \}$. We leave it to the reader.
\end{Remark}

\subsection{A corollary of the shuffle equation : the (usual, adjoint and harmonic) reversal equations\label{paragraph reflexion}}

Let $f \in \Pi_{1,0}(K)$ ; since $f$ satisfies the shuffle equation, denoting by $S$ the antipode of the shuffle Hopf algebra, we have $\hat{S}^{\vee}(f) = f^{-1}$. By writing $f^{-1} = (1-(1-f))^{-1} = \sum\limits_{l\geqslant 0} (1-f)^{l}$, we deduce polynomial equations on the coefficients of $f$ :

\begin{equation} \label{eq: reflexion}
\text{for all words }w,\text{ } 
f[S(w)] = \sum_{\substack{l\geqslant  0 \\ w_{1}\ldots w_{l}=w}} (-1)^{l}\prod_{i=1}^{l}f[w_{i}] .
\end{equation}

\begin{Definition} \label{def reflexion} We call (\ref{eq: reflexion}) the reversal equation.
\end{Definition}

Applying it to $f=\Phi_{\KZ}$ resp. $f=\Phi_{p,\alpha}$ gives a family of polynomial equation on MZV$\mu_{N}$'s, resp. $p$MZV$\mu_{N}$'s. Our terminology ``reversal equation'' in Definition \ref{def reflexion} is motivated by Rosen's ``asymptotic reversal theorem'', \cite{Rosen} which we are going to recover below, as a particular case of Proposition \ref{harmonic reflexion}, as particular cases of our results, and interpret in terms of $\pi_{1}^{\un}(\mathbb{P}^{1} - \{0,1,\infty\})$. We are going to see that the reversal equation has adjoint and harmonic analogues which are given by relatively simple formulas and which are quite natural.

\begin{Proposition} \label{reflexion adjoint} Let $\psi \in K\langle \langle e_{0\cup \mu_{N}}\rangle\rangle$ satisfying the shuffle equation modulo products (equation (\ref{eq:ds adjoint 1})). We have :
\begin{multline} \label{eq:reflexion adjoint}
\text{for any positive integers }
n_{i}, n'_{i'} \text{ and N-th roots of unity } \xi_{i},\xi'_{i},
\\
\sum_{\substack{l \geqslant 0 \\
		l_{i} \geqslant 0 (1 \leqslant i \leqslant d)\\ 
	l+l_{1}+\ldots+l_{d}=L}}
\prod_{i=1}^{d} {-n_{i} \choose l_{i}} \psi\big( \big( (n'_{i})_{d'},(n_{i}+l_{i})_{d};({\xi'}_{i})_{d'},(\xi_{i})_{d},\xi \big),l\big) 
\\ = \sum_{\substack{l' \geqslant 0 \\
		{l'}_{i} \geqslant 0 (1 \leqslant i \leqslant d') \\ 
	l'+{l'}_{1}+\ldots+{l'}_{d'}=L}}
\prod_{i=1}^{d'} {-{n'}_{i} \choose {l'}_{i}}
\psi\big( \big( (n_{i})_{d'},({n'}_{i}+{l'}_{i})_{d};\xi,({\xi}_{i})_{d'},({\xi'}_{i})_{d}\big),l'\big) .
\end{multline}
\end{Proposition}

\begin{proof} The result amounts to the following equality :
\begin{multline} \label{eq: 4 2 5} (-1)^{\sum\limits_{i=1}^{d}n_{i}} \psi[ \frac{1}{1-\Lambda e_{0}}e_{\xi} \frac{e_{0}^{n_{1}-1}}{(1-\Lambda e_{0})^{n_{1}}}e_{\xi_{1}} \ldots
\frac{e_{0}^{n_{d}-1}}{(1-\Lambda e_{0})^{n_{d}}}e_{\xi_{d}}e_{0}^{{n'}_{d'}-1}e_{{\xi'}_{d'}}\ldots e_{0}^{{n'}_{1}-1}e_{{\xi'}_{1}}]  \\ =(-1)^{\sum\limits_{i'=1}^{d'}{n'}_{i'}} 
\psi \big[\frac{1}{1-\Lambda e_{0}}e_{\xi_{1}} \ldots \frac{e_{0}^{{n'}_{1}}-1}{(1-\Lambda e_{0})^{n_{1}}}e_{\xi_{d'}}
\frac{e_{0}^{{n'}_{d'}}-1}{(1-\Lambda e_{0})^{{n'}_{d'}}}e_{\xi_{d}}e_{0}^{n_{d}-1}\ldots e_{\xi_{1}}e_{0}^{n_{1}-1}e_{\xi} \big] .
\end{multline}
(a) By the hypothesis, $\psi$ is primitive for the shuffle coproduct. Thus, let $S$ be the antipode of $\mathcal{O}^{\sh}$ ; we have $\hat{S}^{\vee}(\psi)=-\psi$ ; and the first line of (\ref{eq: 4 2 5}) is equal to
\begin{equation} \label{eq:this} (-1)^{\sum\limits_{j=1}^{d'}{n'}_{i}} \psi\big[e_{{\xi'}_{1}}e_{0}^{{n'}_{1}-1}\ldots e_{{\xi'}_{d'}}e_{0}^{{n'}_{d'}-1}e_{{\xi}_{d}}
\frac{e_{0}^{n_{d}-1}}{(1+\Lambda e_{0})^{n_{d}}} \ldots e_{\xi_{1}} \frac{e_{0}^{n_{1}-1}}{(1+\Lambda e_{0})^{n_{1}}}e_{\xi}\frac{1}{1+\Lambda e_{0}}\big] .
\end{equation}
(b) Let  $f \in K \langle \langle e_{0\cup \mu_{N}} \rangle\rangle$ such that, for all words $w$, we have $f[w\text{ }\sh\text{ }e_{0}] = f[w]f[e_{0}]$ ; let $T,U_{1},\ldots,U_{d}$ formal variables ; we have :
\begin{multline}
\label{eq: shuffle coefficients e0 TU}
f[\frac{e_{0}^{{n''}_{d}-1}}{(1-U_{d}e_{0})^{{n''}_{d}}}e_{{\xi''}_{d}} \ldots
\frac{e_{0}^{{n''}_{1}-1}}{(1-U_{1}e_{0})^{{n''}_{1}}}e_{{\xi''}_{1}} \frac{1}{1 - Te_{0}}]
\\ = 
f[\frac{e_{0}^{{n''}_{d''}-1}}{(1 - (U_{d}-T)e_{0})^{{n''}_{d''}}}e_{{\xi''}_{d''}}\ldots
\frac{e_{0}^{{n''}_{1}-1}}{(1 - (U_{1}-T)e_{0})^{{n''}_{1}}}
e_{{\xi''}_{1}}] e^{f[e_{0}]T} .
\end{multline}
We apply this to the particular case where $f[e_{0}]=0$ and where the first line of (\ref{eq: shuffle coefficients e0 TU}) is (\ref{eq:this}). In that case, the second line of (\ref{eq: shuffle coefficients e0 TU}) becomes equal to the second line of (\ref{eq: 4 2 5}).
\end{proof}

\begin{Definition} We call (\ref{eq:reflexion adjoint}) the adjoint reversal equation.
\end{Definition}

\begin{Proposition} \label{harmonic reflexion}We can prove by the three frameworks $\int_{1,0}$, $\int$ and $\Sigma$ that : 
\begin{equation} \label{eq:first eq of corollary}
\text{for all words }w,w',\text{ }
\text{har}_{\mathcal{P}^{\mathbb{N}}}( (\shft_{\ast}S_{Y})(w')w) = \text{har}_{\mathcal{P}^{\mathbb{N}}}( (\shft_{\ast}S_{Y})(w)w') .
\end{equation}
\end{Proposition}

\begin{proof} (i) In the framework $\int_{1,0}$ : we take $f=\sum\limits_{\xi \in \mu_{N}(K)} \xi^{-p^{\alpha}} {\Phi_{p,\alpha}^{(\xi)}}^{-1}e_{\xi}\Phi_{p,\alpha}^{(\xi)}$ and $\Lambda=1$ in equation (\ref{eq: 4 2 5}) and we use equation (\ref{eq:formula for n=1}).
\newline (ii) In the framework $\smallint$ : this is a direct consequence of the harmonic shuffle equation (equation \ref{eq: DS har int 2}).
\newline (iii) In the framework $\Sigma$ : this is an easy generalization of Rosen's proof of the ``asymptotic reversal theorem'' \cite{Rosen} which corresponds to the $N=1$, $\alpha=1$ and $w'=\emptyset$ case.
\end{proof}

\begin{Definition} We call (\ref{eq:first eq of corollary}) the harmonic reversal equation.
\end{Definition}

We also note that, by choosing $w'$ to be the empty word in (\ref{eq:the last equation}) is deduced by choosing $w'$ equal to the empty word in (\ref{eq:first eq of corollary}), we deduce the simpler equation
\begin{equation} \label{eq:the last equation}
\text{ for all words }w,\text{ }
\text{har}_{p^{\alpha}}((\shft_{\ast}S_{Y})(w)) = \text{har}_{p^{\alpha}}(w) .
\end{equation}

\section{Around associator equations and Kashiwara-Vergne equations\label{double shuffle}}

We review briefly associator equations and Kashiwara-Vergne equations and the known relation between them (\S4.1) then we explain that Kashiwara-Vergne equations can be formulated as a property of adjoint MZV's more naturally than MZV's (\S4.2) and we explain equations satisfied by multiple harmonic values, and more generally harmonic multiple polylogarithms, which are related to Kashiwara-Vergne equations (\S4.3). In this section, we take $N=1$ most of the time, because to our knowledge there is no cyclotomic Kashiwara-Vergne theory.

\subsection{Review on associators and Kashiwara-Vergne equations}

\subsubsection{Associators}

The notion of associators has been introduced in \cite{Drinfeld}. Let $k$ be a field of characteristic $0$ and let $\mu \in k$. By \cite{Furusho pentagon}, the set of associators $M_{\mu}(k)$ is the set of elements $\Phi \in \Pi(k)$ (we are using Notation \ref{la premiere notation}, and assuming $N=1$), such that $\mu=\pm \sqrt{24\Phi[e_{0}e_{1}]}$ and
$$ \phi(e_{12},e_{23}+e_{24})\phi(e_{13}+e_{23},e_{34}) = \phi(e_{23},e_{34}) \phi(e_{12}+e_{13},e_{24}+e_{34}) 
\phi(e_{12},e_{23}) $$
where $e_{ij}$ are the generators of $\Lie(\pi_{1}^{\un,\dR}(\mathcal{M}_{0,5},\omega_{\dR}))$ (\S2.1.3). The definition in \cite{Drinfeld} uses several equations : (2.12), (2.13) and a rescaled version of equation (5.3) of \cite{Drinfeld}.
\newline\indent On the other hand one has the scheme $\GRT_{1}$ defined in \cite{Drinfeld}, (equations (5.12), (5.13), (5.14), (5.15) of \cite{Drinfeld}) ; by \cite{Drinfeld}, Proposition 5.9, $\GRT_{1}$ is isomorphic to $M_{0}$. We note that equation (5.15) of \cite{Drinfeld} is 
\begin{equation} \label{eq:sum of residues} e_{0} + \phi(e_{0},e_{1})^{-1}e_{1}\phi(e_{0},e_{1}) + \phi(e_{0},e_{\infty})^{-1}e_{\infty}\phi(e_{0},e_{\infty}) = 0 
\end{equation}
where $e_{0}+e_{1}+e_{\infty}=0$ : the $e_{x}$'s generate $\Lie(\Pi)$.
\newline\indent The Ihara product (\ref{eq:Ihara}) restricts to a group law on $\GRT_{1}$ (\cite{Drinfeld}, equation (5.16)) and to an action of $\GRT_{1}$ on $M_{\mu}$ which makes $M_{\mu}$ into a $\GRT_{1}$-torsor (\cite{Drinfeld}, Proposition 5.5).
\newline\indent The generating series of multiple zeta values, $\Phi_{\KZ} \in \Pi_{1,0}(\mathbb{R})$ (equation (\ref{eq:Phi KZ})) has been first introduced in \cite{Drinfeld} \S2, and we have $\Phi_{\KZ} \in M_{2i\pi}(\mathbb{R})$ by \cite{Drinfeld}. By \cite{Unver Drinfeld}, combined with the fact that $\Phi_{p,1}$ is in the commutator subgroup of $\Pi_{1,0}(\mathbb{Q}_{p})$, proved in \cite{Furusho 2}, \S3, we have $\GRT_{1}(\mathbb{Q}_{p})$ ; this implies that $\Phi_{p}^{\KZ} \in \GRT_{1}(\mathbb{Q}_{p})$ (\cite{Furusho 2}, proof of Proposition 3.1). By the relations of iteration of the Frobenius (\cite{I-3}, equations (1.11), (1.12), (1.13) and Proposition 1.5.2), which involve the Ihara product, this implies that $\Phi_{p,\alpha} \in \GRT_{1}(\mathbb{Q}_{p})$ for all $\alpha \in \mathbb{Z} \cup \{\pm \infty\} - \{0\}$.

\subsubsection{Review on Kashiwara-Vergne equations, according to Alekseev, Enriquez and Torossian \cite{AT}, \cite{AET}}

Let $k$ be a field of characteristic $0$.
Let $\lie_{n}$ be the free Lie algebra over $k$ on $n$ variables $x_{1},\ldots,x_{n}$. Let $\widehat{\lie}_{n}$ be its degree completion (where the $x_{i}$'s have degree 1).
\newline\indent For any $u_{1},\ldots,u_{n}$ in $\widehat{\lie}_{n}$, let $[[u_{1},\ldots,u_{n}]]$ be the derivation of $\lie_{n}$ defined by $x_{i} \mapsto [x_{i},u_{i}]$ for all $i$. Such a derivation is called tangential and the set of tangential derivations is denoted by $\tder_{n}$.
\newline\indent For any $U_{1},\ldots,U_{n}$ in $\exp(\widehat{\lie}_{n})$, let $[[U_{1},\ldots,U_{n}]]$ be the automorphism of $\widehat{\lie}_{n}$ defined by $x_{i} \mapsto U_{i}x_{i}U_{i}^{-1}$ for all $i$. Such an automorphism is called tangential and the set of tangential automorphism is denoted by $\TAut_{n}$.
\newline\indent Let $A_{n}$ be the universal enveloping algebra of $f_{n}$ and $T_{n}=A_{n}/[A_{n},A_{n}]$. The image of an element $S$ by the map $A_{n} \rightarrow T_{n}$ is denoted by $\langle S \rangle$. Let $\partial_{k} : A_{n} \rightarrow A_{n}$ be defined by $x= x_{0}+\sum_{k=1}^{n} \partial_{k}(x)x_{k}$ et $x_{0} \in k$. Let $j : \text{tder}_{n} \rightarrow \hat{T}_{n}$, $[[u_{1},\ldots,u_{n}]] \mapsto \langle \sum_{k=1}^{n} x_{k} \partial_{k}(u_{k})\rangle$. It integrates into $J: \TAut_{n} \mapsto \hat{T}_{n}$ (\cite{AT}, Proposition 5.1). 
\newline\indent Below, $F_{2}(k)$ be the pro-unipotent completion of the free-group on two generators $X,Y$ and $\sim$ means ``is conjugated to''. Following \cite{AET} \S2.1, let
\begin{multline*} \Sol \KV = 
\\ \{ \mu : F_{2}(k) \simlra \Aut(\exp(\widehat{\lie}_{2})) \text{ }|\text{ }\mu(X) \sim e^{x},\mu(Y) \sim e^{y}, \mu(XY)=e^{x+y}, \exists r \in u^{2}k[[u]], j(a) = \langle r(x+y)-r(x)-r(y)\rangle \} 
\end{multline*}
where $r$ is uniquely determined by $\mu$. Following \cite{AET}, \S2.2, let
$$ \KRV = \{ a \in \Aut(\widehat{\lie}_{2}) \text{ }|\text{ }a(x) \sim x,a(y) \sim y, a(x+y)=x+y,
\\ \exists s \in u^{2}k[[u]], J(a) = \langle s(x+y)-s(x)-s(y)\rangle \} . $$
\indent By \cite{AET}, Theorem 2.1, the Ihara action (\ref{eq:Ihara}) makes $\KRV$ into a group and $\Sol \KV$ a torsor under that group, and there is a morphism of torsors $M_{1}(K) \rightarrow \Sol\KV(K)$, which sends $\Phi$ to 
\begin{equation} \label{eq:XY} 
\mu_{\Phi} : \begin{array}{l}
X \mapsto \Phi(x,-x-y)^{-1}e^{x}\Phi(x,-x-y) \\
Y \mapsto e^{-(x+y)/2}\Phi(y,-x-y)^{-1}e^{y} \Phi(y,-x-y)e^{(x+y)/2}
\end{array} .
\end{equation}
The element $r$ associated with $\mu_{\Phi}$ is equal to $-\log(\Gamma_{\Phi})$ where $\Gamma_{\Phi}= \exp \big(\sum\limits_{n\geq 2} \frac{(-1)^{n}}{n} \zeta_{\Phi}(n)u^{n} \big)$ (\cite{AET}, Proposition 2.2). We note that 
$\zeta_{\Phi}(2n) = -\frac{1}{2}\frac{B_{2n}}{(2n)!}$ for all $n$, which is independent of $\Phi$.

\subsection{The Kashiwara-Vergne equations as a property of adjoint multiple zeta values}

We restrict the study to $N=1$ because to our knowledge there is no cyclotomic Kashiwara-Vergne theory. There is a cyclotomic associator theory \cite{Enriquez}.

We are going to see that Kashiwara-Vergne equations can be naturally viewed as a property of adjoint MZV's rather than as a property of MZV's : this gives much simpler equations. This gives an analogy between the passage from associator equations to Kashiwara-Vergne equations constructed in \cite{AET} and \cite{AT} and the passage from double shuffle equations to adjoint double shuffle equations that we have constructed in \S3.

In particular, our question \ref{question} is an analogue of a conjecture of Alekseev-Torossian \cite{AT} which compares Drinfeld associators and solutions to the Kashiwara-Vergne problem.

\subsubsection{In the $p$-adic case}

We consider that the Kashiwara-Vergne equations are the adjoint version of associator equations. Let us rewrite them - or rather the equations of KRV - in terms of Ad$p$MZV's. By the above discussion, we consider the following equation :
\begin{equation} \label{eq:KV} \langle j(\mu_{\Phi_{p,\alpha}})\rangle = \langle \log(\Gamma_{\Phi_{p,\alpha}}(x)) + \log(\Gamma_{\Phi_{p,\alpha}}(y)) - \log(\Gamma_{\Phi_{p,\alpha}})(x+y) \rangle .
\end{equation}

Moreover, by the isomorphism $\lie_{2} \simeq \Lie \Pi_{0,0}$ given by $(x,y) \leftrightarrow (e_{1},e_{\infty})$ where $e_{0}+e_{1}+e_{\infty}=0$,
$\mu_{\Phi_{p,\alpha}}$ sends $\left\{ \begin{array}{l} 
e_{0} \mapsto e_{0}
\\ e_{1} \mapsto \Phi_{p,\alpha}(e_{0},e_{1})^{-1} e_{1} \Phi_{p,\alpha}(e_{0},e_{1}) 
\\ e_{\infty} \mapsto \Phi_{p,\alpha}(e_{0},e_{\infty})^{-1}  e_{\infty} \Phi_{p,\alpha}(e_{0},e_{\infty}) \end{array} \right.$. This recovers equation (\ref{eq:sum of residues}).
\newline This is obtained by applying 
$\underset{\mu \rightarrow 0}{\lim} \frac{1}{\mu} \frac{d}{d\mu}$ to the automorphism
$\left\{ \begin{array}{l} e^{e_{1}} \mapsto \Phi(e_{0},e_{1})^{-1} e^{\mu e_{1}} \Phi(e_{0},e_{1})
\\ e^{e_{\infty}} \mapsto  e^{\frac{\mu}{2}e_{0}} \Phi(e_{0},e_{\infty})^{-1} e^{\mu e_{\infty}} \Phi(e_{0},e_{\infty}) e^{-\frac{\mu}{2} e_{0}}
\end{array} \right.$ obtained by rescaling  (\ref{eq:XY}) : this is the variant of $\Sol KV$ and of the map $M_{1} \rightarrow \Sol \KV$ obtained by choosing $M_{\mu}$ instead of $M_{1}$.

\begin{Proposition} \label{KV explicit}The equation (\ref{eq:KV}) amounts to explicit linear equations on adjoint $p$MZV's.
\end{Proposition}

\begin{proof} Equation (\ref{eq:KV}) amounts to say that 
$$ j(\mu_{\Phi}) - \log(\Gamma_{\Phi}(e_{1})) - \log(\Gamma_{\Phi}(e_{\infty})) + \log(\Gamma_{\Phi})(e_{1}+e_{\infty}) \rangle $$
has image zero by the quotient map $k\langle\langle e_{1},e_{\infty}\rangle\rangle \rightarrow k
\langle\langle e_{1},e_{\infty} \rangle\rangle / [k
\langle\langle e_{1},e_{\infty} \rangle\rangle,k
\langle\langle e_{1},e_{\infty} \rangle\rangle]$.
\newline\indent Let $\tilde{\partial}_{e_{1}}^{(e_{1},e_{\infty})}, \tilde{\partial}_{e_{\infty}}^{(e_{1},e_{\infty})} : 
k\langle \langle e_{1},e_{\infty}\rangle\rangle \rightarrow k\langle\langle e_{1},e_{\infty} \rangle\rangle$ be defined, by, for any word $w$ on the alphabet $e_{1},e_{\infty}$,
$w = \partial_{e_{1}}^{(e_{1},e_{\infty})}(w)e_{1} + 
\partial_{e_{\infty}}^{(e_{1},e_{\infty})}(w)e_{\infty}$.
\newline One can compute $j(\mu_{\Phi})$ as follows. We have 
\begin{equation} \label{eq:jmuPhi} j(\mu_{\Phi}) = e_{1} \tilde{\partial}_{e_{1}}^{(e_{1},e_{\infty})} \big(\Phi^{-1}e_{1}\Phi\big)(e_{0},e_{1}) + e_{\infty}\tilde{\partial}_{e_{\infty}}^{(e_{1},e_{\infty})}\big(\Phi^{-1}e_{1}\Phi\big)(e_{0},e_{\infty}) 
\end{equation}
and the right-hand side of (\ref{eq:jmuPhi}) can be expressed using what follows. Let $w(e_{0},e_{1}) = e_{0}^{n_{d}-1}e_{1}\ldots e_{1}e_{0}^{n_{0}-1}$.
\newline\indent
If $n_{0}\geq 2$ we have 
$\left\{
\begin{array}{l} \tilde{\partial}_{e_{1}}^{(e_{1},e_{\infty})}(w(e_{1},e_{\infty})) =  e_{0}^{n_{d}-1}e_{1}\ldots e_{1}e_{0}^{n_{0}-1}e_{1}(-e_{1}) 
\\ \tilde{\partial}_{e_{\infty}}^{(e_{1},e_{\infty})}(w(e_{1},e_{\infty})) =  e_{0}^{n_{d}-1}e_{1}\ldots e_{1}e_{0}^{n_{0}-1}e_{1}(-e_{\infty}) 
\end{array} \right.$. 
\newline\indent 
If $n_{0}=1$, we have 
$\left\{
\begin{array}{l} \tilde{\partial}_{e_{1}}^{(e_{1},e_{\infty})}(w(e_{1},e_{\infty})) =  e_{0}^{n_{d}-1}e_{1}\ldots e_{1}e_{0}^{n_{0}-1} 
\\ \tilde{\partial}_{e_{\infty}}^{(e_{1},e_{\infty})}(w(e_{1},e_{\infty})) = 0 \end{array} \right. $
\newline\indent On the other hand, $\log(\Gamma_{\Phi}(e_{1})) + \log(\Gamma_{\Phi}(e_{\infty})) - \log(\Gamma_{\Phi})(e_{1}+e_{\infty}) \rangle$ is equal to :
$$ \sum_{n\geqslant 2} \frac{(-1)^{n}}{n}\zeta_{\Phi}(n) ( e_{1}^{n} + e_{\infty}^{n}-(e_{1}+e_{\infty})^{n})  = \sum_{n\geqslant 2} \frac{(-1)^{n}\zeta_{\Phi}(n)}{n}  e_{1}^{n} +
\sum_{n\geqslant 2} \frac{\zeta_{\Phi}(n)}{n} \big( (e_{0}+e_{1})^{n} - e_{0}^{n} \big) . $$
\indent We have a basis of $k\langle \langle e_{1},e_{\infty}\rangle\rangle$ formed by the words on the alphabet $\{e_{1},e_{\infty}\}$. It is sent to a generating family of $k\langle \langle e_{1},e_{\infty}\rangle\rangle/[k\langle \langle e_{1},e_{\infty}\rangle\rangle,k\langle \langle e_{1},e_{\infty}\rangle\rangle]$ from which we can extract a basis. A simple computation gives a partition of the set of words on $\{e_{1},e_{\infty}\}$ according to their image in $k\langle \langle e_{1},e_{\infty}\rangle\rangle/[k\langle \langle e_{1},e_{\infty}\rangle\rangle,k\langle \langle e_{1},e_{\infty}\rangle\rangle]$.
\newline\indent Finally we use the isomorphism $k\langle \langle e_{0},e_{1}\rangle\rangle \simeq k\langle \langle e_{1},e_{\infty}\rangle\rangle$, $f(e_{0},e_{1}) \mapsto f(e_{1},e_{\infty})$ defined by $e_{\infty}=-e_{0}-e_{1}$.
\end{proof}

\begin{Remark} Equation (\ref{eq:sum of residues}) can be regarded as a part of the Kashiwara-Vergne equations. It also amounts to an equation on $p$MZV's : for all words,
\begin{equation} \label{eq:KV dim 1} \forall w,\text{ } \zeta_{p,\alpha}^{\Ad}(w) + \zeta_{p,\alpha}^{\Ad}(w(e_{0}-e_{1},-e_{1})) = 0 .
\end{equation}
Indeed, we have $(\phi^{-1}e_{1}\phi)(e_{0},e_{\infty}) = \sum\limits_{w\text{ word on }\{e_{0},e_{1}\}} 
(\phi^{-1}e_{1}\phi)[w(e_{0}-e_{1},-e_{1})]w$.
\end{Remark}

\subsubsection{In the complex case}

The automorphism of $\Pi_{0,0}(\mathbb{C})$ which sends $(e^{e_{0}},e^{e_{1}}) \mapsto (e^{2i\pi e_{0}}, \Phi_{\KZ}^{-1}e^{2i\pi e_{1}}\Phi_{\KZ})$ satisfies the Kashiwara-Vergne equations rescaled by $\tau(2\pi i)$ where $\tau$ is as in equation (\ref{eq:tau}). This can be easily generalized to any $N$, using \cite{AKKN}. The computations of \S4 can be repeated with the equations instead of the equations of the group $\KRV$. As concerns the equations of dimension 1, the equation
$e_{0} + (\Phi_{\KZ}^{-1}e_{1}\Phi_{\KZ})(e_{0},e_{1}) + (\Phi_{\KZ}^{-1}e_{1}\Phi_{\KZ})(e_{0},e_{\infty}) \equiv 0 \mod (\zeta(2))$ is
$$ e^{2i\pi e_{0}} . \big( \Phi(e_{0},e_{1})^{-1} e^{2i\pi e_{1}} \Phi(e_{0},e_{1}) \big). 
\big( e^{\frac{\mu}{2}e_{0}} \Phi(e_{0},e_{\infty})^{-1} e^{2i\pi e_{\infty}} \Phi(e_{0},e_{\infty}) e^{-\frac{\mu}{2} e_{0}} \big) = 1 . $$
We leave the details to the reader.

\subsection{Related properties of harmonic multiple polylogarithms}

We find properties of multiple harmonic values and harmonic multiple polylogarithms, in the three frameworks $\int_{1,0}$, $\int$ and $\Sigma$, which are related to the previosu considerations. We prove that the equations arising from $\pi_{1}^{\un}(\mathbb{P}^{1} - \{0,1,\infty\})$ in these three frameworks are equivalent.

\subsubsection{In the framework $\int_{1,0}$}

\begin{Proposition} \label{harmonic duality DR} 
(i) $\Phi^{-1}e_{1}\Phi$ satisfies equation (\ref{eq:KV dim 1}) if and only if
$h(w) = \comp^{\Lambda \Ad,\Ad} (-1)^{\depth}\Phi^{-1}e_{1}\Phi$ satisfies
$$ \forall w,\text{ } h( w(e_{0}+e_{1},-e_{1})) = - \sum_{d'\geqslant  1, \text{ }z = e_{0}^{t_{d'}-1}e_{1}\ldots e_{0}^{t_{1}-1}e_{1}} (-1)^{d'} h(z.w) . $$
(ii) $\Phi^{-1}e_{1}\Phi$ satisfies equation (\ref{eq:KV}) if and only if $\comp^{\Lambda\Ad,\Ad}(\Phi^{-1}e_{1}\Phi)$ satisfies certain equations on $\Lambda$-adjoint $p$MZV's.
\end{Proposition}

Of course, in this statement, we can replace $\comp^{\Lambda \Ad,\Ad}$ by $\comp^{\har,\Ad}$, and $\Lambda$-adjoint $p$MZV's by MHV's.

\begin{proof} (i) We have, for all words $w$ : 
$(1-\Lambda(e_{0}+e_{1}))^{-1}e_{1}w = (1-\Lambda e_{0})^{-1}e_{1}w + (1-\Lambda e_{0})^{-1}\Lambda e_{1}(1-\Lambda (e_{0}+e_{1}))^{-1}e_{1}w$.
Indeed, we have $1 = (1 - \Lambda e_{0} - \Lambda e_{1} )(1 - \Lambda e_{0} - \Lambda e_{1} )^{-1} = 
(1 - \Lambda e_{0} )(1 - \Lambda e_{0} - \Lambda e_{1} )^{-1} - \Lambda e_{1} (1 - \Lambda e_{0} - \Lambda e_{1} )^{-1}$. Left multiplication by $(1 - \Lambda e_{0})^{-1}$ gives $(1-\Lambda(e_{0}+e_{1}))^{-1} = (1-\Lambda e_{0})^{-1} + (1-\Lambda e_{0})^{-1}\Lambda e_{1}(1-\Lambda (e_{0}+e_{1}))^{-1}$. Whence the equality. This implies the result.
\newline (ii) Follows from the definitions and from the translation of equation (\ref{eq:KV}) in the proof of Proposition \ref{KV explicit}.
\end{proof}

\begin{Definition} Let $\M_{\har}$ the ind-scheme defined by the equations among MHV's resp. $\Lambda$-adjoint $p$MZV's obtained in Proposition \ref{harmonic duality DR}.
\end{Definition}

\subsubsection{In the framework $\int$ \label{pre-associator paragraph}}

For any $n \geqslant 4$, any automorphism of $\mathcal{M}_{0,n}$ induces by functoriality an automorphism of  $\pi_{1}^{\un,\dR}(\mathcal{M}_{0,n})$ equipped with $\nabla_{\KZ}$. The functoriality of $\nabla_{\KZ}$ gives a relation of the type 
\begin{equation} \label{eq: LtildeCL} \tilde{L} = C L
\end{equation} 
where $L$ and $\tilde{L}$ are two different branchs of multiple polylogarithms on $\mathcal{M}_{0,n}$, defined as the unique solutions to $\nabla_{\KZ}$ with prescribed asymptotics at a chosen base-point, and $C \in \pi_{1}^{\un,\dR}(\mathcal{M}_{0,5},\omega_{\dR})$

\begin{Definition} Let us call pre-associator equations the equations of the form (\ref{eq: LtildeCL}).
\end{Definition} 

It is sufficient to restrict to $n \in \{4,5\}$.
The associator equations are deduced from the pre-associator equations by writing $C$ in terms of $\Phi_{\KZ}$ and by using that the automorphisms which are involved are of finite order.
\newline\indent Let us denote by $O$ the tangential base-point $(\vec{1}_{0},\vec{1}_{0})$ in cubic coordinates on $\overline{\mathcal{M}}_{0,5}$ as well as its image by $\overline{\mathcal{M}}_{0,5} \rightarrow \overline{\mathcal{M}}_{0,4}$, which we choose as the origin of the paths of integration. Let $\Stab_{O}^{\mathcal{M}_{0,4}}$ and $\Stab_{O}^{\mathcal{M}_{0,5}}$ be the stabilizers of $O$ in $\Aut(\mathcal{M}_{0,4}) = S_{3}$ and $\Aut(\mathcal{M}_{0,5})= S_{5}$. 
\newline\indent The associator equations (duality, hexagon and pentagon) can be obtained by the pre-associator equations associated with the automoprhisms $(z \mapsto 1-z)_{\ast}$, $(z \mapsto \frac{1}{z})_{\ast}$ resp. $\big(\sigma: (x_{1},x_{2},x_{3},x_{4},x_{5}) \mapsto ( x_{5},x_{4},x_{1},x_{3},x_{2})\big)_{\ast} = \big((c_{1},c_{2}) \mapsto (c_{2}, \frac{1 -c_{2}}{1 - c_{1}c_{2}})\big)_{\ast}$ in cubic coordinates, which are elements of $\Aut(\mathcal{M}_{0,4}) - \Stab_{O}^{\mathcal{M}_{0,4}}$, resp. $\Aut(\mathcal{M}_{0,4}) - \Stab_{O}^{\mathcal{M}_{0,5}}$.
\newline\indent We now write consequences of the pre-associator equations associated with elements of $\Stab_{O}^{\mathcal{M}_{0,4}}$ and $\Stab_{O}^{\mathcal{M}_{0,5}}$. We consider  the involution $\sigma : z \mapsto \frac{z}{z-1}$ which is in $\Aut(\mathcal{M}_{0,4})$, and the element of $\Aut(\mathcal{M}_{0,5})$ written in cubic coordinates as $\rho : (c_{1},c_{2}) \mapsto \big( -c_{1} \frac{1-c_{2}}{1-c_{1}}, - c_{2} \frac{1-c_{1}}{1-c_{2}}\big)$.

\begin{Proposition} \label{harmonic duality DR RT} (i) The pre-associator equation associated with $\sigma : z\mapsto \frac{z}{z-1}$ implies :
$$ \har_{\mathcal{P}^{\mathbb{N}}}( w(e_{0}+e_{1},-e_{1})) = - \sum_{\substack{d'\geqslant  1,\text{ }z = e_{0}^{t_{d'}-1}e_{1}\ldots e_{0}^{t_{1}-1}e_{1}}} (-1)^{\depth(z)} \har_{\mathcal{P}^{\mathbb{N}}}(z.w) . $$
\noindent (ii) \label{5 32}The pre-associator equation associated with $\rho : (c_{1},c_{2}) \mapsto \big( -c_{1} \frac{1-c_{2}}{1-c_{1}}, - c_{2} \frac{1-c_{1}}{1-c_{2}}\big)$ implies a family of equations on multiple harmonic values.
\end{Proposition}

\begin{proof} (i) See more generally the proof of Proposition \ref{proposition 5 30}, which we can apply to $(a,c,d)=(1,1,-1)$, knowing that $\displaystyle \sigma_{1,1,-1}^{\ast}\big( \frac{dz}{z} \big) = \frac{dz}{z} - \frac{dz}{z-1}$ and $\displaystyle\sigma_{1,1,-1}^{\ast} \big( \frac{dz}{z-1} \big) = \frac{dz}{z-1}$.
\newline (ii) We have, for any sequence of coefficients $(a_{m,n})$,
$$ \sum_{n,m\geq 0} a_{n,m}
\bigg(-\frac{c_{1}(1-c_{2})}{1-c_{1}}\bigg)^{n}
\bigg(-\frac{c_{2}(1-c_{1})}{1-c_{2}}\bigg)^{m} 
= \bigg( \sum_{n<m} + \sum_{n=m} + \sum_{n>m} \bigg) a_{n,m}
\bigg(-\frac{c_{1}(1-c_{2})}{1-c_{1}}\bigg)^{n}
\bigg(-\frac{c_{1}(1-c_{2})}{1-c_{1}}\bigg)^{m} . $$
The subsum $\displaystyle \sum\limits_{n=m}$ is
$\sum_{n \geq 0}  a_{n,n} \big(c_{1}c_{2}\big)^{n} $. 
The subsum $\displaystyle \sum\limits_{n<m}$ is, after writing $(n,m)=(n_{1},n_{1}+n_{2})$,
$$ \sum_{\substack{0\leq n_{1} \\ 1\leq n_{2}}} a_{n_{1},n_{1}+n_{2}} (-c_{1})^{n_{1}}(-c_{2})^{n_{1}+n_{2}} \big( \frac{1-c_{1}}{1-c_{2}} \big)^{n_{2}} = 
\sum_{\substack{0\leq n_{1}\\ 1\leq n_{2}}}\sum_{\substack{0\leq l_{1} \leq n_{2} \\ 0 \leq l_{2}}} a_{n_{1},n_{1}+n_{2}}
{n_{2} \choose l_{1}}{-n_{2} \choose l_{2}} 
(-c_{1})^{n_{1}+l_{1}}(-c_{2})^{n_{1}+n_{2}+l_{2}} 
$$
and after writing $(N,M)=(n_{1}+l_{1},n_{1}+n_{2}+l_{2})$,
$$ = \sum_{N,M\geq 0} c_{2}^{N}c_{1}^{M} \sum_{\substack{0 \leq l_{1} \leq n_{2} \leq M \\0 \leq l_{2} \leq N}} a_{n_{1},n_{1}+n_{2}} (-1)^{M-n_{1}+n_{2}}
{n_{2} \choose M-n_{1}} {-n_{2} \choose N-n_{1}-n_{2}} $$
\begin{equation} \label{eq:interm} = \sum_{N,M\geq 0} c_{2}^{N}c_{1}^{M} (-1)^{N+M} \sum_{\substack{0 \leq \leq n_{1} \leq N \\0 \leq n_{1}+n_{2} \leq M}} a_{n_{1},n_{1}+n_{2}} 
{n_{2} \choose M-n_{1}}
{N-n_{1}-1 \choose n_{2}-1} .
\end{equation}
On the other hand, for all integers $0 \leqslant a \leqslant b$, we have 
$$ {b \choose a} = \frac{b}{a} \bigg( \frac{b}{1}-1 \bigg) \ldots \bigg( \frac{b}{a-1}-1 \bigg) = \frac{b}{a} 
\sum_{r \geq 0}  b^{r} \frak{h}_{a}(\underbrace{1,\ldots,1}_{r})(-1)^{b-r}$$
where, for any positive integer $m$ and any word $w$, $\frak{h}_{m}(w)$ is defined by $\har_{m}(w)=m^{\weight(w)}\frak{h}_{m}(w)$. 
Whence (\ref{eq:interm}) equals 
\begin{multline*} \sum_{N,M\geq 0} c_{2}^{N}c_{1}^{M} (-1)^{N+M} \sum_{\substack{0 \leq l_{1} \leq n_{2} \leq M \\0 \leq l_{2} \leq N \\ n_{1}=N-l_{1}}} a_{n_{1},n_{1}+n_{2}}
\\
\sum_{R_{1}\geq 0} n_{2}^{R_{1}+1}\frac{(-1)^{n_{2}-R_{1}}}{M-n_{1}} 
\frak{h}_{M-n_{1}}(\underbrace{1,\ldots,1}_{R_{1}})
\sum_{R_{2} \geq 0} \frac{(-1)^{N-n_{1}-1-R_{2}} }{n_{2}-1} (N-n_{1}-1)^{r+1} \frak{h}_{n_{2}-1}(\underbrace{1,\ldots,1}_{R_{2}}) .
\end{multline*}
We now assume that $N=M=p^{\alpha}$ (a more general version of the computation in which $M$ and $N$ are any powers of $p$). In the domain of summation 
$\{(m_{1},\ldots,m_{R_{1}}) \in \mathbb{N}^{r}\text{ }|\text{ }0<m_{1}<\ldots<m_{R_{1}} \}$ of $\frak{h}_{M-n_{1}}(\underbrace{1,\ldots,1}_{R_{1}})$, we make the change of variable $m'_{i}=M-R_{1}$. We conclude by using the property of delocalization of cyclotomic multiple harmonic sums established in \cite{I-2}, Proposition-Definition 4.2.2.
\end{proof}

\begin{Definition} Let $M_{\har}^{\smallint}$ be the affine ind-scheme defined by the equations of Proposition \ref{harmonic duality DR RT}.
\end{Definition}

Let us consider all homographies of $\mathbb{P}^{1}$ which preserve $0$. For $a,c,d \in K$ with $K=\mathbb{C}$ or $K=\mathbb{C}_{p}$, let $\sigma_{a,c,d} : z \in \mathbb{P}^{1}(K)\mapsto \frac{az}{cz+d} \in \mathbb{P}^{1}(K)$.

\begin{Proposition} \label{proposition 5 30}
	The horizontality of the morphisms of the form $(\sigma_{a,c,d})_{\ast} : \pi_{1}^{\un,\dR}(\mathbb{P}^{1} - D) \rightarrow  \pi_{1}^{\un,\dR}(\mathbb{P}^{1} - D')$ with respect to $\nabla_{\KZ}$ gives an explicit family of relations between prime weighted multiple harmonic sums.
\end{Proposition}

\begin{proof} Let $\big((n_{i})_{d}; (z_{i})_{d+1} \big)$ be any index. For $z$ close to $0$, we consider the following equality of iterated integrals :
	\begin{equation} \label{eq:formula 0} \int_{0}^{\sigma_{a,c,d}(z)} \omega_{0}^{n_{d}-1}\omega_{z_{d}} \ldots  \omega_{0}^{n_{1}-1}\omega_{z_{1}} =
	\int_{0}^{z} (\sigma_{a,c,d})^{\ast} \big( \omega_{0}^{n_{d}-1}\omega_{z_{d}} \ldots  \omega_{0}^{n_{1}-1}\omega_{z_{1}}\big) .
	\end{equation}
	\noindent The left-hand side of (\ref{eq:formula 0}) is $\Li \big((n_{i})_{d};(z_{i})_{d+1}\big) \big(\frac{az}{cz+d} \big)$ ; let us write its power series expansion at $0$. For $z \in \mathbb{C}$ on a neighborhood of $0$, we have, for all $n \in \mathbb{N}$ :
	$\big(\frac{az}{cz+d} \big)^{n}
	= \big(\frac{az}{d} \big)^{n} \sum\limits_{l\geqslant  0} {-n_{d} \choose l} \big(\frac{c}{d}z\big)^{l}$. Thus, by the power series expansion of multiple polylogarithms (\ref{eq:Li series bis}),
	\begin{equation} \label{eq:first step} \Li \big( (n_{i})_{d};(z_{i})_{d} \big) \bigg(\frac{az}{cz+d} \bigg) = \sum_{0<m_{1}<\ldots<m_{d}<m} 
	\frac{ (\frac{z_{i_{2}}}{z_{i_{1}}})^{m_{1}} \ldots (\frac{z_{i_{d}}}{z_{i_{d-1}}})^{m_{d-1}} (\frac{a}{c z_{i_{d}}})^{m_{d}}}{m_{1}^{n_{1}} \ldots m_{d}^{n_{d}}} \big(\frac{c}{d}\big)^{m} {m-1 \choose m - m_{d}} .
	\end{equation}
	\noindent For all $m \in \mathbb{N}^{\ast}$, and $\tilde{m} \in \{1,\ldots,m-1\}$, we have 
	$\displaystyle {m - 1 \choose m - \tilde{m}} =  
	\frac{(m-1)(m-2) \ldots (m - (m-\tilde{m}))}{1 \times 2 \times \ldots \times (m-\tilde{m})}
	= \big(\frac{m}{1}-1\big)\big(\frac{m}{2}-1\big) \ldots \big(\frac{m}{m-\tilde{m}}-1\big)$ ; expanding this product gives : $\displaystyle {m-1 \choose m - \tilde{m}} = \sum\limits_{r \geqslant  0} m^{r} (-1)^{m-\tilde{m}-r} \frak{h}_{m-\tilde{m}} (\underset{r}{\underbrace{1,\ldots,1}})$  where the sum over $r$ is finite. By a change of variable, we have, for all $r \in \mathbb{N}$ : $\displaystyle \frak{h}_{m-\tilde{m}} (\underset{r}{\underbrace{1,\ldots,1}}) = \sum\limits_{\tilde{m}<j_{r}< \ldots <j_{1}<m} \frac{(-1)^{r}}{(m-j_{r}) \ldots (m-j_{1})}$. Now, we assume that $m = p^{\alpha}$ with $p$ a prime number and $\alpha \in \mathbb{N}^{\ast}$. By factorizing by $ \frac{1}{(-j_{1})\ldots (-j_{r})}$ and expanding into $p$-adic series factors of the form $\frac{1}{1-x}$ with $|x|_{p}<1$, we obtain 
	\begin{equation} \label{eq:second step} {p^{\alpha} - 1 \choose p^{\alpha} - \tilde{m}} = (-1)^{\tilde{m}-1} 
	\sum_{\substack{ r \geqslant  0 \\  l_{1},\ldots,l_{r} \geqslant  0}} (p^{\alpha})^{r + \sum\limits_{i=1}^{r} l_{i} }\sum_{\tilde{m}<j_{r}<\ldots <j_{1}<p^{\alpha}} \frac{1}{j_{r}^{1+l_{r}} \ldots j_{1}^{1+l_{1}}}
	\end{equation}
	\noindent By (\ref{eq:first step}) and (\ref{eq:second step}) we have :
	\begin{multline} \label{eq: 5 4 4} (p^{\alpha})^{n_{1}+\ldots+n_{d}}\Li \Big( (n_{i})_{d};(z_{i})_{d} \Big) \big(\frac{az}{cz+d} \big)[z^{p^{\alpha}}] \\
	= - \sum_{\substack{ r \geqslant  0 \\ l_{1},\ldots,l_{r} \geqslant  0}} 
	(p^{\alpha})^{\sum\limits_{j=1}^{d}n_{j} + r + \sum\limits_{i=1}^{r} l_{i}}\sum_{0<n_{1}<\ldots<n_{d}<j_{r}<\ldots<j_{1}<p^{\alpha}}
	\frac{ (\frac{z_{i_{2}}}{z_{i_{1}}})^{m_{1}} \ldots (\frac{z_{i_{d}}}{z_{i_{d-1}}})^{m_{d-1}} (-\frac{a}{c z_{i_{d}}})^{m_{d}}}{m_{1}^{n_{1}} \ldots m_{d}^{n_{d}} j_{r}^{l_{r}+1} \ldots j_{1}^{l_{1}+1}} \big(\frac{c}{d} \big)^{p^{\alpha}}
	\\ = -\sum_{\substack{ r \geqslant 0 \\  l_{1},\ldots,l_{r} \geqslant  0}} \har_{p^{\alpha}}
	\bigg( (n_{i})_{d};(1+l_{i})_{d}; (z_{i_{j}})_{d},\big(\frac{-a}{c}\big)_{r},
	\frac{-ad}{c^{2}} \bigg)
	\end{multline}
	\noindent The right-hand side of (\ref{eq:formula 0}) can be expressed in terms of multiple harmonic sums via (\ref{eq:Li series bis}) and the fact that, for all $y\in \mathbb{C}$, we have :
	$$ \sigma_{a,c,d}^{\ast}\bigg(\frac{dz}{z - y}\bigg) = \bigg( \frac{(a- yc)}{(a-y c)z + (-yd)} - \frac{c}{cz+d} \bigg) dz $$
\end{proof}

\subsubsection{In the framework $\Sigma$}

\begin{Proposition} \label{harmonic duality RT} We can prove in the framework $\Sigma$ that, for all words $w$ :
\begin{equation}
\label{eq: harmonic duality RT} 	 
\har_{\mathcal{P}^{\mathbb{N}}}( w(e_{0}+e_{1},-e_{1})) = - \sum_{\substack{d'\geqslant  1 \\ z = e_{0}^{t_{d'}-1}e_{1}\ldots e_{0}^{t_{1}-1}e_{1}}} (-1)^{\depth(z)} \har_{\mathcal{P}^{\mathbb{N}}}(z.w) .
\end{equation}
\end{Proposition}

This is a generalization of a result of Rosen called the ``asymptotic duality theorem'' \cite{Rosen}, which is equivalent to the $\alpha=1$ case of this result, although it is formulated differently.

\begin{proof} In Rosen's proof of the asymptotic duality theorem \cite{Rosen}, which we generalize from $\alpha=1$ to any $\alpha \in \mathbb{N}^{\ast}$, let us modify the last step by writing, for all $n_{d} \in \{1,\ldots,p^{\alpha}-1\}$, ${p^{\alpha} - 1 \choose m_{d} - 1} = {p^{\alpha} - 1 \choose p^{\alpha} - m_{d}}$ and applying the canonical expansion of binomial coefficients in terms of multiple harmonic sums to ${p^{\alpha} - 1 \choose p^{\alpha} - m_{d}}$ instead of ${p^{\alpha} - 1 \choose m_{d}-1}$. Instead of obtaining quantities of the form $\sum_{0<n<p^{\alpha}} \frak{h}_{m}(w') \frak{h}_{m}(w'')$ and linearizing them by the quasi-shuffle formula as in \cite{Rosen}, we obtain quantities of the form $\sum_{0<m<p^{\alpha}} \frak{h}_{m}(w')\frak{h}_{p^{\alpha}-m}(w'')$, and we make a change of variable $m_{i}\mapsto p^{\alpha}-m_{i}$ in the domain of summation of $\frak{h}_{p^{\alpha}-m}(w'')$, which gives an infinite sums of prime weighted multiple harmonic sums $\har_{p^{\alpha}}(w''')$.
\end{proof}

The asymptotic duality theorem in \cite{Rosen} is that for all indices $w$, we have, for all primes $p$,
\begin{equation} \label{eq: remarque 1}
\har_{p}\big( w(e_{0}+e_{1},-e_{1})\big) = \har_{p}\big( w + (w \ast (\frac{1}{1+y_{1}}) \big)
\end{equation}

It is obtained as a $p$-adic lift of a theorem of Hoffman \cite{Hoffman 2}, about multiple harmonic sums $\frak{h}_{p}(w) \mod p$, called the "duality theorem", which relies on the Newton series of multiple harmonic sums, in the sense of \S2.2.4. The proof in \cite{Hoffman 2} and \cite{Rosen} remains true if one replaces $\har_{p}$ by $\har_{p^{\alpha}}$ for all $\alpha \in \mathbb{N}^{\ast}$. Our alternative version (\ref{eq: harmonic duality RT}) does not use the quasi-shuffle equation. However, given the quasi-shuffle relation for $\har_{p^{\alpha}}$, equations (\ref{eq: remarque 1}) and the generalization of (\ref{eq: harmonic duality RT}) to any $\alpha$ are equivalent. Indeed, we have, for all $w = e_{0}^{n_{d}-1}e_{1}\ldots e_{0}^{n_{1}-1}e_{1}$,
\begin{equation} \label{eq: remarque 2}
\har_{p^{\alpha}}(w(e_{0}+e_{1},-e_{1}))
= \har_{p^{\alpha}}(e_{0}^{n_{d}}e_{1} - e_{0}^{n_{d}-1}e_{1})(e_{0}^{n_{d-1}-1}e_{1}\ldots e_{0}^{n_{1}-1}e_{1} \ast \frac{1}{1+e_{1}}) .
\end{equation}
\noindent Finally, for all words $w$ and all $n \in \mathbb{N}$, we have :
$$ \har_{p^{\alpha}}\big( e_{0}^{n}e_{1} ( w \ast \frac{1}{1+e_{1}}) \big)
= - \sum_{\substack{d\geqslant  1 \\ x_{d},\ldots,x_{1} \geqslant  1 \\ z=e_{0}^{x_{d}-1}e_{1}\ldots e_{0}^{n+x_{1}-1}e_{1}}} (-1)^{\depth(z)}\har_{p^{\alpha}}(z.w) . $$

\subsubsection{Comparison between the results of the three frameworks}

We observe that Proposition \ref{harmonic duality DR} (i), Proposition 4.3.4 (i) and Proposition 4.3.6, obtained respectively in the frameworks $\int_{1,0}$, $\int$ and $\Sigma$ are the same results. This is the part of the harmonic associator equations obtained by using $\pi_{1}^{\un}(\mathcal{M}_{0,4})$.

\section{Finite cyclotomic multiple zeta values and finite multiple polylogarithms}

We define ``finite'' analogues of adjoint cyclotomic multiple zeta values and study their algebraic properties. This gives a generalization of the notion of finite multiple zeta values of Kaneko and Zagier and an interpretation of that notion in terms of the crystalline pro-unipotent fundamental groupoid of $\mathbb{P}^{1} - \{0,\mu_{N},\infty\}$. 

\subsection{Review on finite multiple zeta values}

Let $\mathcal{P}$ be the set of prime numbers. The following ring is a $\mathbb{Q}$-algebra, it is denoted by $\mathcal{A}$ by Kaneko and Zagier and is called the ring of integers modulo infinitely large primes \cite{Kontsevich} :

\begin{equation}
\label{eq:integers mod infinitely large}
\mathbb{F}_{p\rightarrow\infty} = 
 \mathbb{Z}/_{p\rightarrow\infty} = \big( \prod_{p \in \mathcal{P}} \mathbb{Z}/p\mathbb{Z} \big) / \big( \bigoplus_{p \in \mathcal{P}} \mathbb{Z}/p\mathbb{Z} \big) .
\end{equation}

Kaneko and Zagier have defined finite multiple zeta values as the following numbers, where $d$ and $n_{i}$ $(1\leqslant i \leqslant d)$ are any positive integers (see also \cite{Hoffman 2} and \cite{Zhao} for earlier almost identical notions) :

\begin{equation} \zeta_{\mathcal{A}} \big((n_{i})_{d} \big)  =  \left(\sum_{0<m_{1}<\ldots<m_{d}<p} \frac{1}{m_{1}^{n_{1}} \ldots m_{d}^{n_{d}}}  \mod p \right)_{p \in \mathcal{P}} \in \mathcal{A} .
\end{equation}

They conjecture that the following formula defines an isomorphism between the $\mathbb{Q}$-algebra generated by fMZV's and the $\mathbb{Q}$-algebra generated by MZV's moded out by the ideal $(\zeta(2))$ : the explicit formula is striking
\begin{equation} \label{eq:Kaneko Zagier iso}
\zeta_{\mathcal{A}}\big( (n_{i})_{d} \big) \mapsto \sum_{d'=0}^{d} (-1)^{n_{d'+1}+\ldots+n_{d}}
\zeta(n_{1},\ldots,n_{d'})\zeta(n_{d},\ldots,n_{d'+1}) \mod \zeta(2) .
\end{equation}
According to the period conjectures on MZV's and $p$MZV's, it is equivalent to make the same conjecture with $p$MZV's instead of MZV's modulo $\zeta(2)$. With Definition \ref{def adjoint}, the $p$-adic analogue of the right-hand side of (\ref{eq:Kaneko Zagier iso}) is $\zeta_{p,\alpha}^{\Ad}\big((n_{i})_{d};0\big)$. Their complex analogues are sometimes called symmetric or symmetrized multiple zeta values \cite{NoteCRAS}, or finite real multiple zeta values.

\subsection{Finite cyclotomic multiple zeta values}

\subsubsection{Variants of the ring of integers modulo infinitely large primes}

We introduce some variants of the ring  $\mathbb{F}_{p\rightarrow \infty}$ (\ref{eq:integers mod infinitely large}) :

\begin{Definition} \label{def mod p infinitely large}(i) For each $\alpha \in \mathbb{N}^{\ast}$, let $\mathbb{F}_{p^{\alpha} \rightarrow \infty}=\big( \prod_{p \in \mathcal{P}} \mathbb{F}_{p^{\alpha}} \big) / \big( \bigoplus_{p \in \mathcal{P}} \mathbb{F}_{p^{\alpha}}\big)$.
	\newline (ii) Let $\overline{\mathbb{F}}_{p\rightarrow \infty} =  \big( \prod_{p \in \mathcal{P}} \overline{\mathbb{F}_{p}} \big) / \big( \bigoplus_{p \in \mathcal{P}} \overline{\mathbb{F}_{p}} \big)$.
	\newline (iii) Let the Frobenius of  $\mathbb{F}_{p^{\alpha}\rightarrow \infty}$ resp.  $\overline{\mathbb{F}}_{p\rightarrow \infty}$ be the automorphism defined as $\sigma : (x_{p})_{p} \mapsto (x_{p}^{p})_{p}$.
\end{Definition}

These rings are $\mathbb{Q}$-algebras. We have an inclusion
$\mathbb{F}_{p^{\alpha}\rightarrow \infty} \hookrightarrow  \overline{\mathbb{F}}_{p\rightarrow \infty}$, and, if $\alpha'$ divides $\alpha''$, we have an inclusion $\mathbb{F}_{p^{\alpha'}\rightarrow \infty} \hookrightarrow \mathbb{F}_{p^{\alpha''}\rightarrow \infty}$, and these inclusions are compatible.

\subsubsection{Finite cyclotomic multiple zeta values}

We now generalize the definition of finite multiple zeta values. Below, for any root of unity $\xi \in \overline{\mathbb{Q}}$ and for any prime $p$, we denote by $\overline{\xi}$ the image of $\xi$ in $\overline{\mathbb{F}_{p}}$.

\begin{Definition} \label{def cycl mzv} Let finite cyclotomic multiple zeta values (fMZV$\mu_{N}$'s) be the following numbers : for any positive integers $d$ and $n_{i}$ ($1\leqslant i \leqslant d$) and roots of unity $\xi_{i}$ ($1\leqslant i \leqslant d+1$),
	\begin{equation} \zeta_{f}
	\big( (n_{i})_{d};(\xi_{i})_{d+1} \big) =  \left(\sum_{0<m_{1}<\ldots<m_{d}<p} \frac{\big( \frac{\bar{\xi}_{2}}{\bar{\xi}_{1}} \big)^{m_{1}} \ldots \big(\frac{\bar{\xi}_{d+1}}{\bar{\xi}_{d}}\big)^{m_{d}} \big(\frac{1}{\bar{\xi}_{d+1}} \big)^{p}}{m_{1}^{n_{1}} \ldots m_{d}^{n_{d}}} \right)_{p \in \mathcal{P}} \in \overline{\mathbb{F}}_{p \rightarrow \infty} .
	\end{equation}
\end{Definition}

\noindent We note that, letting $N$ be the lcm of the orders of the $\xi_{i}$'s as roots of unity, then, for $p$ large enough, $p$ does not divide $N$, and the crystalline realization of $\pi_{1}^{\un}(\mathbb{P}^{1} \setminus \{0,\mu_{N},\infty\})$ is defined.
\newline\indent By extrapolating on Kaneko-Zagier's conjecture (\S6.1), we can speculate an equivalence between the properties of the numbers  $\zeta_{f}(n_{i})_{d};(\xi_{i})_{d+1}$ and the numbers 
\begin{equation} \label{eq: analogue of finite} \zeta_{p,\alpha}^{\Ad} \big( (n_{i})_{d};0;(\xi_{i})_{d}\big) .
\end{equation}

We take again the notations of Definition \ref{def mod p infinitely large}. For any $z \in D \setminus \{0,\infty\}$, we denote by $\bar{z}$ the reduction of $z$ modulo $p$ large enough.

\begin{Definition} Let finite multiple polylogarithms be the numbers :
$$ \Li_{f}
\big( (n_{i})_{d};(z_{i})_{d+1} \big) =  \left(\sum_{0<m_{1}<\ldots<m_{d}<p}
\frac{\big( \frac{\bar{z}_{2}}{\bar{z}_{1}} \big)^{m_{1}} \ldots \big(\frac{\bar{z}_{d+1}}{\bar{z}_{d}} \big)^{m_{d}} \big(\frac{1}{\bar{z}_{d+1}}\big)^{p} }{m_{1}^{n_{1}} \ldots m_{d}^{n_{d}}} \mod p \right)_{p \in \mathcal{P}} \in \overline{\mathbb{F}}_{p\rightarrow \infty} . $$
\end{Definition}

\subsubsection{Equations satisfied by finite cyclotomic multiple zeta values}

We have
$$ \overline{\mathbb{F}}_{p\rightarrow \infty} = \bigg\{ (x_{p})_{p} \in \prod_{p\in \mathcal{P}}\mathbb{C}_{p} \text{ }|\text{ } v_{p}(x_{p}) \geqslant  0 \text{ for p large} \bigg\} \text{ }\bigg/\text{ }\bigg\{ (x_{p})_{p} \in \prod_{p\in \mathcal{P}} \mathbb{C}_{p} \text{ }|\text{ } v_{p}(x_{p}) \geqslant  1 \text{ for p large} \bigg\} . $$
\indent Finite cyclotomic multiple zeta values are expressed as reductions of cyclotomic multiple harmonic values, up to the Frobenius of Definition \ref{def mod p infinitely large} (iii) : for any harmonic word $w$, we have  
\begin{equation} \label{eq:red mod cycl har} \sigma^{\alpha-1}\zeta_{f}\big(w\big) =
(p^{-\weight(w)}\har_{p^{\alpha}}(w) \mod p)_{p,\alpha} . 
\end{equation}
Indeed, with the notations of (\ref{eq:mult har sums}), the subsum of $\har_{p^{\alpha}}\big( (n_{i})_{d};(\xi_{i})_{d+1}\big)$ on the subdomain $(m_{1},\ldots,m_{d}) \in p^{\alpha-1}\mathbb{N}^{\ast} \times \ldots \times p^{\alpha-1}\mathbb{N}^{\ast}$ is equal, by the change of variable $m_{i}=p^{\alpha-1}m'_{i}$, to $\har_{p} \big( (n_{i})_{d};(\xi^{p^{\alpha}-1}_{i})_{d+1}) \big)$, which is of valuation $\geqslant  0$, and the subsum on the complementary domain has $p$-adic valuation $\geqslant 1$.
\newline\indent On the other hand, the numbers (\ref{eq: analogue of finite}) are reductions of the $\Lambda$-adjoint $p$-adic cyclotomic multiple zeta values :
\begin{equation} \label{eq:modulo Lambda}
\zeta_{p,\alpha}^{\Lambda \Ad}\big( (n_{i})_{d};0;(\xi_{i})_{d+1}\big)=  \zeta_{p,\alpha}^{\Lambda \Ad} \big((n_{i})_{d};(\xi_{i})_{d+1} \big) \mod \Lambda^{\weight(w)+1} .
\end{equation}

In the $N=1$ case, Akagi-Hirose-Yasuda have proved an integrality property of $p$MZV's and deduced a finite variant of equation (\ref{eq:formula for n=1}) \cite{AHY} :

\begin{equation} \label{eq: AHY} \zeta_{\mathbb{Z}_{/p\rightarrow \infty}}\big((n_{i})_{d}\big) = \big( p^{-\sum_{i=1}^{d}n_{i}}\zeta_{p,1}^{\Ad}\big((n_{i})_{d};0\big) \mod p \big)_{p\text{ prime}} \in \mathcal{A} .
\end{equation}

By the relations of iteration of the Frobenius (\cite{I-3}, equations (1.11), (1.12), (1.13) and Proposition 1.5.2), equation (\ref{eq: AHY}) remains true if $\zeta_{p,1}^{\Ad}$ is replaced by $\zeta_{p,\alpha}^{\Ad}$. Moreover, this can be generalized to any $N$ by using the integrality property of $p$-adic iterated integrals proved in \cite{Chatzis}.
\newline By equations (\ref{eq:red mod cycl har}), (\ref{eq:modulo Lambda})
and the cyclotomic generalization of (\ref{eq: AHY}), we deduce immediately from \S3 and \S4 some equations satisfied by fMZV$\mu_{N}$'s and their analogues (\ref{eq: analogue of finite}). This leads to the following definition :

\begin{Definition} (i) Let $\DS_{f}$ the scheme of finite double shuffle equations, be the pro-affine scheme defined the term of lowest weight in the equations defining $\DS_{\har}$ of \S3.3.
	\newline (ii) Let $\GRT_{f}$, the scheme of finite associator equations, be the pro-affine scheme defined by the term of lowest weight in the equations defining $\M_{\har}$ of \S4.3.
\end{Definition}

Some of these equations already appeared in the literature. In the $N=1$ case, the finite double shuffle equations appear in \cite{Kaneko}, \cite{Kaneko Zagier} and in \cite{NoteCRAS}. To our knowledge, the method of proof in \cite{Kaneko}, \cite{Kaneko Zagier} uses the framework $\int$.
The one dimensional part of the finite associator equations appears in \cite{Hoffman 2}, where it is called the "duality theorem".
The finite version of the reversal equation (\ref{paragraph reflexion}) in the $N=1$ case is the formula 
$\zeta_{\mathcal{A}}(n_{1},\ldots,n_{d}) = (-1)^{n_{1}+\ldots+n_{d}}\zeta_{\mathcal{A}}(n_{1},\ldots,n_{d}) \mod p$, which appears in \cite{Zhao}, lemma 3.3 and \cite{Hoffman 2}, theorem 4.5.

\subsection{A generalization : finite analogues of adjoint cyclotomic multiple zeta values}

We now associate a ``finite'' analogue to all adjoint cyclotomic multiple zeta values $\zeta_{p,\alpha}^{\Ad}\big((n_{i})_{d};b;(\xi_{i})_{d+1}\big)$, not only the case $b=0$. We are going to use the overconvergent cyclotomic multiple harmonic values (Definition \ref{variant harmonic}) :

\begin{Definition} Let the adjoint finite cyclotomic multiple zeta values (Ad$f$MZV$\mu_{N}$'s) be the numbers
	$$ \zeta^{\Ad}_{f,\alpha} \big( (n_{i})_{d};l;(\xi_{i})_{d+1}  \big) = 
	(p^{-(n_{1}+\ldots+n_{d}+l)}\har_{p,\alpha}^{\dagger}( (n_{i})_{d};l;(\xi_{i})_{d+1})  \mod p)_{p} \in \overline{\mathbb{F}}_{p\rightarrow \infty} .	$$
\end{Definition}

By the integrality of $p$-adic iterated integrals on punctured projective lines proved in \cite{Chatzis}, and by the expression of 
$\har_{p,\alpha}^{\dagger}\big((n_{i})_{d};l;(\xi_{i})_{d+1}\big)$ as an infinite sums of Ad$p$MZV$\mu_{N}$'s (Proposition \ref{prop formula dagger} (i)), the numbers $\zeta^{\Ad}_{f,\alpha} \big( (n_{i})_{d};l;(\xi_{i})_{d+1}  \big)$ are well defined and satisfy a generalization of equation (\ref{eq: AHY}) :
$$ \zeta_{f,\alpha} \big( (n_{i})_{d};l;(\xi_{i})_{d+1}  \big) = \big( \zeta^{\Ad}_{p,\alpha} \big((n_{i})_{d};l;(\xi_{i})_{d+1}\big) \mod p \big) \in \overline{\mathbb{F}}_{p\rightarrow\infty} .
$$
We also have a generalization of equation (\ref{eq:modulo Lambda}), which refers to the overconvergent variant of $\Lambda$Ad$p$MZV$\mu_{N}$'s (Definition \ref{def over adjoint}) :
$$ \zeta^{\Ad}_{p,\alpha} \big((n_{i})_{d};l;(\xi_{i})_{d+1}\big) = \zeta_{p,\alpha}^{\Lambda \Ad \dagger} \big((n_{i})_{d};(\xi_{i})_{d+1}\big) \mod \Lambda^{\weight(w)+1} $$

In the $N=1$ case, the numbers $\zeta_{f,\alpha} \big( (n_{i})_{d};l\big)$ are in the $\mathbb{Q}$-vector space generated by finite multiple zeta values. This follows from equation (\ref{eq: AHY}) and the fact that, by \cite{Yasuda}, that the numbers $\zeta^{\Ad}_{p,\alpha}\big((n_{i})_{d};0\big)$ generate the space of $p$MZV's. We guess that this remains true for any $N$, up to the Frobenius $\sigma$ of $\overline{\mathbb{F}}_{p\rightarrow \infty}$ (Definition \ref{def mod p infinitely large}). We generalize again the conjecture of Kaneko and Zagier :

\begin{Conjecture} The numbers 
	$\zeta^{\Ad}_{f,\alpha} \big( (n_{i})_{d};l;(\xi_{i})_{d+1}\big)$ and $\zeta_{p,\alpha}^{\Ad} \big( (n_{i})_{d};l;(\xi_{i})_{d+1}\big)$ satisfy the same algebraic properties.
\end{Conjecture}

In a summary, starting with MHV$\mu_{N}$'s, one can on the one hand consider the associated graded for the weight filtration which gives Ad$p$MZV$\mu_{N}$'s, and on the other hand consider the associated graded for the filtration defined by the uniform topology on $\underset{p\in\mathcal{P}}{\prod} \mathbb{C}_{p}$ which arises from the $p$-adic topologies (for $p$ infinitely large), which gives Ad$f$MZV$\mu_{N}$'s ; the conjecture says that these two grading operations give equivalent results. Thus, this conjecture states a certain adelic case of equality in the inequality between the slopes of the Frobenius and the Hodge filtration on the log-crystalline cohomology which represents $\pi_{1}^{\un,\crys}(\mathbb{P}^{1} \setminus\{0,\mu_{N},\infty\},\vec{1}_{1},\vec{1}_{0})$. As a conclusion, we have interpreted Kaneko-Zagier's notion of finite multiple zeta values in crystalline terms. We will study it further in a subsequent paper.

\end{document}